\documentclass[a4paper,leqno,11pt]{amsart}
\usepackage{latexsym,amsmath,amsfonts,amscd,amssymb,epic,times}
\numberwithin{equation}{section}

\newtheorem {theorem}[equation]         {Theorem}

\newtheorem {lemma}[equation]           {Lemma}
\newtheorem {proposition}[equation]     {Proposition}

\newtheorem {corollary}[equation]       {Corollary}

\newtheorem {definition}[equation]      {Definition}
\newtheorem {remark}[equation]          {Remark}

\newcommand{\n}{\noindent}
\makeatletter
\newcommand{\bfsubsection}
{\@startsection{subsection}{3}{0pt}{\baselineskip}
 {-\fontdimen2\font}{\normalfont\bfseries}}
\newcommand{\bfsubsubsection}
{\@startsection{subsubsection}{3}{0pt}{\baselineskip}
 {-\fontdimen2\font}{\normalfont\bfseries}}
\makeatother

\newfont{\bfc}{cmbsy10 scaled 1200}  
\newfont{\dr}{msbm10 scaled \magstep1}  
\newfont{\sdr}{msbm8}  
\newfont{\gl}{eufm10 scaled \magstep1}  

\newcommand{\Met}{\mbox{$Met$}}
\newcommand{\IFF}{if and only if}

\DeclareFontFamily{OT1}{rsfs}{}
\DeclareFontShape{OT1}{rsfs}{n}{it}{<->rsfs10}{}
\DeclareMathAlphabet{\curly}{OT1}{rsfs}{n}{it}

 \newcommand{\CC}{{\mathbb C}}

 \newcommand{\PP}{{\mathbb P}}
 \newcommand{\QQ}{{\mathbb Q}}
 \newcommand{\RR}{{\mathbb R}}

 \newcommand{\ZZ}{{\mathbb Z}}

 \newcommand{\bmE}{\mbox{\boldmath$E$}}

 \newcommand{\bmR}{\mbox{\boldmath$R$}}

 \newcommand{\bmV}{\mbox{\boldmath$V$}}

 \newcommand{\bmf}{\mbox{\boldmath$f$}}
 \newcommand{\bmg}{\mbox{\boldmath$g$}}

 \newcommand{\bvarepsilon}{{\boldsymbol{\varepsilon}}}
 \newcommand{\bphi}{{\boldsymbol{\phi}}}
 
 \newcommand{\bpsi}{{\boldsymbol{\psi}}}
 
 \newcommand{\bvarphi}{{\boldsymbol{\varphi}}}
 
 \newcommand{\bsigma}{{\boldsymbol{\sigma}}}

 \newcommand{\btau}{{\boldsymbol{\tau}}}

 \newcommand{\cE}{{\mathcal E}}
 \newcommand{\cF}{{\mathcal F}}
 
 \newcommand{\cH}{{\mathcal H}}
 
 \newcommand{\cK}{{\mathcal K}}

 \newcommand{\cO}{{\mathcal O}}

 \newcommand{\cS}{{\mathcal S}}

 \newcommand{\cW}{{\mathcal W}}

 \newcommand{\bcC}{{\boldsymbol{\mathcal{C}}}}
 
 \newcommand{\bcE}{{\boldsymbol{\mathcal{E}}}}
 \newcommand{\bcF}{{\boldsymbol{\mathcal{F}}}}

 \newcommand{\bcR}{{\boldsymbol{\mathcal{R}}}}
 \newcommand{\bcS}{{\boldsymbol{\mathcal{S}}}}

 \newcommand{\glb}{\hbox{\gl b}} 
 \newcommand{\gld}{\hbox{\gl d}} 
 \newcommand{\glg}{\hbox{\gl g}} 
 \newcommand{\glh}{\hbox{\gl h}} 
 \newcommand{\glj}{\hbox{\gl j}} 
 \newcommand{\glk}{\hbox{\gl k}} 
 \newcommand{\gll}{\hbox{\gl l}} 
 \newcommand{\glp}{\hbox{\gl p}} 
 \newcommand{\glr}{\hbox{\gl r}} 
 \newcommand{\glt}{\hbox{\gl t}} 
 \newcommand{\glu}{\hbox{\gl u}} 
 \newcommand{\glz}{\hbox{\gl z}} 
 \newcommand{\glX}{{\hbox{\gl X}}}   

\newcommand{\CCC}{{\curly C}}
\newcommand{\DDD}{{\curly D}}

\newcommand{\GGG}{{\curly G}}

\newcommand{\NNN}{{\curly N}}

\newcommand{\RRR}{{\curly R}}

\newcommand{\qu}{/\kern-.7ex/}
\newcommand{\exh}{\to\kern-1.8ex\to}

 \newcommand{\ra}{\rightarrow}
 \newcommand{\lra}{\longrightarrow}
 \newcommand{\hra}{\hookrightarrow}

 \newcommand{\wt}{\widetilde} 

 \newcommand{\kahler}{K\"{a}hler}

\newcommand{\imag}{\mathop{{\fam0 \sqrt{-1}}}\nolimits}

\newcommand{\ad}{\operatorname{ad}}

\newcommand{\dbar}{\bar{\partial}}             
\newcommand{\Rep}{\operatorname{Rep}}

\newcommand{\Aut}{\operatorname{Aut}}
\newcommand{\End}{\operatorname{End}}

\newcommand{\Hom}{\operatorname{Hom}}
  
\newcommand{\id}{\operatorname{id}}

\renewcommand{\Im}{\operatorname{Im}}

\newcommand{\rk}{\operatorname{rk}}
\newcommand{\SL}{\operatorname{SL}}

\newcommand{\SU}{\operatorname{SU}}

\newcommand{\tr}{\operatorname{tr}}

\newcommand{\Vol}{\operatorname{Vol}}

\renewcommand{\exp}{\operatorname{exp}}

\newcommand{\hts}{{\rm ht}_{\Sigma }}

 \newcommand{\HKC}{Hitchin--Kobayashi correspondence}
 \newcommand{\kah}{\kahler}

 \newcommand{\holfil}{{0\hra \cF_0\hra \cF_1\hra\cdots\hra  \cF_{m}=\cF}}
 \newcommand{\holsfil}{{0\hra \cF_0'\hra \cF_1'\hra\cdots\hra\cF'_{m}=\cF'}}

\title[Dimensional reduction and quiver bundles]{Dimensional reduction
and quiver bundles} 
\author{Luis \'Alvarez--C\'onsul}  
\address{Department of Mathematics, University of Illinois at
Urbana-Champaign, Urbana, IL 61801, USA}
\email{lalvarez@math.uiuc.edu}
\author{Oscar Garc\'{\i}a--Prada}
\address{Departamento de Matem{\'a}ticas \\
Universidad Aut{\'o}noma de Madrid \\
28049 Madrid, SPAIN}
\email{oscar.garcia-prada@uam.es}

\begin{document}
\begin{abstract}
The so-called Hitchin--Kobayashi correspondence, proved by Donaldson,
Uhlenbeck and Yau, establishes that an indecomposable holomorphic vector 
bundle over a
compact \kah\ manifold admits a Hermitian--Einstein metric if and only
if the bundle satisfies the Mumford--Takemoto stability condition. In
this paper we consider a  variant  of this correspondence for
$G$-equivariant vector bundles on the product of a compact \kah\
manifold $X$ by a flag manifold $G/P$, where $G$ is a complex
semisimple Lie group and $P$ is a parabolic subgroup. The modification
that we consider is determined by a filtration of the vector bundle which
is naturally defined by the equivariance of the bundle.
The study of invariant solutions to the modified
Hermitian--Einstein equation over $X\times G/P$ leads, via dimensional
reduction techniques, to gauge-theoretic equations on $X$. These are
equations for hermitian metrics on a set of holomorphic bundles on $X$
linked by morphisms, defining what we call a {\em quiver bundle} for  
a  quiver with relations whose 
structure is entirely determined by the parabolic subgroup $P$. 
Similarly, the corresponding stability condition for the invariant filtration
over $X\times G/P$ gives rise to a stability condition for the quiver
bundle on $X$, and hence to a Hitchin--Kobayashi correspondence.
In the simplest case, when
the flag manifold is the complex projective line, one recovers the
theory of vortices, stable triples and stable chains,
as studied by  Bradlow, the  authors, and others. 
\end{abstract}

\maketitle

\section*{Introduction}
Let $M$ be a compact K\"ahler manifold and let $\cF$ be a
holomorphic vector bundle over $M$. It is well-known that 
a natural differential equation to consider for a Hermitian
metric $h$ on $\cF$ is the {\em Hermitian--Einstein} equation, also
referred sometimes as the {\em Hermitian--Yang--Mills} equation.
This says that $F_h$, the curvature of the Chern connection of
$h$ must satisfy
$$
\imag\Lambda F_h=\lambda I,
$$
were $\Lambda$ is the contraction with the \kahler\ form of $M$,
$\lambda$ is a real number determined by the topology and $I$ is the
identity endomorphism of $\cF$. A theorem of Donaldson, Uhlenbeck and
Yau \cite{D1,D2,UY}, also known as the {\em Hitchin--Kobayashi} correspondence
establishes that the existence of such metric is equivalent to the
stability of $\cF$ in the sense of Mumford--Takemoto.

Suppose now that a compact Lie group $K$ acts on $M$ by isometries so
that $M/K$ is a smooth \kahler\ manifold and the action on $M$ can be
lifted to an action on $\cF$. One can then apply {\em dimensional reduction}
techniques to study $K$-invariant solutions to the Hermitian--Einstein
equation on $\cF$ and the corresponding stability condition to obtain
a theory expressed entirely in terms of the orbit space $M/K$. Many 
important equations in gauge theory arise in this way (cf. e.g. \cite{AH,FM,H,W}).
 In this paper we carry out this programme for the manifold $M=X\times
G/P$, where $X$ is a compact \kahler\ manifold, $G$ is a connected
simply connected semisimple complex Lie group and $P\subset G$ is a
parabolic subgroup, i.e. $G/P$ is a flag manifold. The group $G$ (and
hence, its maximal compact subgroup $K\subset G$) 
act trivially  on $X$ and in  the standard way on $G/P$. The \kahler\
structure on $X$ together with a $K$-invariant \kahler\  structure on 
$G/P$ define a product \kahler\ structure on $X\times G/P$. In
\cite{AG1} we studied the case in which $G/P=\PP^1$, the complex
projective line, which is obtained as the quotient of 
$G=\SL(2,\CC)$ by the subgroup of lower
triangular matrices 
$$
P=\begin{pmatrix} * & 0 \\ * & * \end{pmatrix},
$$
generalizing previous work by \cite{G1,G2,BG}.

Already in the study of the dimensional reduction on $X\times\PP^1$ \cite{AG1}
one realizes that the Hermitian--Einstein is not quite the appropriate
equation to consider. It turns out that every $\SL(2,\CC)$-equivariant 
holomorphic vector bundle on $X\times \PP^1$ admits an equivariant
holomorphic filtration
$$
\bcF:\, 0\hra \cF_0\hra \cF_1\hra\cdots\hra  \cF_{m}=\cF,
$$
which, in turn, is in one-to-one correspondence with a {\em chain}
$$
\cE_m\stackrel{\phi_m}{\lra }\cE_{m-1}\stackrel{\phi_{m-1}}{\lra} 
        \cdots\stackrel{\phi_1}{\lra}\cE_0,
$$
consisting of holomorphic vector bundles $\cE_i$ on $X$ and
morphisms $\phi_i:\cE_i\ra\cE_{i-1}$. If one considers
$\SU(2)$-invariant solutions to the Hermitian--Einstein 
equation on $\cF$ one obtains certain equations of vortex-type
for Hermitian metrics on $\cE_i$. These involve, of course the 
homomorphims $\phi_i$, which in this context are referred as 
{\em Higgs fields}. 
The key point is that these equations naturally
have as many real parameters as morphisms in the chain. By
weighting  the \kahler\ structure on $X\times \PP^1$ one can accommodate
one parameter but not all, unless the  chains are just one step
 ---the so-called {\em triples}--- \cite{G2,BG}. In the general case 
the filtration on $\cF$ has $m$ steps, and a natural  equation to consider
for a metric $h$ on $\cF$ is
$$
\imag\Lambda F_h=\begin{pmatrix}
                \tau_0 I_0      &               &        &              \\ 
                                & \tau_1 I_1    &        &              \\ 
                                &               & \ddots &              \\ 
                                &               &        & \tau_m I_m
\end{pmatrix},
$$
where the RHS is a diagonal
matrix, with constants $\tau_0,\tau_1,\ldots ,\tau_m\in\RR$, written in
blocks corresponding to the splitting which a hermitian metric $h$ defines
in the filtration $\bcF$. 
This equation reduces of course to the Hermitian--Einstein equation
when $\tau_0=\tau_1=\cdots=\tau_m=\lambda$.
In \cite{AG1} we proved a Hitchin--Kobayashi correspondence 
relating this equation to a stability criterion for the filtration 
$\bcF$ that  depends on the parameters, of which  only $m$ 
are actually independent.

Coming back to the general situation that we will deal with  in this paper,
the first thing we do  to study dimensional  reduction on 
$X\times G/P$ is to analyse the structure 
$G$-equivariant holomorphic vector bundles on $X\times G/P$. As in
the case of $X\times \PP^1$, the basic tool we use is the study of
the holomorphic representations of $P$. This is simply because $G$-equivariant
bundles on $X\times G/P$ are in one-to-one correspondence with 
$P$-equivariant bundles on $X$, where the action of $P$ on $X$ is
trivial. Each fibre of such a bundle is hence  a representation of $P$.
This is carried out in Section  \ref{sec:homvb-quiver}.
A key fact in our study is the existence of a quiver with relations
$(Q,\cK)$ naturally associated to the subgroup 
$P$. A {\em quiver} is a pair $Q=(Q_0,Q_1)$ formed by two sets, where $Q_0$
is the set of {\em vertices} and $Q_1$ is the set of {\em arrows},
together with two maps $t,h:Q_1\ra Q_0$, which to any arrow
$a:\lambda\ra\mu$ associate the {\em tail} $ta=\lambda$ of the arrow
and the {\em head} $ha=\mu$ of the arrow.
A {\em relation} of the quiver is a formal 
complex linear combination of {\em paths} of the quiver.
The set of vertices in this case coincides with the set of irreducible 
representations of $P$. The description of the arrows and relations involves
studying certain isotopical decompositions related to the nilradical
of the Lie algebra of $P$. This construction was initially studied by
Bondal and Kapranov in \cite{BK} and later by Hille
(cf. e.g. \cite{Hl2}). 
The quiver we obtain is the same as in \cite{BK}, but different from
\cite{Hl2}, while the relations differ from
those of the previous authors. This is due to the possibility of
defining {\em different} quivers with relations whose categories of
representations are equivalent (the set of vertices is always the same). 
The quiver and relations that we obtain seem to be arising more
naturally from the point of view of dimensional reduction (that is,
when studying the relation between the moment maps for actions on
the spaces of unitary connections induced by actions on homogeneous
bundles, and the moment maps for actions on quiver representations). 
In Section \ref{sec:homvb-quiver}, we also provide examples of quivers
and relations associated  to $P$, for some parabolic subgroups. 

A representation of a quiver with relation $(Q,\cK)$ consists of a
collection of complex vector spaces $V_\lambda$ indexed by the
vertices $\lambda\in Q_0$, and a collection of linear maps
$\phi_a:V_{ta}\ra V_{ha}$ indexed by the arrows $a\in Q_1$, satisfying
the relations $\cK$. The crucial fact is that holomorphic
representations of $P$ are in one-to-one correspondence with
representations of $(Q,\cK)$ which, in turn, are in one-to-one
correspondence with holomorphic homogeneous vector bundles on $G/P$. In Section 
\ref{sec:equivb-quivers} we extend this result to $G$-equivariant
bundles on $X\times G/P$. To do this, we first introduce the notion 
of {\em quiver bundle} (term due to Alastair King), more precisely $(Q,\cK)$-bundle on $X$. This
is a  collection of holomorphic vector bundles $\cE_\lambda$ on $X$
indexed by $\lambda \in Q_0$ and a collection of homomorphisms 
$\phi_a:\cE_{ta}\ra \cE_{ha}$ indexed by $a\in Q_1$, that satisfy the
relations $\cK$. One has that there is an equivalence of categories 
between $G$-equivariant bundles on $X\times G/P$ and $(Q,\cK)$-bundles
on $X$. By standard techniques in representation theory of complex Lie
groups on Fr\'echet spaces of sections of equivariant coherent sheaves
(with a holomorphic action of the complex Lie group), this equivalence
can be extended to coherent sheaves. 

It turns out that the parabolicity of $P$ implies that the quiver has
no oriented cycles. An important consequence of this is  that a
$G$-equivariant holomorphic vector bundle admits a natural
$G$-equivariant filtration by holomorphic subbundles. This suggests that the natural equation 
for a Hermitian metric to consider is again the deformed
Hermitian--Einstein equation, like in the case of $X\times \PP^1$.
To study this equation and its dimensional reduction we first study 
in Section \ref{sec:hol-str} the correspondence established in Section 
\ref{sec:equivb-quivers} by means of Dolbeault operators. We reinterpret 
in these terms many 
of the ingredients studied  in the previous section. In particular, 
we see that the relations of the quiver correspond to the integrability
of the Dolbeault operator on the corresponding homogeneous vector bundle
on $G/P$. We go on to consider the dimensional reduction of a $K$-invariant
solution of the filtered Hermitian--Einstein equation on  a
$G$-equivariant vector bundle $\cF$ over $X\times G/P$.
We obtain that such solutions are in correspondence with a collection
of metrics satisfying a set of vortex-type equations on the bundles
$\cE_\lambda$ in the $(Q,\cK)$-bundle defined by $\cF$. We also show that
the stability of $\cF$ as an equivariant filtration is equivalent 
to an appropriate stability condition for the $(Q,\cK)$-bundle on $X$.
As a corollary we obtain a Hitchin--Kobayashi correspondence for 
$(Q,\cK)$-bundles. We finish this section by observing that the quiver 
vortex equations that we obtain make actually perfect sense for arbitrary
quivers with relations (i.e. not necessarily related to a parabolic group). 
One can indeed proof a Hitchin--Kobayasi correspondence in this more
general situation --- this is a subject of the paper \cite{AG2}.

In the last section, we use the examples provided in Section
\ref{sec:homvb-quiver}, to illustrate the general theory developed
throughtout this paper. We compute explicitly the dimensional
reduction of the gauge equations and stability conditions. In
particular, when $G/P$ is the projective line  we recover the results
of \cite{AG1,BG} from our general theory. 

A dimensional reduction problem closely related to the one treated  in this 
paper has  been studied  by Steven Bradlow, Jim Glazebrook and Franz Kamber \cite{BGK}.
It would be very interesting to combine the two points of view to
pursue further research in the subject.

\bfsubsection*{Acknowledgements} 

We would like to thank Alastair King for very useful discussions.
The research of L.A. was partially supported by the Comunidad
Aut\'onoma de Madrid (Spain) under a FPI Grant. 
The  authors are members of VBAC (Vector Bundles on Algebraic Curves),
which is partially supported by EAGER (EC FP5 Contract no. HPRN-CT-2000-00099)
and by EDGE (EC FP5 Contract no. HPRN-CT-2000-00101).

\section{Homogeneous bundles and quiver representations}
\label{sec:homvb-quiver}

One goal of this paper is to investigate the dimensional
reduction of natural gauge equations and stability conditions for
equivariant vector bundles on the product of 
a compact \kah\ manifold $X$ by a flag variety $G/P$, where $G$ is a
connected semisimple complex Lie group, and $P\subset G$ is a
parabolic subgroup. As we shall see, the process of dimensional
reduction leads to new objects defined on $X$, called quiver bundles,
which we shall define in \S \ref{sec:equivb-quivers}, 
where we shall also prove that there is a one-to-one relation between
equivariant vector bundles on $X\times G/P$ and quiver bundles on $X$. 
Preparing and leading up to this correspondence, 
this section is concerned with the case where $X$ is
a single point. This section is organised as follows.  After introducing
notation and some preliminaries in \S\S \ref{sec:preliminaries} and
\ref{sec:A-B-psi},  in \S \ref{sec:Q,K,P} we define the quiver with relations
associated to any complex Lie group. In \S \ref{sec:P-modules-quiver-rep} 
we prove the equivalence between the category of holomorphic 
homogeneous vector bundles on the flag variety $G/P$ and the 
category of representations of the quiver with relations associated 
to $P$. Then, in \S \ref{sec:flags} we prove that the quiver 
associated to a parabolic subgroup $P\subset G$ has no oriented
cycles, and as a result, the holomorphic homogeneous vector bundles on
$G/P$ admit certain filtrations by homogeneous subbundles. 

\bfsubsection{Preliminaries and notation}
\label{sec:preliminaries}

Througout this paper, $G$ is a connected simply connected complex
semisimple Lie group, $P$ is a parabolic subgroup of $G$, and $X$ is a
compact \kah\ manifold. Given a complex manifold $M$ with a
holomorphic action of $G$, a $G$-equivariant holomophic vector bundle
on $M$ is a holomophic vector bundle $\pi:\cF\ra M$ on $M$, together
with  a holomorphic $G$-action $\rho:G\times\cF\ra\cF$ on its total
space $\cF$ which commutes with $\pi$, and such that for all $(g,p)\in
G\times M$, the map $\rho_{g,p}:\cF_p\ra\cF_{g\cdot p}$, from
$\cF_p=\pi^{-1}(p)$ into $\cF_{g\cdot p}=\pi^{-1}(g\cdot p)$, induced
by this action, is an isomorphism of vector spaces. In this paper we
are concerned with the $G$-equivariant holomorphic vector bundles on
the $G$-manifold $M=X\times G/P$. The group $G$ acts on $X\times G/P$
by the trivial action on $X$ and the (left) standard action on the
flag variety $G/P$. 

\bfsubsubsection{Induction and reduction}\label{subsub:ind-red}

It is a standard result that there is a one-to-one
correspondence between the isomorphism classes of $G$-equivariant
holomophic vector bundles on $X\times G/P$ and the isomorphism classes
of $P$-equivariant holomorphic vector bundles on $X$, where $P$ acts 
trivially on $X$. The correspondence is defined as follows. A
$G$-equivariant holomophic vector bundle $\cF$ on $X\times G/P$
defines a $P$-equivariant holomorphic vector bundle $i^*\cF$ on $X$
by restriction to the slice $i:X\cong X\times P/P\hra X\times G/P$. 
Conversely, a $P$-equivariant holomorphic vector bundle $\cE$ on 
$X$ defines by induction a $G$-equivariant holomophic vector bundle 
$G\times_P\cE$ on $X\times G/P$. This holomorphic vector bundle is by
definition the quotient of $G\times \cE$ by the action of $P$, defined
by $p\cdot (g,e)=(gp^{-1},p\cdot e)$ for $p\in P$ and $(g,e)\in
G\times\cE$. The bundle projection is
$$
G\times_P\cE\ra X\times G/P,  \quad [g,e]\mapsto (x,gP)
$$ 
for $e\in E_x,x\in X$. The vector bundle $G\times_P E\ra X\times G/P$
has an equivariant action of $G$ given by 
$$
G\times (G\times_P\cE)\ra G\times_P\cE, \quad (g',[g,e])\mapsto [g'g,e]. 
$$
This construction defines equivalence functors between the categories
of $G$-equivariant holomophic vector bundle on $X\times G/P$ and
$P$-equivariant holomorphic vector bundles on $X$. 

\bfsubsubsection{Homogeneous vector bundles}
\label{sub:induction-homogeneous}

A holomorphic homogeneous vector bundle on $G/P$ is
a $G$-equivariant holomorphic vector bundle on $G/P$. If we take
$X$ to be a point in the previous correspondence, we get an equivalence
between the category of holomorphic homogeneous vector bundles on
$G/P$ and the category of holomorphic representations of $P$. Now,
$X\times G/P$ is a `family' of flag varieties $\{ x\}\times G/P\cong
G/P$, parametrised by the points $x\in X$. Let $\cF$ be a
$G$-equivariant holomophic vector bundle on $X\times G/P$, and
$\cE=i^*\cF$ the corresponding holomorphic $P$-equivariant vector
bundle on $X$. Given $x\in X$, the restriction $\cF_x:=i_x^*\cF$ of
$\cF$ to the flag variety $i_x:G/P\cong\{ x\}\times G/P\hra X\times
G/P$ is a holomorphic homogeneous vector bundle on $G/P$. So $\cF$ is
also a `family' of holomorphic homogeneous vector bundles $\cF_x$ on
$G/P$. The holomorphic homogeneous vector bundle $\cF_x$ is in
correspondence with the holomorphic representation $\cE_x$ of $P$.

\bfsubsubsection{Irreducible representations of $P$}
\label{subsub:irred-P}

It is apparent from the constructions above that a detailed study of
the holomorphic representations of $P$ is important in
classifying the $G$-equivariant holomophic vector bundles on $X\times
G/P$. To study the holomorphic representations of $P$, it is
convenient to start with the irreducible ones. Throughout this paper,
$U$ is the unipotent radical of $P$, and $L$ is a (reductive) Levi
subgroup of $P$, so there is a semidirect decomposition $P=U\ltimes
L$, coming from a short exact sequence
$$
\begin{CD}
1  @>>>  U  @>>>  P @>>>  L  @>>>  1  .
\end{CD}
$$
It is an immediate consequence of Engels theorem that
a holomorphic representation $V$ of $P$ is irreducible if and only if
the action of $U$ on $V$ is trivial and $V$ is irreducible when
considered as a holomorphic representation of $L$. Therefore there is
a one-to-one correspondence between irreducible representations of $P$
and irreducible representations of its Levi subgroup $L$. 

\bfsubsubsection{Notation} 
\label{subsub:notation-groups}
Throughout this paper, $H$ is a Cartan
subgroup of $G$ such that $H\subset L$, and $\glp, \glu, \gll, \glh$
are the Lie algebras of $P,U,L,H$, respectively. The lattice of
integral weights of $H$, which parametrises the (isomorphism classes
of) irreducible representations of $H$, is denoted $\Lambda\subset
\glh^*$. We fix a fundamental chamber $\Lambda^+_P$ in $\Lambda$ for
the reductive Lie group $L$, so its elements are called integral {\em
dominant} weights of $P$. The fundamental chamber parametrises the
(isomorphism classes of) irreducible representations of $L$ (or $P$,
cf. \S \ref{subsub:irred-P}). (The fundamental chamber $\Lambda^+_P$
only depends on $P$, since any two Levi subgroups are conjugate, by
Malcev's theorem, cf. e.g. \cite{OV}.) Thus, given an integral
dominant weight for $P$, $\lambda\in\Lambda^+_P$, we fix an
irreducible representation $M_\lambda$ of $P$ (or $L$) (cf. \S
\ref{subsub:irred-P}), with highest weight $\lambda$. 

A holomorphic representation $\rho:P\ra GL(V)$ of $P$ on $V$ restricts
to holomorphic representations $\sigma:=\rho|_L:L\ra GL(V)$,
$\tau:=\rho|_U:U\ra GL(V)$ of $L$ and $U$ on $V$. Let $\sigma:L\ra
GL(V)$ and $\tau:U\ra GL(V)$ be holomorphic representations of $L$ and
$U$ on a vector space $V$. To determine when there exists a
holomorphic representation $\rho$ of $P$ whose restrictions to $L$ and
$U$ are $\sigma$ and $\tau$, it is convenient to use the
representation theory of quivers. Good references for details on
quivers and their representations are \cite{ARS,GR}.

\bfsubsubsection{Representations of quivers} 

A quiver, or directed graph, is a pair of sets $Q=(Q_0,Q_1)$
together with two maps $h,t:Q_1\ra Q_0$. The elements of $Q_0$
(resp. $Q_1$) are called the vertices (resp. arrows) of the
quiver. For each arrow $a\in Q_1$, the vertex $ta$ (resp. $ha$) is
called the tail (resp. head) of the arrow $a$, and the arrow $a$ is 
sometimes represented by $a:v\to v'$, with $v=ta, v'=ha$. 
We will not require $Q_0$ to be finite, but will require $Q$ to be
locally finite, i.e. $h^{-1}(v),t^{-1}(v)$ should be finite, for 
$v\in Q_0$.  A (non-trivial) path in $Q$ is a sequence $p=a_0\cdots
a_m$ of arrows $a_i\in Q_1$ which compose, i.e. with $ta_{i-1}=ha_i$
for $1\leq i\leq m$: 
\begin{equation}\label{eq:path}
p: \quad 
\stackrel{tp}{\bullet }\stackrel{a_m}{\lra }\bullet\stackrel{a_{m-1}}{\lra }\cdots\stackrel{a_0}{\lra }\stackrel{hp}{\bullet }
\end{equation}
The vertices $tp:=ta_m$ and $hp:=ha_0$ are called the tail and the
head of the path $p$.  The trivial path $e_v$ at $v\in Q_0$ consists
of the vertex $v$ with no arrows. A (complex) relation of a
quiver $Q$ is a formal finite sum $r=c_1 p_1 + \cdots + c_l p_l$ of
paths $p_1,\ldots ,p_l$ with coefficients $c_i\in\CC$. A {\em quiver
with relations} is a pair $(Q,\cK)$, where $Q$ is a quiver and $\cK$
is a set of relations of $Q$.
A linear representation $\bmR=(\bmV,\bvarphi)$ of
$Q$ is given by a collection $\bmV$ of complex vector spaces $V_v$, for
each vertex $v\in Q_0$, together with a collection $\bvarphi$ of
linear maps $\varphi_a:V_{ta}\ra V_{ha}$, for each arrow $a\in
Q_1$. We also require that $V_v=0$ for all but finitely many $v$. A
morphism $\bmf:\bmR\ra\bmR'$ between two representations
$\bmR=(\bmV,\bvarphi)$ and $\bmR'=(\bmV', \bvarphi')$ is given by
morphisms $f_v:V_v\ra V'_v$, for each $v\in Q_0$, such that
$\varphi'_a\circ f_{ta} =f_{ha}\circ\varphi_a$ for all $a\in Q_1$. 
Given a quiver representation $\bmR=(\bmV ,\bvarphi)$, a (non-trivial)
path (\ref{eq:path}) induces a linear map 
$$
\varphi(p):=\varphi_{a_0}\circ\cdots\circ\varphi_{a_m}:\, V_{tp}\lra V_{hp}.
$$
The linear map induced by the trivial path $e_v$ at $v\in Q_0$ is
$\varphi(e_v)=\id :V_v\ra V_v$. A linear representation
$\bmR=(\bmV,\bvarphi)$ of $Q$ is said to satisfy the relation $r=\sum_i
c_i p_i$ if $\sum_i c_i \varphi (p_i)=0$. Given a set $\cK$ of
relations of $Q$, a $(Q,\cK)$-module is a linear representation
of $Q$ satisfying the relations in $\cK$. The category of
$(Q,\cK)$-modules is clearly abelian.

\bfsubsection{The twisted nilradical representations}
\label{sec:A-B-psi}

This subsection provides the key technical point of the understanding of
representations of a complex Lie group $P$ as quiver representations. 
The basic ingredients are vector spaces $A_{\mu\lambda }, 
B_{\mu\lambda }$ and linear maps $\psi_{\mu\nu\lambda },
\psi_{\mu\lambda }$, for $\lambda,\mu,\nu$ integral dominant 
weights of $P$.

As a motivation for the definitions below, note that, if $V$ is a
$P$-module, with action of its Lie algebra $\glp$ given by
$\rho:\glp\ra\Aut(V)$, then the action of $L$ is specified by the
decomposition of $V$, as an $L$-module, into isotopical components
(since $L$ is reductive), and the rest of the $P$-action is
determined by a linear map $\tau=\rho|_{\glu}:\glu\ra\End(V)$. 
If $\sigma:\gll\ra\End(V)$ is the representation of $\gll$ induced by
$\rho$, it is clear that $\tau$ satisfies the following commutation relation
 $\tau([f,e])=[\sigma(f),\tau(e)]$  for $f\in\gll$ and
$e\in\glu$; and $\tau([e,e'])=[\tau(e),\tau(e')]$ for
$e,e'\in\glu$. The former condition means that $\tau$ is 
$L$-equivariant, and motivates the definition of the linear spaces
$A_{\mu\lambda }$, while the latter is more delicate, and motivates
the definition of the linear spaces $B_{\mu\lambda }$ and the linear
maps $\psi_{\mu\nu\lambda },\psi_{\mu\lambda }$.

\bfsubsubsection{Isotopical decompositions}\label{subsub:A-B}

Consider the nilradical $\glu\subset\glp$ as a representation
of $L$. Decompose the $L$-modules obtained by twisting the
exterior powers of the dual nilradical representation with the
irreducible $P$-modules,  $\glu^*\otimes M_\mu$,
$\wedge^2\glu^*\otimes M_\mu$, for  $\mu\in \Lambda^+_P$, into
irreducible components, as a representation of $L$:
$$
\glu^*\otimes M_\mu=\bigoplus_{\lambda\in \Lambda^+_P} A_{\mu\lambda }\otimes M_\lambda ,\quad
A_{\mu\lambda }:=(\glu^*\otimes\Hom(M_{\lambda },M_{\mu }))^L,
$$
$$
\wedge^2\glu^*\otimes M_\mu=\bigoplus_{\lambda\in \Lambda^+_P} B_{\mu\lambda }\otimes M_\lambda ,\quad
B_{\mu\lambda }:=(\wedge^2\glu^*\otimes\Hom(M_{\lambda },M_{\mu }))^L.
$$
\bfsubsubsection{Linear maps}\label{sub:linear-maps}
(a) For any $\lambda ,\mu ,\nu\in \Lambda^+_P$, define the
linear map
$$
\begin{array}{rcl}
\psi_{\mu\nu\lambda} : \left(\glu^*\otimes\Hom(M_\mu ,M_\nu )\right) \otimes \left(\glu^*\otimes\Hom(M_\nu ,M_\lambda)\right) &
\lra & \wedge^2\glu^*\otimes\Hom(M_{\lambda},M_{\mu})\\
a'\otimes a &\longmapsto & a'\wedge a
\end{array}
$$
(b) For any $\lambda ,\mu\in \Lambda^+_P$, define the linear map
$$
\psi_{\mu\lambda} :\glu^*\otimes\Hom(M_\lambda,M_\mu) 
\lra  \wedge^2\glu^*\otimes\Hom(M_{\lambda},M_{\mu})
$$
by sending $a\in\glu^*-\otimes\Hom(M_\lambda ,M_\mu)$ into the
$\Hom(M_\lambda ,M_\mu)$-valued 2-form on $\glu$ given by
$$
\psi_{\mu\lambda }(a)(e,e')=-a([e,e'])\quad {\rm ~for~}e,e'\in\glu.
$$
By the exterior product $a'\wedge a$ of part (a), of course we mean 
$\psi_{\mu\nu\lambda }(a'\otimes a)=a'\wedge a=(s'\wedge s)\otimes 
(f'\circ f)$ for $a=s\otimes f, a'=s'\otimes f'$, with
$s,s'\in\glu^*$ and $f\in\Hom(M_\lambda ,M_\nu), f'\in\Hom(M_\nu
,M_\mu)$.

The following lemma is straightforward. We shall need it in the
definition \S \ref{subsub:P-relations} of the relations $\cK$
associated to the complex Lie group $P$.

\begin{lemma}\label{lemma:psi}
If $\lambda ,\mu ,\nu\in \Lambda^+_P$, then $\psi_{\mu\nu\lambda
}(A_{\mu\nu }\otimes A_{\nu\lambda })\subset B_{\mu\lambda }$ and
$\psi_{\mu\lambda }(A_{\mu\lambda })\subset B_{\mu\lambda }$.
\end{lemma}

\bfsubsection{Quiver and relations associated to the Lie group $P$}
\label{sec:Q,K,P}

\bfsubsubsection{Quiver} \label{subsub:P-quiver}
Given $\lambda,\mu\in \Lambda^+_P$, let $\{
a^{(i)}_{\mu\lambda}|i=1,\ldots ,n_{\mu\lambda }\}$ be a basis of
$A_{\mu\lambda }$, with $n_{\mu\lambda }:=\dim A_{\mu\lambda }$. The
{\em quiver $Q$ associated to $P$} has as
vertex set, $Q_0=\Lambda^+_P$, i.e. the set of irreducible
representations of $P$, and as arrow set $Q_1=\{ a^{(i)}_{\mu\lambda
}|\lambda ,\mu\in Q_0,1\leq i\leq n_{\mu\lambda }\}$, the set of all
basis elements of the vector spaces $A_{\mu\lambda }$. The tail and
head maps $t,h:Q_1\ra Q_0$ are defined by
$$
t(a^{(i)}_{\mu\lambda })=\lambda ,\quad h(a^{(i)}_{\mu\lambda })=\mu .
$$
\bfsubsubsection{Relations} \label{subsub:P-relations}
Let $Q=(Q_0,Q_1)$ be as above. Let $\{ b^{(p)}_{\mu\lambda }| 1\leq
p\leq m_{\mu\lambda }\}$ be a basis of $B_{\mu\lambda }$, with
$m_{\mu\lambda }:=\dim B_{\mu\lambda }$, for $\lambda ,\mu\in
Q_0$. Expand $\psi_{\mu\nu\lambda}(a^{(j)}_{\mu\nu}\otimes a^
{(i)}_{\nu\lambda})\in B_{\mu\lambda }$ and
$\psi_{\mu\lambda}(a_{\mu\lambda}^{(k)})\in B_{\mu\lambda}$ in this
basis (cf. Lemma \ref{lemma:psi}), for $\lambda ,\mu ,\nu\in Q_0$:
\begin{equation}\label{eq:coefficientes-relations}
\psi_{\mu\nu\lambda}(a^{(j)}_{\mu\nu}\otimes a^{(i)}_{\nu\lambda})=\sum_{p=1}^{m_{\mu\lambda}}
c^{(j,i,p)}_{\mu\nu\lambda} ~b^{(p)}_{\mu\lambda},\quad
\psi_{\mu\lambda }(a^{(k)}_{\mu\lambda})=\sum_{p=1}^{m_{\mu\lambda}} c^{(k,p)}_{\mu\lambda} ~b^{(p)}_{\mu\lambda}.
\end{equation}
The {\em set of relations of the quiver $Q$ associated to $P$} is
$\cK=\{ r_{\mu\lambda}^{(p)}| \lambda,\mu\in Q_0, 1\leq p\leq
m_{\mu\lambda}\}$, where
$$
r_{\mu\lambda}^{(p)}=
\sum_{\nu\in Q_0}\sum_{i=1}^{n_{\nu\lambda}}\sum_{j=1}^{n_{\mu\nu}}c^{(j,i,p)}_{\mu\nu\lambda}  a^{(j)}_{\mu\nu} a^{(i)}_{\nu\lambda} +\sum_{k=1}^{n_{\mu\lambda}} c^{(k,p)}_{\mu\lambda} a^{(k)}_{\mu\lambda}.
$$
Note that in this definition $a^{(j)}_{\mu\nu} a^{(i)}_{\nu\lambda}$
does not mean the composition of
$a^{(i)}_{\nu\lambda}\in\glu^*\otimes\Hom(M_\lambda ,M_\nu)$ with
$a^{(j)}_{\mu\nu}\in\glu^*\otimes\Hom(M_\nu ,M_\lambda )$, but the
path that these two arrows define.

\bfsubsection{Correspondence between group and quiver representations} 
\label{sec:P-modules-quiver-rep}

\begin{theorem}\label{thm:equivalence-P-modules-quivers}
Let $Q$ and $\cK$ be the quiver and the set of relations associated
to the group $P$. There is an equivalence of
categories 
$$
\left\{\begin{array}{c}
{\rm finite~dimensional} \\
{\rm holomorphic~representations} \\
{\rm of~the~Lie~group~} P
\end{array}\right\}
\longleftrightarrow
\left\{\begin{array}{c}
{\rm finite~dimensional}\\
{\rm representations~of~the~quiver~} Q \\
{\rm ~satisfying~the~relations~} \cK
\end{array}\right\}.
$$
\end{theorem}

This equivalence was first proved by Bondal and Kapranov \cite{BK}
when $P$ is a Borel subgroup or the simple components of $G$
are in the series $A, D, E$. They were actually only interested in
homogeneous bundles over projective spaces, so these cases were enough
for their purposes. Bondal and Kapranov also gave a very simple and
explicit description of the quiver and the relations in these cases
(we collect their results about the quivers in Propositions
\ref{prop:Q-K-Borel} and \ref{prop:dim-arrows}).  Unfortunately, their
theorem is not always true with the simple definition of the relations
given in \cite{BK}, as shown by a counterexample found by Hille
\cite{Hl2}. Our definition of the quiver is precisely as in \cite{BK},
but the relations have been appropriately corrected. Hille also defined a quiver with
relations for any parabolic subgroup, and proved the corresponding
theorem of equivalence of categories in \cite{Hl2}. Hille uses a
different definition of the quiver associated to $P$, obtained by
removing certain arrows from the quiver in \cite{BK}  (although they
coincide when $U$ is abelian so the relations are quadratic). As
mentioned in the introduction, the quiver and relations that we obtain
seem to arise more naturally from the point of view of dimensional
reduction (cf. Theorem \ref{thm:corresp-equations}).\\

\proof
To prove the theorem, we shall define an equivalence functor from the
category of representations of $P$ into the category of
$(Q,\cK)$-modules. Given a representation $\rho:P\ra\Aut(V)$ of $P$
on a (finite-dimensional) complex vector space $V$, we obtain by
restriction a representation $\sigma=\rho|_L:L\ra\Aut(V)$ of $L\subset
P$ on $V$. Since $L$ is reductive, $\sigma$ decomposes into a direct sum
\begin{equation}\label{eq:module-decomposition}
V=\bigoplus_{\lambda\in Q_0}V_\lambda\otimes M_\lambda ,
\quad V_\lambda=\Hom_L(M_\lambda ,V).
\end{equation}
The vector spaces $V_\lambda$ are the multiplicity spaces, and have
trivial $L$-action. Let us fix a representation $\sigma :L\ra\Aut(V)$
of $L$ on $V$ as in (\ref{eq:module-decomposition}). We consider the
set $\Rep_P(\sigma)$ of (isomorphism classes of) representations
$\rho:P\ra\Aut(V)$ of $P$ on $V$ whose restriction to $L$ is
$\rho|_L=\sigma$. Apart from $L$, the rest of the $P$-module
structure of $V$ is given by the action $\tau=\rho|_{U}:U\ra
GL(V)$ of $U$ on $V$. 
Now, for the unipotent complex Lie group $U$, the map
$\exp:\glu\ra U$ is a $P$-equivariant isomorphism (of algebraic
varieties, cf. e.g. \cite[\S $3.3.6^\circ$, Theorem 7]{OV}), and for
any element $u=\exp e$ of $U$, with $e\in\glu$, $\tau(e)$ is a
nilpotent operator, so its exponential is a finite sum 
\begin{equation}\label{eq:tau-exp}
\tau(u)=\exp(\tau (e))=\sum_i \frac{1}{i!}\tau(e)^i,
\end{equation}
hence a polynomial in $\glu$. 
(We use the same symbol $\tau$ for the representation $d\tau
:\glu\ra\End(V)$ of the Lie algebra $\glu$ of $U$.)
The conditions for a linear map $\tau:\glu\ra\End(V)$ to define a
representation $\rho\in\Rep_P(\sigma)$ consist of two kinds of
commutation relations:

\noindent
{(CR1)}  $\quad \tau([f,e])=[\sigma(f),\tau(e)]$ for $f\in\gll$ and
$e\in\glu$;\\
\noindent
{(CR2)} $\quad \tau([e,e'])=[\tau(e),\tau(e')]$ for
$e,e'\in\glu$.

If we consider $\glu$ as an $L$-module, and $\End(V)$ as the
$L$-module obtained from the isotopical decomposition
(\ref{eq:module-decomposition}), then condition (CR1) simply means
that $\tau:\glu\ra\End(V)$ is a morphism of $L$-modules. In
other words, (CR1) is satisfied \IFF\ $\tau$ belongs to the
$L$-invariant part $W_1^{L}$ of the $L$-module
$$
W_1=\glu^*\otimes\End(V).
$$
(The subindex `1' accounts for condition (CR1).) The isotopical
decomposition \S \ref{subsub:A-B} applies to give
$$
W_1\cong\bigoplus_{\lambda ,\mu ,\nu\in Q_0} A_{\nu\lambda }
\otimes\Hom(V_{\lambda },V_{\mu })\otimes\Hom(M_{\nu },M_{\mu }),
$$
while Schur's lemma implies
$$
W_1^{L}\cong\bigoplus_{\lambda ,\mu\in Q_0}A_{\mu\lambda
}\otimes\Hom(V_{\lambda },V_{\mu }).
$$
Let $\bmV$ be a collection of linear spaces $V_\lambda$, for each
$\lambda\in Q_0$. There is a linear isomorphism between $W^{L}_1$ and
the space of quiver representations {\em into} $\bmV$, 
$$
\RRR(Q,\bmV)=\bigoplus_{a\in Q_1}\Hom(V_{ta},V_{ha}),
$$
i.e. the space of representations $\bmR=(\bmV,\bvarphi)$ with fixed
$\bmV$. This isomorphism takes any 
\begin{equation}\label{eq:varphi}
\tau=\sum_{\lambda ,\mu\in Q_0}\sum_{i=1}^{n_{\mu\lambda }} a^{(i)}_{\mu\lambda
}\otimes\varphi^{(i)}_{\mu\lambda}\in \bigoplus_{\lambda
,\mu\in Q_0} A_{\mu\lambda }\otimes\Hom(V_{\lambda },V_{\mu }),
\end{equation}
into the representation $\bmR=(\bmV ,\bvarphi)$ 
given by the morphisms $\varphi_a=\varphi^{(i)}_{\mu\lambda
}:V_{\lambda }\ra V_{\mu }$, for $a=a^{(i)}_{\mu\lambda }$.

Next we consider how condition (CR2) translates into
$\RRR(Q,\bmV)$. Define the $L$-module
$$
W_2:=\wedge^2\glu^*\otimes\End(V)
$$
(the subindex `2' accounts for condition (CR2)) and the (non-linear) map
$$
\psi :W_1 \lra  W_2
$$
by
$$
\psi(\tau)(e,e')=[\tau(e),\tau(e')]-\tau([e,e']),
\quad {\rm for~} \tau\in W_1, ~ e,e'\in\glu.
$$
Then $\tau$ satisfies condition
(CR2) \IFF\ $\psi(\tau)=0$, so the
space of representations $\rho$ of $P$ with $\rho|_L\cong\sigma$ is
$\Rep_P(\sigma)\cong\psi^{-1}(0)\cap W_1^{L}$; note that this is
contained in $W_1^{L}\cong\RRR(Q,\bmV)$. In order to express this
condition in terms of the relations $r^{(p)}_{\mu\lambda}$, we rewrite
$\psi$ using the following linear maps $\psi_1,\psi_2$. The first
linear map is 
$$
\begin{array}{rccl}
\psi_1 : & W_1\otimes W_1 & \lra & W_2 \\
& \tau'\otimes\tau &\longmapsto & \tau'\wedge\tau
\end{array}
$$
where by the exterior product $\tau'\wedge\tau$, we mean
$\psi_1(\tau'\otimes\tau)=\tau'\wedge\tau=(s'\wedge s)\otimes
(f'\circ f)$ for $\tau=s\otimes f,\tau'=s'\otimes f'$, with
$s,s'\in\glu^*$ and $f,f'\in\End(V)$ (this map shoud be compared
with \ref{sub:linear-maps}(a)). The second linear
map
$$
\psi_2:W_1\lra W_2,
$$
is given by
$$
\psi_2(\tau)(e,e')=-\tau([e,e']), \quad {\rm for~} \tau\in W_1, ~~ e,e'\in\glu^*
$$
(compare with \ref{sub:linear-maps}(b)). Obviously
\begin{equation}\label{eq:psi}
\psi(\tau)=\psi_1(\tau\otimes\tau)+\psi_2(\tau), \quad {\rm for~} \tau\in W_1.
\end{equation}
Note that 
$$
W_2^{L}= (\wedge^2\glu^*\otimes\End(V))^{L}
\cong \bigoplus_{\lambda,\mu\in Q_0} B_{\mu\lambda}\otimes\Hom(V_{\lambda},V_{\mu}).
$$

When $\tau\in W_1^{L}$ is given by (\ref{eq:varphi}), linearity of
$\psi_1,\psi_2$ allows one to obtain $\psi(\tau)$:
\begin{equation}\label{eq:psi_1}\begin{split}
\psi_1(\tau\otimes\tau)&=\sum_{\lambda,\nu,\nu',\mu\in Q_0}
\sum_{i=1}^{n_{\nu\lambda}}\sum_{j=1}^{n_{\mu\nu'}}
\psi_1((a^{(j)}_{\mu\nu'}\otimes\varphi^{(j)}_{\mu\nu'})
\otimes(a^{(i)}_{\nu\lambda}\otimes\varphi^{(i)}_{\nu\lambda})) \\
&=\sum_{\lambda,\nu,\mu\in Q_0}\sum_{i=1}^{n_{\nu\lambda
}}\sum_{j=1}^{n_{\mu\nu }}
\psi_{\mu\nu\lambda}(a^{(j)}_{\mu\nu}\otimes
a^{(i)}_{\nu\lambda})\otimes(\varphi^{(j)}_{\mu\nu}\circ\varphi^{(i)}_{\nu\lambda})\\
&=\sum_{\lambda,\mu\in
Q_0}\sum_{p=1}^{m_{\mu\lambda}}b^{(p)}_{\mu\lambda}\otimes
\left(\sum_{\nu\in Q_0}\sum_{i=1}^{n_{\nu\lambda}}\sum_{j=1}^{n_{\mu\nu}}c^{(j,i,p)}_{\mu\nu\lambda}\varphi^{(j)}_{\mu\nu}\circ\varphi^{(i)}_{\nu\lambda}\right),
\end{split}\end{equation}
and 
\begin{equation}\label{eq:psi_2}\begin{split}
\psi_2(\tau)
&=\sum_{\lambda,\mu\in Q_0}\sum_{i=1}^{n_{\mu\lambda}}
\psi_2(a_{\mu\lambda}^{(i)}\otimes\varphi^{(i)}_{\mu\lambda}) \\
&=\sum_{\lambda,\mu\in Q_0}\sum_{i=1}^{n_{\mu\lambda}}
\psi_{\mu\lambda}(a_{\mu\lambda}^{(i)})\otimes\varphi^{(i)}_{\mu\lambda}
=\sum_{\lambda,\mu\in Q_0}\sum_{p=1}^{m_{\mu\lambda}} b^{(p)}_{\mu\lambda}\otimes 
\left(\sum^{n_{\mu\lambda}}_{k=1} c^{(k,p)}_{\mu\lambda}\varphi^{(k)}_{\mu\lambda}\right),
\end{split}\end{equation}
so that (\ref{eq:psi}), (\ref{eq:psi_1}) and (\ref{eq:psi_2}) imply
\begin{equation}\label{eq:CR2-psi}
\psi(\tau)=
\sum_{\lambda,\mu\in Q_0}\sum_{p=1}^{m_{\mu\lambda}}b^{(p)}_{\mu\lambda}\otimes
\left(\sum_{\nu\in Q_0}\sum_{i=1}^{n_{\nu\lambda}}\sum_{j=1}^{n_{\mu\nu}}c^{(j,i,p)}_{\mu\nu\lambda}\varphi^{(j)}_{\mu\nu}\circ\varphi^{(i)}_{\nu\lambda} +\sum_{k=1}^{n_{\mu\lambda}} c^{(k,p)}_{\mu\lambda}\varphi^{(k)}_{\mu\lambda}\right).
\end{equation}
This equation gives the relations of the quiver $Q$ which realise
condition (CR2).
\qed

It is worth remarking that the inclusion
$\wedge^2\glu^*\subset\glu^*\otimes\glu^*$ gives 
\begin{multline*}
B_{\mu\lambda}
\subset
\left( M_{\mu}\otimes\glu^*\otimes\glu^*\otimes M^*_{\lambda }\right)^{L}
 \cong \left( \bigoplus_{\nu\in Q_0} A_{\mu\nu }\otimes M_{\nu}\otimes \bigoplus_{\nu'\in Q_0}A_{\nu'\lambda}\otimes M^*_{\nu'} \right)^{L}
\cong \bigoplus_{\nu\in Q_0} A_{\mu\nu}\otimes A_{\nu\lambda}.
\end{multline*}
Therefore, given $\lambda,\mu\in Q_0$, a necessary condition for
$B_{\mu\lambda}\neq 0$ is that there exists a $\nu\in Q_0$ with
$A_{\mu\nu}\otimes A_{\nu\lambda}\neq 0$, i.e. there must be a path
$\lambda\ra\nu\ra\mu$ in the quiver $Q$.

\begin{remark}{\rm
The quiver with relations $(Q,\cK)$, as defined in \S \ref{sec:Q,K,P},
can be associated to any algebraic group $P$ (not necessarily a
parabolic subgroup) over any algebraically closed field of
characteristic zero. Actually, any such $P$ admits a semidirect
decomposition $P=U\ltimes L$, as in \S \ref{subsub:irred-P}, where $U$
is the unipotent radical of $P$, and $L$ is any (reductive) Levi subgroup
$P$ (a proof of the existence of $L$ can be found e.g. in
\cite[Chapter 6]{OV}). Of course, Theorem
\ref{thm:equivalence-P-modules-quivers} generalises to any such $P$ as
well, so that the category of representations of $(Q,\cK)$ is
equivalent to the category of rational representations of $P$.  
}\end{remark}

As a result of \S \ref{sub:induction-homogeneous} and Theorem
\ref{thm:equivalence-P-modules-quivers}, it follows that:

\begin{corollary}\label{coro:equivalence-homvb-quivers}
Let $G$ be a connected complex Lie group and $P\subset G$ a
parabolic subgroup. Let $Q$ and $\cK$ be the quiver and the set of
relations associated to $P$. There is an equivalence of categories 
$$
\left\{\begin{array}{c}
{\rm holomorphic~homogeneous} \\
{\rm vector~bundles~on~} G/P
\end{array}\right\}
\longleftrightarrow
\left\{\begin{array}{c}
{\rm finite~dimensional}\\ 
{\rm representations~of~the~quiver~} Q \\
{\rm ~satisfying~the~relations~in~} \cK
\end{array}\right\}.
$$
\end{corollary}

\bfsubsection{Holomorphic filtrations of homogeneous bundles}
\label{sec:flags}

Let $Q$ be the quiver associated to $P$. This quiver has no {\em
oriented cycles}, i.e. there are no paths of length $>0$ in $Q$ whose
tail and head coincide. To prove this, we first introduce the notion
of $\Sigma$-height, adapting \cite[\S 10.1]{Hu} (see also
\cite{CSS}). An essential ingredient will be that $P$ is a {\em
parabolic} subgroup of $G$. We should mention that Hille \cite{Hl2} 
has defined a function similar to our $\Sigma$-height, that he
calls a level funtion, for the quivers that he associates to the
parabolic subgroups. \\

\bfsubsubsection{Notation} \label{subsub:notation-S-Sigma}
Choose a system $\cS$ of simple roots of $\glg$ with respect to
$\glh$, such that all the {\em negative} roots of $\glg$ with respect
to $\glh$ are roots of $\glp$. 
Let $\Sigma$ be the set of non-parabolic simple roots, i.e. it
consists of those simple roots of $\glg$ which are not 
roots of $\glp$ (in terms of $\gll$, this means that 
$\cS\setminus\Sigma$ is a system of simple roots of $\gll$, while the
roots in $\Sigma$ are not roots of $\gll$).  
Given a $P$-module $V$, let $\Delta(V)$ be its set of
weights with respect to $\glh$, so 
\begin{equation}\label{eq:V-weight-decomposition}
V=\bigoplus_{\lambda\in\Delta(V)} V^{\lambda },
\end{equation}
where $V^{\lambda }:=\{ v\in V| h\cdot v=\lambda(h)v {\rm ~for~}
h\in\glh\}$ for $\lambda\in\Lambda$, and $\Delta (V):=\{\lambda\in 
\Lambda|V^\lambda\neq 0\}$. 

\bfsubsubsection{$\Sigma$-height} \label{subsub:sigma-height}
Any integral weight $\mu\in\Lambda$ admits a decomposition
$\mu=\sum_{\alpha\in\cS } n_{\alpha }\alpha$, with $n_{\alpha }
\in\QQ$. We define the $\Sigma$-{\em height} of $\mu$ as the rational
number 
$$
\hts(\mu)=\sum_{\alpha\in\Sigma }n_{\alpha}.
$$ 

\begin{lemma}\label{lemma:sigma-height}
If the quiver $Q$ has an arrow $\lambda\ra\mu$, then
$\hts(\lambda)>\hts(\mu)$. Therefore, the quiver $Q$ is directed.
\end{lemma}

\proof
If the irreducible $L$-modules $M_{\lambda }^*, M_{\mu }, \glu^*$
have weight space decompositions 
$$
M_{\mu }=\bigoplus_{\mu'\in\Delta(M_{\mu})} M^{\mu'}_{\mu} ,\quad
M^*_{\lambda }=\bigoplus_{\lambda'\in\Delta(M_{\lambda})}
(M_{\lambda}^*)^{-\lambda'} ,
\quad\glu^*\cong\bigoplus_{\gamma\in\Delta(\glu)}\glg^{-\gamma } ,
$$
(see \eqref{eq:V-weight-decomposition} above), then $H\subset L$ implies
$$
A_{\mu\lambda }=(M_{\mu }\otimes\glu^*\otimes M_{\lambda }^*)^L \subset 
(M_{\mu }\otimes\glu^*\otimes M_{\lambda }^*)^H \cong
\bigoplus_{(\mu',\lambda')\in\Delta(M_\mu ,M_{\lambda })}
M^{\mu'}_{\mu}\otimes\glg^{-(\mu'-\lambda')}\otimes (M^*_{\lambda
})^{-\lambda'},
$$
where $\Delta(M_\mu ,M_\lambda )$ is the set of pairs
$(\mu',\lambda')$, with $\mu'\in\Delta(M_\mu)$,
$\lambda'\in\Delta(M_\lambda)$ such that $\mu'-\lambda'\in\Delta(\glu)$. 
If $A_{\mu\lambda }\neq 0$, there are
$\gamma\in\Delta(\glu),\lambda'\in\Delta(M_{\lambda}),\mu'\in\Delta(M_{\mu})$
with $\mu'-\lambda'=\gamma$. But
$\hts(\lambda')=\hts(\lambda),\hts(\mu')=\hts(\mu)$, so 
$\hts(\mu)-\hts(\lambda)=\hts(\mu-\lambda)=\hts(\gamma)<0$.
\qed

\bfsubsubsection{A total order in $Q_0$}\label{subsub:total-order} 
By Lemma \ref{lemma:sigma-height}, defining $\lambda>\mu$ if
$\hts(\lambda)>\hts(\mu)$, for each $\lambda,\mu\in Q_0$ provides a
partial order in $Q_0$. We now enlarge this partial order to get a
total order in $Q_0$: For each $q\in\QQ$ with
$ht^{-1}_{\Sigma}(q)\neq\emptyset$, we choose a total order $(<)_q$ on
$ht^{-1}_{\Sigma}(q)$, and define a total order in $Q_0$ by  
saying that $\mu<\lambda$, for each $\lambda,\mu\in Q_0$, if either
$\hts(\mu)<\hts(\lambda)$, or $\lambda,\mu\in ht^{-1}_\Sigma(q)$ for
some $q\in\QQ$ and $\mu(<)_q\lambda$. Thus, if there is an arrow
$\lambda\ra\mu$, then $\lambda>\mu$.

\begin{proposition}\label{prop:isotopical-flag}
Let $V$ be a $P$-module with isotopical decomposition, as an
$L$-module,
\begin{equation}\label{eq:V-weight-2}
V=\bigoplus_{\lambda\in Q_0(V)}V_\lambda\otimes M_\lambda , \quad
V_\lambda=\Hom_{L }(M_\lambda ,V),
\end{equation}
where $Q_0(V)\subset Q_0$ is a finite set. Let us list the set of
vertices in $Q_0(V)$ in ascending order as  $Q_0(V)=\{
\lambda_0,\lambda_1,\ldots ,\lambda_m\}$,
$\lambda_0<\lambda_1<\cdots<\lambda_m$. Then $V$ admits a flag 
of $P$-submodules with completely reducible quotients:
\begin{equation}\begin{gathered}\label{eq:flag-P-modules}
\bmV_{(\leq\bullet)}: 0\subset V_{(\leq 0)}\subset V_{(\leq 1)}\subset\cdots\subset V_{(\leq m)}=V\\
V_{(\leq s)}/V_{(\leq s-1)}\cong V_{\lambda_s}\otimes M_{\lambda_s}
\end{gathered}\end{equation}
\end{proposition}

(In fact, one can prove that if $V$ is indecomposable as a
$P$-module, then $\hts(\lambda_s)-\hts(\lambda_{s-1})$ is one or zero
for all $1\leq s\leq m$.) 

\proof From the proof of Theorem
\ref{thm:equivalence-P-modules-quivers}, we see that the $L$-module
structure of $V$ is given by an isotopical decomposition, as an
$L$-module,
$
V=\bigoplus_{s=0}^m V_{\lambda_s}\otimes M_{\lambda_s} 
$,
where $V_{\lambda_s}\neq 0$ for $0\leq s\leq m$, while the
$\glu$-structure is given by an $L$-invariant morphism  
$$
\tau=\sum_{0\leq s,s'\leq m}\sum_{i=1}^{n_{\lambda_{s'}\lambda_s }}
a^{(i)}_{\lambda_{s'}\lambda_s }\otimes\varphi^{(i)}_{\lambda_{s'}\lambda_s}
$$
of $W=\glu^*\otimes\End(V)$. 
Define the flag $\bmV_{(\leq \bullet)}$ of vector spaces as in
\eqref{eq:flag-P-modules}, where 
$
V_{(\leq s)}=\bigoplus_{s'=0}^s V_{\lambda_{s'}}\otimes M_{\lambda_{s'}}
$.
This is obviously a flag of $L$-modules. To prove that it is a flag of
$P$-modules as in \eqref{eq:flag-P-modules}, it is enough to see that
$\glu$ takes $V_{(\leq s)}$ into $V_{(\leq s-1)}$, for $1\leq s\leq
m$. Let $e\in\glu$. The action of $e$ on $V$ is given by
$\tau(e)\in\End(V)$, where 
$$
\tau(e)=\sum_{0\leq s,s'\leq m}\sum_{i=1}^{n_{\lambda_{s'}\lambda_s }}
a^{(i)}_{\lambda_{s'}\lambda_s }(e)\otimes\varphi^{(i)}_{\lambda_{s'}\lambda_s}
$$
with $a^{(i)}_{\lambda_{s'}\lambda_s}(e)\in\Hom(M_{\lambda_{s'}},
M_{\lambda_s})$ for $0\leq s,s'\leq m$. By part (c) of Lemma
\ref{lemma:sigma-height}, this is zero unless $s'<s$, so that
$\tau(e)$ takes $V_{\lambda_s}\otimes M_{\lambda_s}$ into 
$\oplus_{s'<s} V_{\lambda_{s'}}\otimes M_{\lambda_{s'}}$.
\qed

\bfsubsubsection{Notation}\label{subsub:notation-O-lambda} 

Given an irreducible representation $M_\lambda$ of $P$, corresponding
to an integral dominant weight $\lambda\in\Lambda^+_P$,
$\cO_{\lambda}:=G\times_P M_\lambda$ is the induced irreducible
holomorphic homogeneous vector bundle on $G/P$.

\begin{corollary}
Any holomorphic homogeneous vector bundle $\cF$ on $G/P$ admits a
filtration of holomorphic homogeneous vector subbundles $\cF_s$, with
completely reducible quotients with respect to the $G$-action,
\begin{equation}\begin{gathered}\label{eq:flag-hom-vb}
\bcF: 0\hra \cF_{0}\hra \cF_{1}\hra\cdots\hra\cF_{m}=\cF,\\
\cF_{s}/\cF_{s-1}\cong V_{\lambda_s}\otimes\cO_{\lambda_s},
\quad 1\leq s\leq m,
\end{gathered}\end{equation}
for some dominant integral weights $\lambda_s$ with
$\lambda_0<\lambda_1<\cdots<\lambda_m$, where $V_{\lambda_s}$ are
vector spaces determined by $\cF$.  
\end{corollary}

\proof 
This result follows from \S \ref{sub:induction-homogeneous} and
Proposition \ref{prop:isotopical-flag}. First, $\cF\cong G\times_P V$,
where $V_o=\cF_o$, the fibre at the base point $o=P\in G/P$, is a
representation of $P$. Then $V$ admits the isotopical decomposition
\eqref{eq:V-weight-2} and the flag of $P$-submodules
\eqref{eq:flag-P-modules}. Let $\cF_s=G\times_P V_{(\leq s)}$ be the
homogeneous holomorphic vector subbundle of $\cF$ induced by $V_{(\leq
s)}$, for $0\leq s\leq m$. Then $V_{(\leq s)}/V_{(\leq s-1)}\cong
V_{\lambda_s}\otimes M_{\lambda_s},$ implies $\cF_{s}/\cF_{s-1}\cong
V_{\lambda_s}\otimes \cO_{\lambda_s}$.
\qed

\subsection{Examples}
\label{sub:examples-hb}

In this subsection we shall present explicit expressions for the
quiver and the relations corresponding to some flag varieties.
We would like to remark that the problem of classification of the
quivers with relations associated to all parabolic subgroups is a
subtle one. At present a complete classification is not known, and it
would require a deep study which is out of the scope of this paper.
To obtain the quiver one only has to evaluate dimension formulas for
the vector spaces $A_{\mu\lambda}$. However, to obtain the relations,
as given in \S \ref{sec:Q,K,P}, is more difficult, 
since one has to choose bases of the vector spaces $A_{\mu\lambda }$
and $B_{\mu\lambda }$, and express the linear maps $\psi_{\mu\lambda}$,
$\psi_{\mu\nu\lambda}$ in these bases, as in \S \ref{subsub:P-relations}. 

\bfsubsubsection{Quiver and relations for Borel subgroups}
\label{subsub:Borel}

Let us assume that $B=P$ is a Borel subgroup of $G$. The set of
integral dominant weights of $B$ is precisely the weight lattice
$\Lambda\cong\ZZ^{{\rm rank}(G)}$ of integral weights (cf. \S \ref{subsub:notation-groups}),
for the Cartan subgroup $H$ is a Levi subgroup of $B$. 
Let $\Delta$ be the set of roots of $(\glg,\glh)$, and for
$\alpha\in\Delta$, let $\glg^{\alpha}$ be the root subspace of $\glg$
corresponding to $\alpha$. We choose the sets $\Delta_+,
\Delta_-\subset\Delta$ of positive and negative roots, so that
the Lie algebra of $B$ is $\glb=\glh\oplus\oplus_{\alpha\in\Delta_-} 
\glg^\alpha$. For each $\alpha\in\Delta$, let $e_\alpha\in\glg^\alpha$
be the corresponding Chevalley generator (see e.g. \cite[\S 25.21]{Hu}). 
Let $N_{\alpha \beta }\in\ZZ$, for $\alpha,\beta\in\Delta$, be the
coefficients defined by the condition $[e_\alpha ,e_\beta ]=
N_{\alpha\beta }e_{\alpha+\beta }$ if $\alpha+\beta\in\Delta$, and
$N_{\alpha\beta }=0$ if $\alpha+\beta\notin\Delta$. 

\begin{proposition}\label{prop:Q-K-Borel}
The quiver with relations $(Q,\cK)$ associated to $B$ is given as follows. 
\begin{enumerate}\item[(1)]
The vertex set $Q_0$ is the weight lattice $\Lambda\cong\ZZ^{{\rm rank}(G)}$. 
\item[(2)] The arrow set $Q_1$ consists of the arrows
$a^{(\gamma)}_{\mu\lambda}:\lambda\to\mu$, for 
$\lambda,\mu\in\Lambda$, with $\gamma=\mu-\lambda \in\Delta_-$. 
\item[(3)]
The set $\cK$ consists of the relations
$r^{(\gamma,\gamma')}_{\lambda\mu}=a^{(\gamma')}_{\mu\nu}a^{(\gamma)}_{\nu\lambda}
-a^{(\gamma)}_{\mu\nu'}a^{(\gamma')}_{\nu'\lambda}-N_{\gamma \gamma' }
a^{(\gamma+\gamma')}_{\mu,\lambda }$, for  $\lambda,\mu\in \Lambda$
and $\gamma,\gamma'\in\Delta_-$, with $\gamma+\gamma'=\mu-\lambda$ and
$\gamma\neq\gamma'$, where $\nu:=\lambda+\gamma$ and
$\nu':=\lambda+\gamma'$.  
\end{enumerate}\end{proposition}

(The term $a^{(\gamma')}_{\mu\nu}a^{(\gamma)}_{\nu\lambda }$ (resp. $a^{(\gamma)}_{\mu\nu'}
a^{(\gamma')}_{\nu'\lambda}$; $N_{\gamma ,\gamma'}a^{(\gamma+\gamma')}_{\mu ,\lambda}$)  is ignored
in the definition of $r^{(\gamma,\gamma')}_{\lambda\mu}$, 
whenever the basis vectors $a^{(\gamma')}_{\mu\nu}$ or $a^{(\gamma)}_{\nu\lambda }$
(resp. $a^{(\gamma)}_{\mu\nu'}$ or $a^{(\gamma')}_{\nu'\lambda}$; $a^{(\gamma+\gamma')}_{\mu ,\lambda}$) do
not make sense.)

\proof
Part (1) is obvious. 
The nilpotent radical of $\glb$ is $\glu=\oplus_{\gamma\in\Delta_-}
\glg_\gamma$, with basis $\{ e_\gamma| \gamma\in\Delta_-\}$. Let $\{
e^\gamma|\gamma\in\Delta_-\}$ be its dual basis.
Since the Levi subgroup $H\subset B$ is abelian, its
irreducible representations $M_\lambda$, for $\lambda\in\Lambda$,
are one dimensional. Let $v_\lambda$ be a basis vector of $M_\lambda$,
and $v^\lambda\in M^*_\lambda$ be its dual basis vector.
Let $A_{\mu\lambda}=(\glu^*\otimes\Hom(M_\lambda,M_\mu))^{H}$, for
$\lambda,\mu\in\Lambda$, as in \S \ref{subsub:A-B}. The weight of
$e^\gamma$ with respect to $\glh$ is $-\gamma$. Thus, if
$\mu-\lambda\not\in\Delta_-$, then $A_{\mu\lambda}=0$, while if
$\gamma:=\mu-\lambda\in\Delta_-$, then $A_{\mu\lambda}$ is one
dimensional, with basis vector 
$$
a^{(\gamma)}_{\mu\lambda}:=e^\gamma\otimes v_\mu\otimes v^\lambda.
$$
This proves part (2). To get part (3), let $B_{\mu\lambda}=(\wedge^2\glu^*\otimes
\Hom(M_\lambda,M_\mu))^{H}$, as in \S \ref{subsub:A-B}. Let $<$ be a
total order for the set $\Delta_-$. A basis of $\wedge^2\glu^*$ is $\{
e^\gamma\wedge e^{\gamma'} | \gamma,\gamma'\in\Delta_-,\gamma<\gamma'\}$. 
The weight of $e^\gamma\wedge e^{\gamma'}$ with respect
to $\glh$ is $-\gamma-\gamma'$. Thus, if $\mu-\lambda\neq\gamma+\gamma'$ 
for any $\gamma,\gamma'\in\Delta_-$, then
$B_{\mu\lambda}=0$, while if $\mu-\lambda=\gamma+\gamma'$ for some 
$\gamma,\gamma'\in\Delta_-$ with $\gamma<\gamma'$, then a basis of
$B_{\mu\lambda}$ is $\{ b^{(\gamma,\gamma')}_{\mu\lambda} |
\gamma,\gamma'\in\Delta_-,\gamma<\gamma', \gamma+\gamma'
=\mu-\lambda\}$, where 
$$
b^{(\gamma,\gamma')}_{\mu\lambda}:=(e^\gamma\wedge e^{\gamma'})\otimes
v_\mu\otimes v^\lambda.  
$$
To express $\psi_{\mu\nu\lambda}$ and
$\psi_{\mu\lambda}$ in the bases $\{a^{(\gamma)}_{\mu\lambda}\}$ and $\{
b^{(\gamma,\gamma')}_{\mu\lambda}\}$, let $\lambda,\mu\in\Lambda$ 
be such that $B_{\mu\lambda}\neq 0$. Hence, there are
$\gamma,\gamma'\in\Delta_-$, with $\gamma<\gamma'$, 
$\mu-\lambda=\gamma+\gamma'$. Let $\nu=\lambda+\gamma$,
$\nu'=\lambda+\gamma'$. Then  
$$
\psi_{\mu\nu\lambda}(a^{(\gamma')}_{\mu\nu}\otimes a^{(\gamma)}_{\nu\lambda})
=b^{(\gamma,\gamma')}_{\mu\lambda},
\quad
\psi_{\mu\nu'\lambda}(a^{(\gamma)}_{\mu\nu'}\otimes a^{(\gamma')}_{\nu'\lambda})
=-b^{(\gamma,\gamma')}_{\mu\lambda}.
$$
Let $\lambda,\mu\in\Lambda$ with $A_{\mu\lambda}\neq 0$.
Then $a^{(\phi)}_{\mu\lambda}=e^\phi \otimes v_\mu\otimes v^\lambda$,
where $\phi:=\mu-\lambda\in\Delta_-$, so 
\begin{equation}\label{eq:Borel-psi_mu-lambda}
\psi_{\mu\lambda}(a^{(\phi)}_{\mu\lambda})(e_\epsilon,e_{\epsilon'})
=-a^{(\phi)}_{\mu\lambda}(N_{\epsilon,\epsilon'}e_{\epsilon+\epsilon'})
=-\delta_{\epsilon+\epsilon'}^{\phi}
N_{\epsilon,\epsilon'}v_\mu\otimes v^\lambda ,
\quad {\rm for~} \epsilon,\epsilon'\in\Delta_-,
\end{equation}
($\delta_{\epsilon+\epsilon'}^{\phi}$ is Kronecker's delta). If
$\gamma,\gamma'\in\Delta_-$, with $\gamma<\gamma'$, are such that
$\mu-\lambda=\gamma+\gamma'$, then 
$
-N_{\gamma\gamma'}b^{(\gamma,\gamma')}_{\mu\lambda}(e_\epsilon,e_{\epsilon'})
=-N_{\epsilon\epsilon'} \delta^\gamma_\epsilon
\delta^{\gamma'}_{\epsilon'} ~ v_\mu\otimes v^\lambda.
$
Comparing this with \eqref{eq:Borel-psi_mu-lambda}, we see that 
$$
\psi_{\mu\lambda}(a^{(\phi)}_{\mu\lambda})
= -\sum_{(\gamma,\gamma')\in\Delta_{\mu\lambda}}
N_{\gamma\gamma'} b^{(\gamma,\gamma')}_{\mu\lambda},
$$
where $\Delta_{\mu\lambda}:=\{(\gamma,\gamma')| \gamma,\gamma'\in
\Delta_-, ~\gamma+\gamma=\mu-\lambda,\gamma<\gamma'\}$. 
It follows from \S \ref{subsub:P-relations} that the relations
$r^{(\gamma,\gamma')}_{\mu \lambda}$ in $\cK$ are as given
in part (3), for $\lambda,\mu\in\Lambda$ and $(\gamma,\gamma')\in
\Delta_{\mu\lambda}$. For $\gamma>\gamma'$ we get the negative of
these $r^{(\gamma,\gamma')}_{\mu\lambda}=-r^{(\gamma,\gamma')}_{\mu
\lambda}$, since $N_{\gamma',\gamma}=-N_{\gamma,\gamma'}$. 
So the relations corresponding to $\gamma>\gamma'$ do not provide 
more constraints on the representations. This proves part (3).
\qed

The previous proposition has also been proved in \cite{BK}. Since
different  authors describe generally the quiver with relations in  different ways,
we include a proof based on our definitions. 

\bfsubsubsection{Homogeneous bundles on products of complex projective lines} 
\label{subsub:hb-prod-P1}

The product of $N$ complex projective lines $\PP^1$ can be written as
a quotient $(\PP^1)^N=G/P$ of groups
$$
G=\prod_{\alpha=1}^N G_\alpha,\quad P=\prod_{\alpha=1}^N
P_\alpha,\quad 
{\rm where~} G_\alpha=SL(2,\CC), {\rm ~and~}
P_\alpha=\begin{pmatrix} * & 0 \\ * & * \end{pmatrix}\subset G_\alpha .
$$
Since $P\subset G$ is a Borel subgroup, its associated quiver with
relations $(Q,\cK)$ is given by Proposition \ref{prop:Q-K-Borel}. Hence, 
the vertex set is $Q_0=\ZZ^N$. 
We easily see that the roots of $\glu$ are $\gamma_k=-2L_k$, for $1\leq
k\leq N$, where $\{ L_1,\ldots,L_N\}$ is the dual of the standard basis of
$\glh\cong\CC^N$. Thus, the arrow set is $Q_1=\{ a^{(i)}_\lambda|\lambda\in 
\ZZ^{N}, 1\leq i\leq N\}$, with $a^{(i)}_\lambda:\lambda\to\lambda-2L_i$. 
For any arrow $a^{(i)}_{\lambda}:\lambda\to\mu=\lambda-2 L_i$, 
we have $\sum_{i=1}^N\mu_i=\sum_{i=1}^N\lambda_i -2$. Thus, the quiver
decomposes into two full subquivers $Q^{(h)}$, for $h=0,1$, whose
vertices $\lambda$ satisfy $\sum_{i=1}^N\lambda_i\equiv h\mod{2}$.  For the 
example $G/P=\PP^1\times\PP^1$, the picture of any connected
component is given in Fig. 1 below.
Since $\glu$ is abelian, the coefficients $N_{\gamma,\gamma'}$, in
Proposition \ref{prop:Q-K-Borel}(3), are zero. Thus, the relations
$r^{(j,i)}_\lambda$ of $P$ are parametrised by a vertex $\lambda$ and
a basis vector of $\wedge^2\CC^N$. They are given by 
$$
r^{(j,i)}_\lambda=a^{(j)}_{\lambda-2 L_i}
a^{(i)}_\lambda-a^{(i)}_{\lambda-2 L_j} a^{(j)}_\lambda ,
\quad {\rm for~} 1\leq i< j\leq N.
$$
Summarizing, {\em the category of $G$-equivariant holomorphic
vector bundles on $G/P=(\PP^1)^N$ is equivalent to the category of
commutative diagrams on the quiver $Q$}.
$$
\begin{array}{ccc}
\setlength{\unitlength}{.8mm}
\begin{picture}(62,42)(-1,-1)
\thicklines
\matrixput(0,0)(10,0){7}(0,10){5}{\circle*{2}}
\matrixput(1,0)(10,0){6}(0,10){5}{\line(1,0){8}}
\matrixput(7,0)(10,0){6}(0,10){5}{{\vector(1,0){0}}}
\matrixput(4,-4)(10,0){1}(0,10){1}{{\scriptsize $a^{(1)}$}}
\matrixput(0,1)(10,0){7}(0,10){4}{\line(0,1){8}}
\matrixput(0,7)(10,0){7}(0,10){4}{{\vector(0,1){0}}}
\matrixput(-6,4)(10,0){1}(0,10){1}{{\scriptsize $a^{(2)}$}}
\matrixput(-1,0)(10,0){1}(0,10){5}{\line(-1,0){4}}
\matrixput(61,0)(10,0){1}(0,10){5}{\line(1,0){4}}
\matrixput(0,41)(10,0){7}(0,10){1}{\line(0,4){4}}
\matrixput(0,-1)(10,0){7}(0,10){1}{\line(0,-4){4}}
\end{picture}
&
\qquad \qquad 
&
\setlength{\unitlength}{0.8mm}
\begin{picture}(62,60)(-1,-1)
\thicklines
\multiput(0,0)(10,0){7}{\circle*{2}}
\multiput(10,10)(10,0){6}{\circle*{2}}
\multiput(20,20)(10,0){5}{\circle*{2}}
\multiput(30,30)(10,0){4}{\circle*{2}}
\multiput(40,40)(10,0){3}{\circle*{2}}
\multiput(50,50)(10,0){2}{\circle*{2}}
\multiput(1,0)(10,0){6}{\line(1,0){8}}
\multiput(11,10)(10,0){5}{\line(1,0){8}}
\multiput(21,20)(10,0){4}{\line(1,0){8}}
\multiput(31,30)(10,0){3}{\line(1,0){8}}
\multiput(41,40)(10,0){2}{\line(1,0){8}}
\multiput(51,50)(10,0){1}{\line(1,0){8}}
\multiput(7,0)(10,0){6}{{\vector(1,0){0}}}
\multiput(17,10)(10,0){5}{{\vector(1,0){0}}}
\multiput(27,20)(10,0){4}{{\vector(1,0){0}}}
\multiput(37,30)(10,0){3}{{\vector(1,0){0}}}
\multiput(47,40)(10,0){2}{{\vector(1,0){0}}}
\multiput(57,50)(10,0){1}{{\vector(1,0){0}}}
\multiput(14,-3.9)(10,0){1}{{\scriptsize $a^{(2)}$}}
\matrixput(10,1)(10,0){1}(0,10){1}{\line(0,1){8}}
\matrixput(20,1)(10,0){1}(0,10){2}{\line(0,1){8}}
\matrixput(30,1)(10,0){1}(0,10){3}{\line(0,1){8}}
\matrixput(40,1)(10,0){1}(0,10){4}{\line(0,1){8}}
\matrixput(50,1)(10,0){1}(0,10){5}{\line(0,1){8}}
\matrixput(60,1)(10,0){1}(0,10){5}{\line(0,1){8}}
\matrixput(10,7)(10,0){1}(0,10){1}{{\vector(0,1){0}}}
\matrixput(20,7)(10,0){1}(0,10){2}{{\vector(0,1){0}}}
\matrixput(30,7)(10,0){1}(0,10){3}{{\vector(0,1){0}}}
\matrixput(40,7)(10,0){1}(0,10){4}{{\vector(0,1){0}}}
\matrixput(50,7)(10,0){1}(0,10){5}{{\vector(0,1){0}}}
\matrixput(60,7)(10,0){1}(0,10){5}{{\vector(0,1){0}}}
\multiput(11.6,4)(10,0){1}{{\scriptsize $a^{(1)}$}}
\dottedline{2}(-5,-5)(55,55)
\matrixput(61,0)(10,0){1}(0,10){6}{\line(1,0){4}}
\matrixput(0,-1)(10,0){7}(0,10){1}{\line(0,-1){4}}
\matrixput(60,51)(10,0){1}(0,10){1}{\line(0,1){4}}
\end{picture}
\\ & & \\ 
\textnormal{Fig. 1:~} G/P=\PP^1\times\PP^1.
& & 
\textnormal{Fig. 2:~} G/P=\PP^2.
\end{array}
$$

\bfsubsubsection{Quiver for the series $A, D, E$ and complex projective plane}
\label{subsub:hb-P2}

\begin{proposition}\label{prop:dim-arrows}
Let $\Delta(\glu)$ be the set of weights of $\glu$ with respect to
$\glh$, $Q$ be the quiver associated to $P$, and 
$\lambda,\mu\in Q_0$. If all the simple components of $G$ are in one
of the series $A, D, E$, then the number of arrows from $\lambda$ to
$\mu$ is 
$$
n_{\mu\lambda}= \left\{\begin{array}{ll}        
                1 & \quad {\rm if~}\mu-\lambda\in\Delta(\glu);\\
                0 & \quad {\rm otherwise.}
                \end{array}\right.
$$
\end{proposition}

\proof
The proof is as in \cite[Proposition 2]{BK}, where an isotopical
decomposition of $\glu\otimes M_{\mu}$, as an $L$-module, is
obtained by applying Weyl's character formula, and the fact that, for
the series $A, D, E$, the off-diagonal elements of the Cartan matrix
are equal to $0$ or $-1$. We only have to adapt the proof to get an
isotopical decomposition of $\glu^*\otimes M_{\mu}$ instead.
\qed

As an example, in the rest of \S \ref{subsub:hb-P2}, we apply
Proposition \ref{prop:dim-arrows} to describe the quiver corresponding
to the {\bf complex projective plane} $\PP^2$, slightly modifying the
treatment in \cite{BK}. Let $U=U'\oplus U''$ be a $3$-dimensional
vector space, where $U'$ is $2$-dimensional and $U''$ is 
$1$-dimensional. Thus, $\PP^2=\PP(U)=G/P$, where $G$ is $SL(3,\CC)=SL(U)$
and $P$ is its parabolic subgroup 
$$
P=\begin{array}{cc}
\begin{pmatrix}
        * & * & 0     \\
        * & * & 0     \\
        * & * & *     \\
\end{pmatrix} &
\end{array}
$$
of block lower triangular matrices, i.e. automorphisms preserving
$U'\subset U$. A Levi subgroup $L$ is the subgroup of determinant one
matrices in $GL(U')\times GL(U'')$. 
Let $\{ u_1, u_2, u_3\}$ be a basis of $U$, with $u_1,u_2\in U'$, and
$u_{3}\in U''$. Let $\{ E_{i,j}|1\leq i,j\leq 
3\}$ be the basis of $\End U$ given by $E_{i,j}\cdot u_k=\delta_{jk}u_i$.
Let $\glh_\star\subset\End U$ be the subspace with basis vectors $H_{\star
i}:=E_{i,i}$, and let $L_{\star 1},L_{\star 2},L_{\star 3}$ be its dual
basis. The subspace $\glh\subset\glh_\star$ of traceless vectors $c_1
H_{\star 1}+c_2 H_{\star 2}+c_3 H_{\star 3}$, i.e. with
$c_1+c_2+c_3=0$, is a Cartan subalgebra of $\glg$, with dual 
$\glh^*= \glh^*_\star/\gld$, where $\gld=\CC\cdot (L_{\star
1}+L_{\star 2} +L_{\star 3})$. Let $L_i$ be the image of $L_{\star i}$
under the projection $\glh_\star^*\to\glh^*$. Then 
$$
\glp=\glu\ltimes\gll, 
\quad
{\rm where}~~
\glu=\CC E_{1}\oplus\CC E_{2},
~~
\gll=\glh\oplus\CC E_{1,2}\oplus\CC E_{2,1}, 
\quad 
{\rm with~} E_1:=E_{3,1}, E_2:=E_{3,2}.
$$
The sets of roots of $\glg$, $\gll$ and $\glu$ are 
$\Delta=\{\alpha_{i,j}|0\leq i\neq j\leq 3\}$,
$\Delta(\gll)=\{\alpha_{1,2},\alpha_{2,1}\}$,
$\Delta(\glu)=\{\gamma_1,\gamma_2\}$, resp., with 
$\alpha_{ij}=L_i-L_j$, $\gamma_k=\alpha_{3,k}$. 
We choose $\alpha_i:=\alpha_{i,i+1}$, for $i=1,2$, as simple roots
of $(\glg,\glh)$, so $\Sigma=\{\alpha_{2}\}$ (cf. \S \ref{subsub:notation-S-Sigma}).   
The fundamental weights of $\glg$ (resp. $\gll$) are
$\lambda_{\alpha_1}=L_1$, $\lambda_{\alpha_2}=L_1+L_2$, 
(resp. $\lambda_{\alpha_1}=L_1$). Thus, $\Lambda_P^+$ is the set of weights 
$\lambda\in\glh^*$ which are integral for $\glg$ and dominant for
$\gll$, i.e. $\lambda=l_1\lambda_{\alpha_1}+l_2\lambda_{\alpha_2}$, 
with $l_1,l_2\in\ZZ$, and $l_1\geq 0$. Expressing them in the basis
$\{ L_1,L_2\}$, 
$
\lambda=\lambda_1 L_1+\lambda_2 L_2
$ 
with $\lambda_1=l_1+l_2$, $\lambda_2=l_2$, so the condition $l_1\geq
0$ is equivalent to $\lambda_{1}\geq\lambda_{2}$. 
To get a nice picture of this quiver, we use the vectors 
$\epsilon_1=-\frac{1}{3}L_1-\frac{2}{3}L_2, \epsilon_2
=-\frac{2}{3}L_1-\frac{1}{3}L_2$ as the standard basis of $\glh^*$, so
$L_1=(1,-2), L_2=(-2,1)$. Let $\Lambda$ be the lattice generated by
$L_1,L_2$. Any $\lambda=\lambda_1 L_1+ \lambda_2 L_2$ can be written
as $\lambda=x_1\epsilon_1+x_2\epsilon_2$, with
$x_1=\lambda_1-2\lambda_2$, $x_2=-2\lambda_1+\lambda_2$. The condition
$\lambda_1\geq\lambda_2$ is equivalent to $x_1\geq x_2$, for
$x_1-x_2=3(\lambda_1-\lambda_2)$. So the vertex set is 
$$
Q_0=\{(x_1,x_2)\in\Lambda | x_1\geq x_2\}. 
$$
To get the arrows, we see that $L_{3}=-L_1-L_2$, so 
$\gamma_1=L_{3}-L_1=-2L_1-L_2=3\epsilon_2$,
$\gamma_2=L_{3}-L_2=-L_1-2L_2=3\epsilon_1$.
Applying Proposition \ref{prop:dim-arrows}, the arrows are
$$
a^{(1)}_x:x=(x_1,x_2)\to (x_1,x_2+3),\quad a^{(2)}_x:x=(x_1,x_2)\to (x_1+3,x_2),
$$
where $x\in Q_0$. Given an arrow $a^{(i)}_x:(x_1,x_2)\ra (y_1,y_2)$,
with $i=1,2$, we see that $y_1+y_2=x_1+x_2+3$;  thus, the quiver
decomposes into three connected components, say 
$Q^{(0)},Q^{(1)}, Q^{(2)}$, with  $Q_0^{(h)}=\{(x_1,x_2)\in Q_0|
x_1+x_2\equiv -h\mod{3}\}$. The picture of any connected component is
given in Fig. 2 above. 
Although we shall not obtain the relations associated to
$P$, we can easily prove that they are quadratic. 
In fact, $\glu$ is an abelian algebra, for $[E_i,E_j]=0$, so the
linear maps $\psi_{\mu\lambda}$ of \S \ref{sub:linear-maps} are zero. 
Now, when the relations are quadratic, our definition of the quiver
and the relations associated to $P$ coincide with those given by
Hille, who has worked out explicitly the case $n=2$
(cf. \cite{Hl1,Hl3}), showing that the relations are
$r_x=a^{(2)}_{x-3\epsilon_1} a^{(1)}_x-a^{(1)}_{x-3\epsilon_2}
a^{(2)}_x$, for $x\in Q_0$, i.e. commutative diagrams. 

\section{Equivariant bundles, equivariant sheaves  and quivers}
\label{sec:equivb-quivers}

The object of this section is to generalize the results of \S
\ref{sec:homvb-quiver} to equivariant vector bundles and sheaves on
$X\times G/P$. To do this, in \S\S \ref{sec:quiver-bundles} and
\ref{sec:equiv-fil}, we define the categories of holomorphic quiver
bundles and $G$-equivariant holomorphic filtrations, as well as the
corresponding categories of sheaves. 
The definition of a quiver bundle applies to any quiver,
while the definition of a $G$-equivariant holomorphic filtration
applies to any complex $G$-manifold. In \S \ref{sec:corresp-cat}, we
extend the results of \S\S \ref{sec:P-modules-quiver-rep} and \ref{sec:flags} to
equivariant vector bundles and sheaves on $X\times G/P$. Thus, we prove that there
is an equivalence between the category of holomorphic equivariant
vector bundles (resp. coherent equivariant sheaves) on $X\times G/P$
and the category of holomorphic quiver bundles (resp. quiver sheaves)
on $X$. We also show that such an equivariant holomorphic vector
bundle or sheaf admits a natural equivariant filtration. 

\subsection{Quiver bundles and quiver sheaves}
\label{sec:quiver-bundles}

The notion of quiver bundle generalises previous concepts of vector
bundles with additional structure. In this subsection,
$Q$ is a (locally finite) quiver.

\begin{definition}
A {\em $Q$-sheaf} $\bcR=(\bcE,\bphi)$ on $X$ is given by a
collection $\bcE$ of coherent sheaves $\cE_v$, for each vertex $v\in Q_0$,
together with a collection $\bphi$ of morphisms
$\phi_a:\cE_{ta}\ra\cE_{ha}$, for each arrow $a\in Q_1$, such that $\cE_v=0$
for all but finitely many $v\in Q_0$. 
\end{definition}

Given a $Q$-sheaf $\bcR=(\bcE ,\bphi)$ on $X$, every (non-trivial)
path $p=a_0\cdots a_m$ in $Q$ induces a morphism of sheaves
$
\phi(p):=\phi_{a_0}\circ\cdots\circ\phi_{a_m}:\, \cE_{tp}\to\cE_{hp}.
$
The trivial path $e_v$ at a vertex $v$ induces $\phi(e_v)=\id :\cE_v\to\cE_v$. 
A $Q$-sheaf $\bcR=(\bcE ,\bphi)$ satisfies a relation
$r=\sum_i c_i p_i$ if $\sum_i c_i \phi (p_i)=0$. 
Let $\cK$ be a set of relations of $Q$. 
A $Q$-sheaf with relations $\cK$, or a $(Q,\cK)$-{\em sheaf},
is a $Q$-sheaf satisfying the relations in $\cK$. 
A {\em holomorphic $Q$-bundle} is a $Q$-sheaf $\bcR=(\bcE,\bphi)$
such that all the sheaves $\cE_v$ are holomorphic vector bundles,
(i.e. locally free sheaves). 
A {\em morphism} $\bmf:\bcR\ra\bcS$ between two
$Q$-sheaves $\bcR=(\bcE,\bphi)$, $\bcS=(\bcF,\bpsi)$, is given by
morphisms $f_v:\cE_v\to\cF_v$, for each $v\in Q_0$, such that $\psi_a\circ
f_{ta}=f_{ha}\circ\phi_a$, for each $a\in Q_1$. 
It is immediate that $(Q,\cK)$-sheaves form an abelian
category. Important concepts in relation to (semi)stability
(cf. \S \ref{subsub:Q-stability}) are the notions of $Q$-subsheaves
and quotient $Q$-sheaves, as well as indecomposable and simple
$Q$-sheaves, which  are defined as in any abelian category. 

\subsection{Equivariant filtrations}
\label{sec:equiv-fil}

Let $M$ be a complex $G$-manifold. To simplify the notation,
throughout this paper a {\em $G$-equivariant coherent sheaf} on $M$
will mean a coherent sheaf $\cF$ on $M$ together with a {\em
holomorphic} $G$-equivariant action on $\cF$
(cf. \S \ref{subsub:Frechet}, or e.g. \cite{Ak} for more details). 

\begin{definition}
A {\em $G$-equivariant sheaf filtration} on $M$ is a finite sequence
of $G$-invariant coherent subsheaves of a $G$-equivariant coherent
sheaf $\cF$ on $M$, 
\begin{equation}\label{eq:hfil}
\bcF:\,\holfil .
\end{equation}

We say that $\bcF$ is a {\em holomorphic filtration} if the sheaves
$\cF_i$ are locally free.
\end{definition}

When $G$ is the trivial group, a $G$-equivariant coherent sheaf 
will be referred to as a {\em sheaf filtration}. 

The $G$-equivariant sheaf filtrations on $M$ form an abelian
category, whose morphisms are defined as follows. Let $\bcF$, given by
\eqref{eq:hfil}, and $\bcF'$, given by
\begin{equation}\label{eq:hsfil}
\bcF':\, \holsfil ,
\end{equation}
be $G$-equivariant sheaf filtrations. A $G$-equivariant morphism 
from $\bcF'$ to $\bcF$ is a morphism of coherent $G$-equivariant sheaves
$f:\cF'\to\cF$, such that $f(\cF'_i)=\cF_i\cap\Im(f)$ for $0\leq i\leq m$. 
In particular, a subobject of a $G$-equivariant sheaf
filtration \eqref{eq:hfil} is a sheaf filtrations \eqref{eq:hsfil}, where
$\cF'$ is a $G$-invariant coherent subsheaf of $\cF$, and 
$\cF'_i=\cF_i\cap\cF'$, for $0\leq i\leq m$.

\subsection{Correspondence between equivariant sheaves, equivariant
filtrations, and quiver sheaves}
\label{sec:corresp-cat}

We now prove the main results of this section.
 Throughout \S \ref{sec:corresp-cat}, $(Q,\cK)$ is the quiver with relations
associated to $P$. 

\begin{theorem}\label{thm:equivalence-categories} 
There is an equivalence of categories 
$$
\left\{\begin{array}{c}
{\rm coherent~} G {\rm -equivariant} \\
{\rm sheaves~on~} X\times G/P
\end{array}\right\}
\longleftrightarrow
\left\{\begin{array}{c}
(Q,\cK){\rm -sheaves~on~} X
\end{array}\right\}.
$$
The holomorphic $G$-equivariant vector bundles on $X\times G/P$ and the
holomorphic $(Q,\cK)$-bundles on $X$ are in correspondence by this
equivalence. 
\end{theorem}

\begin{proposition} \label{prop:isotopical-equiv-fil-1}
Let us fix a total order in the set $Q_0$, as in \S
\ref{subsub:total-order}. Any coherent $G$-equivariant sheaf 
$\cF$ on $X\times G/P$ admits a $G$-equivariant sheaf filtration
\begin{equation}\begin{gathered}\label{eq:equi-hol-fil}
\bcF:\, \holfil ,\\
\cF_s/\cF_{s-1}\cong p^*\cE_{\lambda_s}\otimes q^*\cO_{\lambda_s},
\quad 0\leq s\leq m,
\end{gathered}\end{equation}
where $\{\lambda_0,\lambda_1,\ldots ,\lambda_m\}$ is a
finite subset of $Q_0$, listed in ascending order as
$\lambda_0<\lambda_1<\cdots<\lambda_m$, and $\cE_0,\ldots,\cE_m$ are
non-zero coherent sheaves on $X$, with trivial $G$-action. 
If $\cF$ is a holomorphic $G$-equivariant vector bundle, 
then $\cE_0,\ldots,\cE_m$ are holomorphic vector bundles. 
\end{proposition}

The proof of Theorem \ref{thm:equivalence-categories} will be
given in two steps, by (i) reduction of a coherent $G$-equivariant
sheaf $\cF$ on $X\times G/P$ to the slice 
$i:X\cong X\times P/P\hra X\times G/P$, obtaining coherent
$P$-equivariant sheaf $\cE=i^*\cF$ on $X$, following \S
\ref{subsub:ind-red}; and (ii) given $\cE$, the construction of a
$(Q,\cK)$-sheaf $\bcR=(\bcE,\bphi)$ on $X$, following the proof of
Theorem \ref{thm:equivalence-P-modules-quivers}. The first step is
Lemma \ref{lemma:induction-restriction}. A preliminary result needed
in the second step is Lemma \ref{lemma:sheaf-decomposition}, where the
$L$-equivariant constant sheaf on $X$ associated to the $L$-module
$M_\lambda$ is also denoted $M_\lambda$, for each $\lambda\in Q_0$.
In Lemma \ref{lemma:induction-restriction} (resp. Lemma
\ref{lemma:sheaf-decomposition}), the $P$-action (resp. $L$-action) 
on $X$ is trivial.

\begin{lemma}\label{lemma:induction-restriction}
There is an equivalence of categories
$$
\left\{\begin{array}{c}
{\rm coherent~} G {\rm -equivariant} \\
{\rm sheaves~on~}  X\times G/P
\end{array}\right\}
\longleftrightarrow
\left\{\begin{array}{c}
{\rm coherent~} P {\rm -equivariant} \\
{\rm sheaves~on~} X
\end{array}\right\}.
$$
\end{lemma}

\begin{lemma}\label{lemma:sheaf-decomposition}
A coherent $L$-equivariant sheaf $\cE$ on $X$ admits an isotopical
decomposition 
\begin{equation}\label{eq:sheaf-decomposition}
\cE\cong\bigoplus_{\lambda\in Q_0'}\cE_\lambda\otimes M_\lambda
\end{equation}
where $Q_0'\subset Q_0$ is finite, and $\cE_\lambda$ is a 
coherent sheaf with trivial $L$-action, for each $\lambda\in
Q'_0$. 
\end{lemma}

Since we are dealing with equivariant sheaves, an ingredient in
the proof of the previous lemmas, which did not appear for 
homogeneous bundles on $G/P$, will be standard
techniques of representation theory of Lie groups on the 
Fr\'echet spaces of sections of equivariant coherent sheaves. 
The necessary preliminaries are explained in \S
\ref{subsub:Frechet}-\ref{subsub:sheaves-invsec}. 
The proofs of Lemmas \ref{lemma:induction-restriction} and
\ref{lemma:sheaf-decomposition} are in \S\S
\ref{subsub:proof-lemma1}, \ref{subsub:proof-lemma2}, and the
proof of Theorem \ref{thm:equivalence-categories} is completed
in \S \ref{subsub:proof-thm-equivalence}. 
The proof of Proposition \ref{prop:isotopical-equiv-fil-1} is the
content of \S \ref{subsub:proof-prop-equiv-fil}.

\bfsubsubsection{Equivariant holomorphic vector bundles} 
\label{subsub:equiv-hol-vb-arguments}

To understand the proof of Theorem \ref{thm:equivalence-categories},
let us first prove the equivalence between $G$-equivariant
holomorphic vector bundles on $X\times G/P$ and holomorphic
$(Q,\cK)$-bundles on $X$. 
In \S \ref{subsub:ind-red} we defined an
equivalence functor between the category of $G$-equivariant 
holomorphic vector bundles on $X\times G/P$ and the category of
$P$-equivariant holomorphic vector bundles on $X$. We now define an
equivalence functor between the category of $P$-equivariant
holomorphic vector bundles on $X$ and the category of holomorphic
$(Q,\cK)$-bundles on $X$, generalising 
Theorem \ref{thm:equivalence-P-modules-quivers}. Let $\cE$ be a
$P$-equivariant holomorphic vector bundle on $X$. The equivalence
functor of Theorem \ref{thm:equivalence-P-modules-quivers} associates
a $(Q,\cK)$-module $\bcR_x=(\bcE_x,\bphi_x)$, to the holomorphic
representation $\cE_x$ of $P$, for each $x\in X$. The point is to
show that these $(Q,\cK)$-modules `vary holomorphically in $X$',
i.e. that they are the fibers of a holomorphic $(Q,\cK)$-bundle on $X$.
To do this, one can use standard techniques in representation
theory. First, there is a (unique) holomorphic projection 
operator $\Pi:\cE\to\cE$ onto the
$L$-invariant part of $\cE$, and its image $\Pi\cE$ and
kernel $(\id-\Pi)\cE$ have induced structures of holomorphic vector
subbundles of $\cE$. A proof is as follows. Let $\Pi$ be
defined by $\Pi(e)=\int_J g\cdot e~dg$, where $dg$ is the Haar measure
of $J$. The map $\Pi$ is obviously a $J$-invariant smooth
projection operator onto the $J$-invariant part. One can prove
that, since $L$ is the universal complexification of $J$ (for $L$ is
reductive), $\Pi$ is actually $L$-invariant and, moreover,
holomorphic. To prove that $\Pi\cE$ and $(\id-\Pi)\cE$ 
have induced structures of holomorphic vector subbundles of $\cE$, one
follows e.g. as in \cite[Lemma (1.4)]{AB}, where a similar result is
proved for a smooth projection operator on a smooth vector bundle
---one simply changes the word `smooth' by `holomorphic' in that
proof. The image of the projection operator $\Pi$ picks out the
$L$-invariant part of $\cE$, so the isotopical decomposition of $\cE$
is obtained as in the familiar case of finite dimensional
representations: 
$$
\cE\cong\bigoplus_{\lambda\in Q'_0} \cE_\lambda\otimes M_{\lambda},
\quad \cE_\lambda=(\Hom(M_{\lambda},\cE))^L,
$$
where $M_\lambda$ is the $L$-equivariant vector bundle $X\times
M_\lambda$, for $\lambda\in Q_0$, $\Hom(M_{\lambda},\cE)$ is the
holomorphic bundle of endomorphisms, $(-)^L:=\Pi(-)$, $\cE_\lambda$
is a holomorphic vector bundle with trivial $L$-action, and
$Q'_0\subset Q_0$ is the subset of weights $\lambda$ with $\cE_\lambda\neq
0$. This isotopical decomposition defines the action of $L$ on
$\cE$. The extension of the $L$-action to a $P$-action 
is defined by means of holomorphic morphisms $\phi_a:\cE_{ta}
\to\cE_{ha}$ satisfying the relations in $\cK$. The proof follows 
as in Theorem \ref{thm:equivalence-P-modules-quivers}. 
The holomorphic vector bundles $\cE_\lambda$ and the holomorphic
morphisms $\phi_a$ define the holomorphic $(Q,\cK)$-bundle
$\bcR=(\bcE,\bphi)$ associated to $\cE$.
\qed

It is worth mentioning that the equivalence proved in \S
\ref{subsub:equiv-hol-vb-arguments} can be stated, and proved, in
terms of $\dbar$-operators on smooth equivariant vector bundles, as
we shall do in \S \ref{sec:hol-str}. 

\bfsubsubsection{Holomorphic actions on coherent sheaves}
\label{subsub:Frechet}

Throughout \S \ref{subsub:Frechet}, $\Gamma$ is a complex Lie group (later
on, it will be either $G$, $P$, or $L$), and $M$ is a complex $\Gamma$-manifold. 
As mentioned in \S \ref{sec:equiv-fil}, to simplify the notation,
throughout this paper a {\em $\Gamma$-equivariant coherent sheaf} on 
$M$ is a coherent sheaf $\cF$ on $M$ together with a {\em holomorphic}
$\Gamma$-equivariant action on $\cF$, following the definitions given
e.g. in \cite{Ak}. This means that $\cF$ is a coherent sheaf on $M$,
hence a covering space with projection map $\pi:\cF\to M$, together
with a $\Gamma$-action on $\cF$, so each $g\in \Gamma$ acts on the stalk
$\cF_y$, for $y\in M$, by an isomorphism $\rho_g:\cF_y\to\cF_{g\cdot
y}$ of $\cO_{g\cdot y}$-modules, such that the $\Gamma$-action commutes with
$\pi$, and is holomorphic in the following sense: 
Since $\cF$ is coherent, the spaces of sections $\cF(B)$, for 
$B\subset M$ open, have canonical Fr\'echet topologies 
(cf. e.g. \cite[\S V.6]{GrR}); the $\Gamma$-action on $\cF$ is
holomorphic if for all open subsets $B,B'\subset X$ and $W\subset \Gamma$
with $W\cdot B\subset B'$, and all sections $s\in\cF(B')$, the map  
$W\to\cF(B)$, $g\mapsto g^{-1}\cdot(s|_{g\cdot B})$,
is holomorphic with respect to the Fr\'echet topology on
$\cF(B)$. 

\bfsubsubsection{Sheaves of invariant sections} 
\label{subsub:sheaves-invsec}

Let $\Gamma$ and $M$ be as in \S \ref{subsub:Frechet}. Let us assume that
the $\Gamma$-action on $M$ is such that 
$M/\Gamma$ is a complex manifold. The structure sheaf of $M/\Gamma$ is given by
$\cO_{M/\Gamma}(B)=\cO_M(\pi^{-1}(B))^{\Gamma}$, for $B\subset M/\Gamma$ open,
where $\pi:M\to M/\Gamma$ is projection. We now define a functor 
from the category of coherent $\Gamma$-equivariant sheaves 
on $M$ to that of (not necessarily coherent) sheaves on $M/\Gamma$, by
`taking invariant sections'. Given a coherent $\Gamma$-equivariant sheaf
$\cF$ on $M$, the sheaf $\cF^{\Gamma}$ associates to $B\subset M/\Gamma$
open, the subspace $\cF^{\Gamma}(B)=\cF(\pi^{-1}(B))^{\Gamma}$ of
$\Gamma$-invariant sections; the restriction maps of $\cF^{\Gamma}$ are those
of $\cF$ restricted to the spaces of invariant sections. This defines a
presheaf, which is in fact a sheaf, called the sheaf of invariant
sections of $\cF$. Given a $\Gamma$-equivariant morphism
$f:\cF_1\to\cF_2$ of coherent $\Gamma$-equivariant sheaves on $M$, the
morphism $f^{\Gamma}:\pi^{\Gamma}_*\cF_1\to\pi^{\Gamma}_*\cF_2$ is defined by
$f^{\Gamma}(B)=f(\pi^{-1}(B))|{ \cE(\pi^{-1}(B))^{\Gamma}}$ for any open set
$B\subset M/\Gamma$. 

\bfsubsubsection{Proof of Lemma \ref{lemma:induction-restriction}}
\label{subsub:proof-lemma1}

Let $\pi_X:X\times G\to X$, $\pi_G:X\times G\to G$, $\pi:X\times
G\to X\times G/P$ be projections. 
The $P$-action on $G$, $r:P\times G\ra G$, $(p,g)\mapsto
r_{p}g=gp^{-1}$, induces a $P$-action on $\cO_G$, 
$\cO_G(B)\to\cO_G(r_p(B))$, $f\mapsto f\circ r_{g^{-1}}$, for $p\in
P$ and $B\subset G$ open, which is holomorphic (cf. e.g. \cite[\S
4.1]{Ak}). This action induces another $P$-action on $\cO_{X\times
G}= \pi^*_X\cO_X\otimes\pi^*_G\cO_G$, obtained from the trivial
$P$-action on the first factor, and the induced $P$-action on the
second factor; thus, $\cO_{X\times G/P}=\cO^P_{X\times G}$. 
Let $\cE$ be a coherent $P$-equivariant sheaf on $X$. 
The coherent sheaf $\cH=\pi_X^*\cE \otimes\pi_G^*\cO_G$ on
$X\times G$ has a holomorphic $P$-equivariant action, obtained from
the induced $P$-actions on the first and the second factor.
We now prove that the sheaf $\cF=\cH^P$ on $X\times G/P$ is
coherent. The set $S(\cF)$ of points $y\in X\times
G/P$ where $\cF$ is not coherent, i.e. where it does not admit a
presentation 
$$
\begin{CD}
        \cO^{r'}_{X\times G/P, y} @>>> \cO^r_{X\times G/P, y} 
        @>>> \cF_{y} @>>> 0, 
\end{CD}
$$
is $G$-invariant, i.e. $G\cdot S(\cF)=S(\cF)$, and $G\cdot y=\{
x\}\times G/P$ for $y=(x,gP)\in X\times G/P$, so $S(\cF)=S(\cE)\times
G/P$, where $S(\cE)$ is the set of points $x\in X$ where $\cE$ is not
coherent. Since $\cE$ is 
coherent, $\cF$ is coherent. Analogously, one shows that $\cF$ is
locally free if $\cE$ is locally free.  

The $G$-action on $G$, $l:G\times G\to G$, $(g,g')\mapsto
l_{g}g'=gg'$, induces a holomorphic $G$-action on $\cO_G$, 
$\cO_G(W)\ra\cO_G(l_{g}(W)),f\mapsto f\circ l_{g^{-1}}$, 
for $g\in G$, $W\subset G$ open, which 
induces a holomorphic $G$-action on $\cO_{X\times G}$
and on the tensor product $\cH=\pi_X^*\cE \otimes\pi_G^*\cO_G$: $G$
acts trivially on the first factor and in the induced way in the
second factor. Thus, for all open subsets $B,B'\subset X\times G$, $W\subset
G$ with $W\cdot B\subset B'$, and $s\in\cH(B')$, the
map $W\to\cH(B)$, $g\mapsto g^{-1}\cdot(s|_{g\cdot B})$, 
is holomorphic. Hence, if $B_1,B_1'\subset X\times G/P$ are
open, $W\cdot B_1\subset B_1'$, and $s\in\cF(B'_1)=\cH(B')^P$,
with $B'=\pi^{-1}(B_1')$, then the map $W\to \cF(B_1)=\cH(B)^P$, $g\mapsto
g^{-1}\cdot(s|_{g\cdot B})$, is holomorphic, where $B=\pi^{-1}(B_1)$. 
Therefore, $\cF$ is a coherent $G$-equivariant sheaf on $X\times G/P$.
\qed

\bfsubsubsection{Proof of Lemma \ref{lemma:sheaf-decomposition}}
\label{subsub:proof-lemma2}

Since $L$ acts trivially on $X$, a holomorphic $L$-equivariant
action on a coherent sheaf $\cE$ on $X$, as defined in \S
\ref{subsub:Frechet}, is simply an $L$-action on 
$\cE$ such that for each $B\subset X$ open and each
$s\in\cE(B)$, the map $L\to\cE(B)$, $g\mapsto g\cdot s$, 
is holomorphic (because the maps $\cE(B)\to\cE(B')$, $s\mapsto
s|_{B'}$, for $B'\subset B\subset X$ open, and $L\to L$, $g\mapsto
g^{-1}$, are holomorphic).
To prove Lemma \ref{subsub:proof-lemma2}, we shall need the following
three lemmas, where it is crucial that $L$ is reductive. The first one
is an equivariant analogue of a well-known property of coherent
sheaves. 

\begin{lemma}\label{lemma:equiv-presentation}
A coherent $L$-equivariant sheaf on $X$ is a sheaf $\cE$ of
$\cO_X$-modules, together with an $L$-action on the space of
sections $\cE(B)$ by automorphisms of $\cO(B)$-modules, 
for each $B\subset X$, with commute with the restriction maps
$\rho_{B'B}:\cE(B)\to\cE(B')$ for $B'\subset B\subset X$ open, and
such that for each $x\in X$, there exists a neighbourhood $B$ of $x$,
finite-dimensional complex representations $V,W$ of $L$ and a
$L$-equivariant exact sequence 
\begin{equation}\label{eq:coh-presentation}
\begin{CD}
        \cO_B\otimes W @>g>> \cO_B\otimes V @>f>> \cE|_B @>>> 0, 
\end{CD} 
\end{equation}
such that, for any $B'\subset B$ open, each $s\in\cE(B')$ is
$L$-finite, i.e. the $\CC$-linear span of the orbit $L\cdot s$ is
finite-dimensional. 
\end{lemma}

\proof
The existence of $B, B', W$ and the $L$-equivariant exact sequence
is a consequence of a result of Roberts \cite[Proposition
2.1]{Ro}. Actually, Roberts states a very close result for complex
$L$-spaces which can be covered by $L$-stable Stein open sets, and
he also proves the existence of a finite set
$s_1,\ldots,s_r\in\cE(B)$ of $L$-finite sections such that the
restrictions $s_1|_{B'},\ldots,s_r|_{B'}$ generate $\cE(B')$ as an 
$\cO_X(B')$-module, for each $B'\subset B$ open. Moreover, 
the representation $V$ in the lemma is a finite dimensional
$L$-invariant complex subspace of $\cE(B)$ generated by
$s_1,\ldots,s_r$ over $\CC$. We now show that for $B'\subset B$ open,
each $s\in\cE(B')$ is $L$-finite. Let $B'\subset B$ open and
$s_i'=s_i|_{B'}$ for $1\leq i\leq r$, so $s'_1,\ldots,s'_r$  
generate $\cE(B')$ over $\cO_X(B')$, and they are also
$L$-finite. Let $B'\subset\cE(B')$ be a finite dimensional
$\CC$-vector subspace containing $s'_1,\ldots,s'_r$. By adding more
generators if necessary, we can assume that the $\CC$-linear span of
$s_1,\ldots,s_r$ is $V'$. Thus, there are functions $f^i_{~j}:L\ra\CC$
such that $g\cdot s_j'=\sum_{i=1}^n f^i_{~j}(g) s_i'$, for each $g\in
L$. To show that each $\xi\in\cE(B')$ is $L$-finite, we first 
expand it in the set of $\cO(B')$-generators $s'_i$, 
$\xi=\sum_{j=1}^n \xi^j s'_j$, with $\xi^j\in\cE(B')$. 
Let $V''\subset\cE(B')$ be the finite-dimensional $\CC$-vector
subspace generated by the subset $\{ \xi^j s'_i| 1\leq i,j\leq n\}$ of
$\cE(B')$. If $g\in G$, then 
$g\cdot \xi=\sum_{j=1}^n \xi^j ~ g\cdot s'_j
=\sum_{i,j=1}^n f^i_{~j}(g) ~ (\xi^j s_i')\in V''$, so $\xi$ is $L$-finite.
\qed 

In the following two Lemmas we shall use the functor which takes a
coherent $L$-equivariant sheaf $\cE$ into the sheaf $\cE^L$ of
invariant sections, as defined in \S \ref{subsub:sheaves-invsec}; we
notice that $\cE$ is a sheaf on $X$, since $X/L=X$.
A projection operator on a coherent $L$-equivariant sheaf $\cE$ on
$X$, i.e. a morphism $\Pi:\cE\to\cE$ with $\Pi^2=\Pi$, is called {\em
$L$-invariant}, if $\Pi(g\cdot s)=\Pi(s)$ for $g\in L, s\in
\cE(B)$, $B\subset X$ open. 
Lemma \ref{lemma:invariant-operator-sheaves} is an immediate
consequence of \cite[Proposition 2.2]{Ro}. 

\begin{lemma}\label{lemma:invariant-operator-sheaves} 
\begin{enumerate}
\item[(a)]
Let $\cE$ be a coherent $L$-equivariant sheaf on $X$. There is a
unique $L$-invariant projection operator $\Pi:\cE\ra\cE$ onto
$\cE^L$. If $B'\subset B$ are two open subsets of $X$ and
$\rho_{B'B}:\cE(B)\ra\cE(B')$ is the restriction map, then
$\Pi_{B'}\circ\rho_{B'B}=\rho_{B'B}\circ\Pi_{B}$.
\item[(b)]
The $L$-invariant projection operator commutes with
$L$-equivariant homomorphisms, i.e. if $f:\cE_1\ra\cE_2$ is an 
$L$-equivariant morphism of coherent $L$-equivariant sheaves on $X$,
then $\Pi \circ f=f^L\circ \Pi$.
\end{enumerate}
\end{lemma}

\begin{lemma}\label{prop:coherent-invariant}
\begin{enumerate}
\item[(a)]
The functor which takes a coherent $L$-equivariant
sheaf $\cE$ into $\cE^L$ is exact.
\item[(b)]
If $\cE$ is a coherent (resp. locally free) $L$-equivariant sheaf on
$X$, then $\cE^L$ is coherent (resp. locally free). 
\end{enumerate}
\end{lemma}

\proof
(a) This functor is obviously left exact. To see that it is right
exact, we have to prove that if $f:\cE_1\ra\cE_2$ is an
$L$-equivariant epimorphism then $f^L:\cE^L_1\ra\cE^L_2$ is
surjective. Given $x\in X$, $s_{2,x}\in\cE^L_{2,x}$, there exists
$s_{1,x}\in\cE_{1,x}$ with $f_x(s_{1,x})=s_{2,x}$. By Lemma
\ref{lemma:invariant-operator-sheaves}, $f^L_x(\Pi_x(s_{1,x}))
=s_{2,x}$, so $f^L$ is surjective.\\ 
(b) 
Let $x\in X$. Assume first that $\cE$ is coherent. By Lemma
\ref{lemma:equiv-presentation}, there exists a neighbourhood $B$ of $x$,
finite-dimensional complex representations $V,W$ of $L$ and a
$L$-equivariant exact sequence \eqref{eq:coh-presentation}.
By part (a), the induced sequence
$$
\begin{CD}
        (\cO_B\otimes W)^L @>g^L>> (\cO_B\otimes V)^L @>f^L>> (\cE|_B)^L @>>> 0
\end{CD} 
$$ 
is exact. Obviously $ (\cO_B\otimes V)^L=\cO_B\otimes V^L,
(\cO_B\otimes W)^L=\cO_B\otimes W^L$ and $(\cE|_B)^L=(\cE^L)|_B$. This
proves (b) for coherent sheaves. For locally free sheaves the argument
is analogous.
\qed

We now prove Lemma \ref{lemma:sheaf-decomposition}.
Since $M_\lambda^{\vee}\otimes\cE$ is a coherent (resp. locally free)
$L$-equivariant sheaf, for $\lambda\in Q_0$, the `multiplicity sheaf'
$\cE_{\lambda}=(M_\lambda^{\vee}\otimes\cE)^L$ is coherent
(resp. locally free). Then
$\cE_{\lambda}(B)=(M_\lambda^{\vee}\otimes\cE(B))^L$, for $B\subset X$
open (for taking invariant sections obviously commutes with
restriction). Let $x\in X$. Let $B\subset X$ open as in Lemma 
\ref{lemma:equiv-presentation}. Since all the
sections of $\cE(B')$ are $L$-finite, for $B'\subset B$ open, and $L$
is reductive, there is an isomorphism 
$$
\cE(B')\cong\bigoplus_{\lambda\in Q_0'}
(M_\lambda^{\vee}\otimes\cE(B'))^L\otimes M_\lambda =
\bigoplus_{\lambda\in Q_0'}
\cE_{\lambda}(B')\otimes M_\lambda ,
$$
where $Q'_0$ is the set of vertices $\lambda\in Q_0$ such that
$\cE_\lambda\neq 0$. Taking direct limits in $B'\ni x$, we get
$\cE_x\cong\oplus_{\lambda\in Q_0'}\cE_{\lambda,x}\otimes M_\lambda $, 
which proves \eqref{eq:sheaf-decomposition}. We now prove that $Q_0'$
is a finite set. 
Since $X$ is compact, it is enough to see that there are isomorphisms 
$\cE|_B\cong\oplus_{\lambda\in Q_0'}\cE_{\lambda}|_B\otimes M_\lambda
$, for the open sets $B\subset X$ satisfying the conditions in Lemma
\ref{lemma:equiv-presentation}, where $Q'_0\subset Q_0$ is finite. Let
$V$ be defined as in Lemma \ref{lemma:equiv-presentation}, so there is
an $L$-equivariant epimorphism $f:\cO_B\otimes V\to\cE|_B$. By Lemma
\ref{prop:coherent-invariant}, tensoring by $M_{\lambda}^{\vee}$ and 
taking the $L$-invariant, we get another epimorphism 
\begin{equation}\label{eq:finite-direct-sum}
\begin{CD}
                \cO_B\otimes V_{\lambda} @>(f\otimes \id)^L>> \cE_{\lambda}|_B @>>> 0.
\end{CD} 
\end{equation}
where $(\cO_B\otimes M^{\vee}_{\lambda}\otimes V)^L=\cO_B\otimes
V_{\lambda}$, with $V_{\lambda }:=(M^{\vee}_{\lambda}\otimes V)^L$, and
$(M^{\vee}_{\lambda}\otimes\cE|_B)^L=\cE_{\lambda}|_B$.
Now $V$ is a finite-dimensional complex representation of $L$, so it
has an isotopical decomposition
$V\cong\oplus_{\lambda\in Q_0'} V_{\lambda}\otimes M_\lambda$,
where $ Q_0'\subset Q_0$ is finite, i.e. 
$V_{\lambda}=0$ for all but finitely many $\lambda\in Q_0$. By
\eqref{eq:finite-direct-sum}, $\cE_{\lambda}=0$ for all
but finitely many $\lambda\in Q_0$, so $Q'_0$ is finite.
\qed

\bfsubsubsection{Proof of Theorem \ref{thm:equivalence-categories}}
\label{subsub:proof-thm-equivalence}

Our proof is similar to  that  of Theorem
\ref{thm:equivalence-P-modules-quivers}, together with 
the previous lemmas. By Lemma 
\ref{lemma:induction-restriction}, to prove the theorem we 
define an equivalence functor from the category of coherent
$P$-equivariant sheaves on $X$, to the category of $(Q,\cK)$-sheaves
on $X$. Let $\cE$ be a coherent $P$-equivariant sheaf on $X$. By
restriction to the Levi subgroup $L\subset G$, we obtain a coherent
$L$-equivariant sheaf on $X$, whose $L$-action is given by an
isotopical decomposition \eqref{eq:sheaf-decomposition}. Obviously, the
action of the unipotent radical $U\subset P$ defines an action
$\tau:\glu\to\End_X(\cE)$ of its Lie algebra, satisfying conditions
similar to (CR1), (CR2) in the proof of Theorem
\ref{thm:equivalence-P-modules-quivers}. Condition (CR1) means that
$\tau$ is a global section of the $L$-invariant part of the coherent
$L$-equivariant sheaf $\cW_1=\glu^*\otimes\cE nd_{\cO_X}(\cE)$,
i.e. $s\in\cW_1^L(X)$. Now, $\cE$ has an isotopical decomposition 
\eqref{eq:sheaf-decomposition}, so 
$$
\cW_1\cong\bigoplus_{\lambda ,\mu ,\nu\in Q_0'} A_{\nu\lambda }
\otimes\cH om_{\cO_X}(\cE_{\lambda },\cE_{\mu }) \otimes\Hom(M_{\nu },M_{\mu }),
$$
hence Schur's lemma implies
$
\cW_1^{L}\cong\oplus_{\lambda ,\mu\in Q_0'}
A_{\mu\lambda }\otimes\cH om_{\cO_X}(\cE_{\lambda },\cE_{\mu }).
$
Thus, $\tau$ is in 
$$
\cW_1^{L}(X)\cong\bigoplus_{\lambda ,\mu\in Q_0'}
A_{\mu\lambda }\otimes\Hom_X(\cE_{\lambda },\cE_{\mu }).
$$
The collection $\bcE$ of coherent sheaves $\cE_\lambda$, together with
the collection $\bphi$ of morphisms $\phi^{(i)}_{\mu\lambda}$, define
a $Q$-sheaf $\bcR=(\bcE,\bphi)$. Condition (CR2) means that $\bcR$
satisfies the relations $\cK$. 
\qed

\bfsubsubsection{Proof of Proposition \ref{prop:isotopical-equiv-fil-1}}
\label{subsub:proof-prop-equiv-fil}

Let $\cF$ be a coherent $G$-equivariant sheaf on $X\times G/P$. The
corresponding coherent $P$-equivariant sheaf $\cE$ on $X$, given 
by Lemma \ref{prop:isotopical-equiv-fil-1}, has an isotopical
decomposition \eqref{eq:sheaf-decomposition}, defining the 
$L$-action on $\cE$. Let $Q_0'=\{\lambda_0,\lambda_1,\ldots
,\lambda_m\}$ be listed in ascending
order as $\lambda_0<\lambda_1<\cdots<\lambda_m$. The
coherent subsheaves of $\cE$ defined by
$$
\cE_{(\leq s)}=\bigoplus_{j=0}^s \cE_j\otimes M_{\lambda_j},
\quad {\rm with~} \cE_s:=\cE_{\lambda_s}, {\rm ~for~} 0\leq s\leq m,
$$ 
are $P$-invariant, so there is a $P$-equivariant filtration, and
$P$-equivariant isomorphisms
$$
0\hra\cE_{\leq 0}\hra\cE_{\leq 1}\hra\cdot\hra\cE_{\leq m}=\cE,
\qquad \cE_{(\leq s)}/_{(\leq s-1)}\cong  \cE_s\otimes M_{\lambda_s},
{\rm ~for~} 1\leq s\leq m 
$$
Applying Lemma \ref{prop:isotopical-equiv-fil-1} to each 
$P$-equivariant sheaf $\cE_{(\leq s)}$, to each map $\cE_{(\leq
s-1)}\hra\cE_{(\leq s)}$, and to each isomorphism $\cE_{(\leq s)}/_{(\leq
s-1)}\cong \cE_s\otimes M_{\lambda_s}$, gives a $G$-equivariant
filtration \eqref{eq:equi-hol-fil}.
\qed

\section{Invariant holomorphic structures and quiver bundles}
\label{sec:hol-str}

In Theorem \ref{thm:equivalence-categories}  we proved an
equivalence between the category of $G$-equivariant holomorphic
vector bundles on $X\times G/P$ and the category of holomorphic
$(Q,\cK)$-bundles on $X$, where $(Q,\cK)$ is the quiver with
relations associated to $P$. 
Our purpose in this section is to study the same equivalence in
terms of {\em invariant $\dbar$-operators} on smooth equivariant
vector bundles on $X\times G/P$. The main result, stated as Proposition
\ref{prop:corresp-inv-hol-str}, will be used in the next section to
obtain the dimensional reduction of the gauge equations and of the
stability criteria for equivariant holomorphic bundles on $X\times G/P$.
Throughout this section, $(Q,\cK)$ is the quiver with relations
associated to $P$.

\subsection{Invariant holomorphic structures on $G$-manifolds} 
\label{sub:inv-dbar-general}

A $K$-equivariant smooth complex vector bundle on a smooth
$K$-manifold $M$ is a smooth complex vector bundle $F$ on $M$ 
together with a smooth lifting of the $K$-action to an action on  $F$. 
Let $F$ be such a smooth $K$-equivariant vector bundle on a compact
$K$-manifold $M$. Since $G$ is the universal complexification of $K$,
the $K$-actions on $M$ and $F$ lift to unique smooth $G$-actions,
and hence $F$ is a smooth $G$-equivariant bundle over the $G$-manifold $M$.
The group $G$ acts naturally on the space $\DDD$ of
$\dbar$-operators on $F$, by $\gamma(\dbar_F)=\gamma\circ
\dbar_F\circ\gamma^{-1}$ for $\gamma\in G$ and $\dbar_F\in\DDD$.  
This action of $G$ leaves invariant the subset $\CCC$ of
$\dbar$-operators $\dbar_F$  with $(\dbar_F)^2=0$, which is in 
bijection with the space of holomorphic structures on $F$. The group
$G$ also acts naturally on the complex gauge group
$\GGG^c=\Omega^0(\Aut(F))$ of $F$, by $\gamma(g)=\gamma\circ
g\circ\gamma^{-1}$ for $\gamma\in G$ and $g\in\GGG^c$. The space of
fixed points $\CCC^G$ is in bijection with the space of holomorphic
structures on $F$ such that the action of $G$ is holomorphic.

\subsection{Preliminaries on smooth equivariant vector bundles on
$X\times G/P$} 
\label{subsub:prelim-equiv-smoothvb}

Since $P\subset G$ is a parabolic subgroup, the natural map $K/J\to
G/P$ is a diffeomorphism. Since $G/P$ is a projective variety, this
map induces the same structure on $K/J$. In particular, $K/J$ is a 
compact \kah\ homogeneous manifold with a symplectic action of $K$.

Let $F$ be a $K$-equivariant smooth vector bundle on $X\times K/J$. 
The $K$-action on $E$ lifts to a unique smooth $G$-action, 
giving  $F$ 
the structure of  a smooth $G$-equivariant bundle over 
$X\times K/J\cong X\times G/P$. 
Now, there is an equivalence, similar to \S \ref{subsub:ind-red},
between smooth $K$-equivariant vector bundles on $X\times K/J$ and
smooth $J$-equivariant vector bundles on $X$: a smooth $J$-equivariant 
bundle $E$ on $X$ induces a smooth $K$-equivariant bundle $F=K\times_J E$. 
Furthermore, any smooth equivariant $J$-action on $E$ 
extends uniquely to an $L$-action on $E$ (since $L$ is
the universal complexification of $J$). Thus, if $\cE$ is a holomorphic
$P$-equivariant vector bundle on $X$ whose underlying smooth
$J$-equivariant vector bundle is $E$, then the induced holomorphic
$G$-equivariant vector bundle $\cF=G\times_P\cE$ on $X\times K/J\cong
X\times G/P$ has underlying smooth $K$-equivariant $F=K\times_J
E$. This means that there is a one-to-one correspondence between (i) 
$G$-invariant holomorphic structures on the smooth $K$-equivariant
vector bundle $F=K\times_J E$ and (ii) $L$-invariant holomorphic
structures on the smooth $J$-equivariant vector bundle $E$ 
together with extensions of the $L$-action on $E$ to a holomorphic 
$P$-action on
$E$. As a preliminary step to describe this correspondence
in terms of $\dbar$-operators and quiver bundles, we introduce some
notation and state the following lemma, whose proof is standard (see
e.g. \cite{Se}).

\bfsubsubsection{Notation}\label{subsub:notation-H-lambda} 
Given an irreducible representation $M_\lambda$ of $J$, corresponding
to an integral dominant weight $\lambda$, $H_{\lambda}:=K\times_J
M_\lambda$ is the induced irreducible smooth homogeneous vector
bundle on $K/J$. It is the smooth $K$-equivariant vector bundle
underlying $\cO_\lambda=G\times_P M_\lambda$ (cf. \S
\ref{subsub:notation-O-lambda}). The $\dbar$-operator corresponding
to the holomorphic structure $\cO_\lambda$ is denoted by
$\dbar_{H_\lambda}$. The maps $p:X\times K/J\to X$ and $q:X\times
K/J\to K/J$ are the canonical projections. 

\begin{lemma}\label{lemma:decomp-ECCVB}
Every smooth $K$-equivariant complex vector bundle $F$ on
$X\times K/J$ can be equivariantly decomposed, uniquely up to
isomorphism, as
\begin{equation}\label{eq:equiv-decomp}
F\cong\bigoplus_{\lambda\in Q_0'} F_\lambda,\quad
F_\lambda:=p^* E_\lambda\otimes q^* H_\lambda,
\end{equation}
for some finite collection $\bmE$ of smooth complex vector bundles
$E_\lambda$ on $X$, with trivial $K$-action, where $Q_0'\subset Q_0$
is the set of vertices with $E_\lambda\neq 0$, which is of course  a finite set.
\end{lemma}

The following lemma, proved in \S \ref{subsub:trivial-but-useful}, is
needed before  introducing some more notation. 

\begin{lemma} \label{lemma:isom-a-eta_a}
There is a natural isomorphism
$\Omega^{0,1}(\Hom(H_{\lambda},H_{\mu}))^K\cong A_{\mu\lambda}$.
\end{lemma}

\bfsubsubsection{Notation}\label{subsub:notation-eta_a}

Let $\{\eta_a|a\in Q_0, ta=\lambda, ha=\mu\}$ be a basis of 
$\Omega^{0,1}(\Hom(H_\lambda,H_\mu))^K$ corresponding to a basis
$\{ a^{(i)}_{\mu\lambda }|i=1,\ldots ,n_{\mu\lambda }\}$ of
$A_{\mu\lambda }$ by Lemma \ref{lemma:isom-a-eta_a}, for each
$\lambda,\mu\in Q_0$. 

\subsection{Spaces of $\dbar$-operators, complex gauge groups, and
quiver bundles} 
\label{sub:dbar-cxgauge-grp-quivers}

Let $F$ be a smooth $K$-equivariant vector bundle on $X\times K/J$,
whith equivariant decomposition \eqref{eq:equiv-decomp}. 
The $K$-action on $F$ lifts to a unique smooth $G$-action, so $F$ is a
smooth $G$-equivariant over $X\times K/J$ (cf. \S
\ref{subsub:prelim-equiv-smoothvb}).
Let $\DDD$ (resp. $\CCC$) be the space of $\dbar$-operators (resp. 
$\dbar$-operators with square zero) on $F$, and let $\GGG^c$ be the
complex gauge group of $F$, with the $G$-actions on these spaces
defined in \S \ref{sub:inv-dbar-general}. Let $Q'=(Q_0', Q_1')$ be  
the full subquiver of $Q$ whose vertices $\lambda$ are defined by the
condition $E_\lambda\neq 0$ (the arrows $a$ are defined by the
conditions $E_{ta}\neq 0$ and $E_{ha}\neq 0$). Let $\DDD_\lambda$
(resp. $\CCC_\lambda$) be the space of $\dbar$-operators
(resp. $\dbar$-operators with square zero) on $E_\lambda$, and let
$\GGG^c_\lambda$ be the complex gauge group of $E_\lambda$, for each
$\lambda\in Q_0'$. The group
$$
\GGG^{\prime c}=
\prod_{\lambda\in Q_0'}\GGG^c_\lambda
$$
acts on the space $\DDD'$ of $\dbar$-operators, and on the
representation space $\RRR(Q',\bmE)$, defined by
$$ 
\DDD'=\prod_{\lambda\in Q_0'}\DDD_\lambda ,\quad 
\RRR(Q',\bmE)=\bigoplus_{a\in Q'_1} \Omega^0(\Hom(E_{ta},E_{ha})).
$$
An element $\bmg\in\GGG^{\prime c}$ is a collection of elements 
$g_\lambda\in\GGG^c_\lambda$, for each $\lambda\in Q'_0$, and an element
$\dbar_{E}\in\DDD'$ (resp. $\bphi\in\RRR(Q',\bmE)$) is a 
collection of $\dbar$-operators $\dbar_{E_\lambda}\in\DDD_\lambda$
(resp. smooth morphisms $\phi_a:E_{ta}\to E_{ha}$), for each
$\lambda\in Q'_0$ (resp. $a\in Q'_1$). 
The $\GGG^{\prime c}$-actions on $\DDD'$ and $\RRR(Q',\bmE)$ are
given by $(\bmg(\dbar_{E}))_\lambda=g_\lambda\circ \dbar_{E_\lambda}\circ  
g_\lambda^{-1}$, and $(\bmg\cdot\bphi)_a=g_{ha} \circ\phi_a\circ
g_{ta}^{-1}$, respectively. 
The induced $\GGG^{\prime c}$-action on the product $\DDD'\times
\RRR(Q',\bmE)$ leaves invariant the subset $\NNN$ of pairs
$(\dbar_{E},\bphi)$ such that $\dbar_{E_\lambda}\in\CCC_{\lambda}$,
for each $\lambda\in Q_0'$, $\phi_a:E_{ta}\to E_{ha}$ is
holomorphic with respect to $\dbar_{ta}$ and $\dbar_{ha}$, for each
$a\in Q_0$, and the holomorphic $Q$-bundle $\bcR=(\bcE,\bphi)$,
defined by these holomorphic structures and morphisms, satisfies the
relations in $\cK$.

\begin{proposition} \label{prop:corresp-inv-hol-str}
\begin{enumerate}
\item[(a)]
There is a one-to-one correspondence between $\DDD^G$ and
$\DDD'\times\RRR(Q',\bmE)$ which, to any $(\dbar_{E},\bphi)\in
\DDD'\times\RRR(Q',\bmE)$, associates, the $\dbar$-operator
$\dbar_F\in\DDD^G$ given by 
\begin{equation}\label{eq:corresp-hol-str}
\dbar_F=\sum_{\lambda\in Q_0'} \dbar_{F_\lambda }\circ\pi_\lambda
+\sum_{a\in Q_1'}\beta_a\circ\pi_{ta}.
\end{equation}
Here $\dbar_{F_\lambda}$ is the $\dbar$-operator of
$F_\lambda$ given by $\dbar_{F_\lambda}=
p^*\dbar_{E_\lambda}\otimes\id + \id\otimes q^*\dbar_{H_{\lambda }}$,
for each $\lambda\in Q'_0$, and $\beta_a:=p^*\phi_a\otimes q^*\eta_a
\in\Omega^{0,1}(\Hom(F_{ta}, F_{ha}))$ for each $a\in
Q_1'$.\\  
\item[(b)]
The previous correspondence restricts to a one-to-one
correspondence between $\CCC^G$ and $\NNN$.\\  
\item[(c)]
There is a one-to-one correspondence between
$({\GGG^c})^G$ and $\GGG^{\prime c}$ which, to any 
$\bmg\in\GGG^{\prime c}$, associates $g=\sum_{\lambda\in
Q_0'}\wt{g}_\lambda\circ \pi_\lambda\in({\GGG^c})^G$, with
$\wt{g}_\lambda=p^*g_\lambda\in\Omega^0 (\Aut(F_\lambda))
\cong\Omega^0(\Aut(p^*E_\lambda))$.\\ 
\item[(d)]
These correspondences are compatible with the actions of the
groups of {\em (c)} on the sets of {\rm (a)} and {\em (b)}, hence there is a one-to-one 
correspondence between $\CCC^{G}/({\GGG^c})^G$ and $\NNN/\GGG^{\prime c}$.
\end{enumerate}
\end{proposition}

\subsection{Preliminaries on smooth homogeneous vector bundles}
To prove the previous proposition, in this subsection we collect several
preliminary results about the homogeneous space $K/J$ and the
homogeneous bundles on $K/J$. We first recall several standard results
(\S\S \ref{subsub:normal-cx-str}-\ref{subsub:tilde},
cf. e.g. \cite{Be}), which we adapt to the notation used throughout
this paper. We then prove other results (\S\S
\ref{subsub:an-invariant-connection}-\ref{subsub:invariant-forms}),
which are elementary but necessary to relate invariant connections on
homogeneous bundles to our definition of the quiver with relations
associated to $P$. 

\bfsubsubsection{The canonical complex structure on $K/J$}
\label{subsub:normal-cx-str}

The {\em canonical complex structure} on $K/J$ is the complex structure
on $K/J$ induced by the complex structure on the projective variety
$G/P$ and by the diffeomorphism $K/J\cong G/P$.
Let $T=H\cap K$ be a maximal torus of $K$, and let $\glt$ be its Lie algebra.
Let $\glr\subset\glk$ be the Lie subalbebra, which is also a
$J$-submodule, given by the isomorphisms of $J$-modules
$$
\glk=\glj\oplus\glr,\quad\glr\cong\glk/\glj.
$$
That is, $\glr$ is the direct sum of the even-dimensional real subspaces of
$\glk$ on which the spectrum for the action of $\glt$ is $\pm\imag\lambda$, 
for $\lambda\in\Delta(\glr)$, where $\Delta(\glr)$ is the set of roots
of $\glg$ with respect to $\glh$ which are not roots of $\gll$. 
The complexification of $\glr$ is $\glr_{\CC}=\bar{\glu}
\oplus\glu$. 
We define a $J$-invariant complex structure of the $J$-module $\glr$ by the
condition that 
$$
\glr^{1,0}=\bar{\glu}, \quad \glr^{0,1}=\glu
$$ 
are the $(1,0)$- and $(0,1)$-subspaces of $\glr_{\CC}$,
respectively. Then there are natural isomorphisms  
\begin{equation}\label{eq:normal-cx-str}
\Lambda^{i,j}T^*(K/J)\cong K\times_J\Lambda^{i,j}\glr.
\end{equation}
To prove this, we first note that the holomorphic and antiholomorphic
cotangent bundles on $G/P$ are isomorphic to $G\times_P\glu$ and
$G\times_P\glu^*$, respectively, as holomorphic $G$-equivariant
vector bundles. This follows from: (i) the holomorphic tangent bundle
on $G/P$ is isomorphic to $G\times_P(\glg/\glp)$ as a holomorphic
$G$-equivariant vector bundle; (ii) the Killing form on $\glg$
induces isomorphisms $\glg/\glp\cong\glu^*\cong \bar{\glu}$ of
representations of $P$.
Second, by the arguments of \S \ref{subsub:prelim-equiv-smoothvb}, the
underlying smooth homogeneous vector bundles of $G\times_P\bar{\glu}$
and $G\times_P\glu$ are $K\times_J\glr^{1,0}$ and
$K\times_J\glr^{0,1}$, respectively. This proves \eqref{eq:normal-cx-str}. 
It is worth mentioning that the canonical complex structure can be defined
directly on $K/J$, without the use of the diffeomorphism $K/J\cong
G/P$ (see e.g. \cite{Be}); however we shall be interested in the
relation between the holomorphic structures on smooth homogeneous
vector bundles on $K/J$ and the representation theory of $P$, so we
shall need the isomorphisms \eqref{eq:normal-cx-str}.

\bfsubsubsection{Spaces of invariant forms on $K/J$}
\label{subsub:trivial-but-useful}

Given a smooth $K$-equivariant vector bundle $F=K\times_J W$ on
$K/J$, induced by a representation $W$ of $J$, there is a
bijection between the space of $K$-invariant sections of $F$ and the
$J$-invariant subspace $W^J$ of $W$. If $W=\wedge^{i,j}\glr^*\otimes V$, 
where $V$ is another complex representation of $J$, then 
\eqref{eq:normal-cx-str} implies $\wedge^{i,j}T^*(K/J)\otimes E\cong
K\times_J(\wedge^{i,j}\glr^*\otimes V)$ as smooth $K$-equivariant vector
bundles, where $E:=K\times_J V$. Therefore, the spaces of
$K$-invariant $E$-valued $r$- and $(i,j)$-forms on $K/J$ are
\begin{equation}\label{eq:trivial-but-useful}
\Omega^{r}(E)^K\cong C^{r}(V), \quad {\rm and} \quad
\Omega^{i,j}(E)^K\cong C^{i,j}(V), 
\end{equation}
respectively, where 
\begin{equation}\label{eq:C-ij(V)}
C^{r}(V):=(\wedge^{r}\glr^*\otimes V)^J, \quad {\rm and} \quad
C^{i,j}(V):=(\wedge^{i,j}\glr^*\otimes V)^J.
\end{equation}
In particular, if $V=\Hom(M_\lambda,M_\mu)$, for
$\lambda,\mu\in Q_0$, and $i=1,2$, we get (see \S \ref{subsub:A-B})
\begin{equation}\label{eq:omega-A-B}
\Omega^{0,1}(\Hom(H_{\lambda},H_{\mu}))^K\cong A_{\mu\lambda},\quad
\Omega^{0,2}(\Hom(H_{\lambda},H_{\mu}))^K\cong B_{\mu\lambda}.
\end{equation}

Let $V$ be a representation of $J$, and let $E=K\times_J V$ be the
induced homogeneous vector bundle on $K/J$. A {\em basic} $V$-valued
$r$-form on the principal $J$-bundle $K\to K/J$ is an element of 
$\Omega^r_K(V):=\Omega^r_K\otimes V$ which is
$J$-invariant and horizontal (cf. e.g. \S 1.1 of \cite{BGV} for
definitions). The space $\Omega^r_b(V)$ of basic $V$-valued
$r$-forms on $K\to K/J$ is isomorphic to the space $\Omega^r(E)$ of
$E$-valued $r$-forms on $K/J$. 
The space $\Omega^r(E)$ is naturally a $K$-module and, since the
principal $J$-bundle $K\to K/J$ is $K$-equivariant (with the
canonical left $K$-action), the space $\Omega^{r}_b(V)$ is a
$K$-module as well. Moreover, the previous isomorphism
$\Omega^{r}_b(V)\cong\Omega^{r}(E)$ is an isomorphism of
$K$-modules. Therefore their $K$-invariant parts are isomorphic, so
\eqref{eq:trivial-but-useful} gives an isomorphism
$\Omega^r_b(V)^K \cong C^r(V)$.  

\bfsubsubsection{Spaces of invariant forms on $K$}\label{subsub:tilde} 

There is a natural isomorphism of Lie algebras between $\glk$ and
the space $\glX(K)^K$ of (left) $K$-invariant vector fields on
$K$. Given $x\in\glk$, let $\wt{x}$ be the corresponding
$K$-invariant vector field on $K$ generated by $x$, i.e. 
$\wt{x}(k)=\frac{d}{dt}(k\exp(tx))\big|_{t=0}$, for each $k\in K$.  
Let $V$ be a representation of $J$. The left action of $J$ on $K$
together with the action of $J$ on $V$ define a structure of
$J$-module on the space $\Omega^r_K(V)=\Omega^r_K\otimes V$ of
$V$-valued $r$-forms on $K$, while the right action of $K$ on itself
induces naturally a structure of $K$-module on
$\Omega^r_K(V)$. Moreover,  the $K$-invariant subspace
$\Omega^r_K(V)^K$ is a $J$-submodule of $\Omega^r_K(V)$.  Given
$a\in\wedge^r\glk^*\otimes V$, let $\wt{a}\in\Omega_K^r(V)^K$ be
defined by 
$$
\wt{a}(\wt{x}_1(k),\ldots , \wt{x}_r(k))=a(x_1,\ldots
,x_r), \quad
{\rm for~} k\in K,~ x_1,\ldots ,x_r\in\glk.
$$ 
This map defines an isomorphism of representations of $J$:
\begin{equation}\label{eq:tilde} 
\wedge^r\glk^*\otimes V\cong\Omega^r_K(V)^K.
\end{equation}

\bfsubsubsection{An invariant connection for smooth homogeneous
bundles} 
\label{subsub:an-invariant-connection}
Let $A':\glk\ra\glj$ be the canonical projection from
$\glk=\glj\oplus\glr$ onto $\glj$. Thus $A'\in\glk^*\otimes\glj$ is
$J$-equivariant, so it induces a $K$-invariant $\glj$-valued
one-form $\wt{A}'\in\Omega^1(\glj)^K$ on $K$, by \eqref{eq:tilde}.

\begin{lemma}\label{lemma:connection-one-form}
\begin{enumerate}
\item[(a)]
The one-form $\wt{A}'\in\Omega^1_K(\glj)^K$ is a $K$-invariant
connection one-form on the smooth principal $J$-bundle $K\to K/J$.\\
\item[(b)]
Its curvature $\wt{F}'\in\Omega_{K/J}^{1,1}(\ad K)$, with $\ad
K:=K\times_J\glj$, is given by means of its isomorphic
element $F'\in C^{1,1}(\glj)$, defined by
\begin{equation}\label{eq:connection-value}
F'(x,x')=-A'([x,x'])\quad {\rm ~for~} x,x'\in\glr.
\end{equation}
\end{enumerate}
\end{lemma}

\proof 
Part (a) is straightforward. 
(b) By \S \ref{subsub:trivial-but-useful}, the curvature is defined by
its isomorphic element $\wt{F}'=d\wt{A}'+\frac{1}{2}[\wt{A},\wt{A}]$, 
in $\Omega^2_b(\glj)^K\subset\Omega^2_K(\glj)^K$, which is given by  
$$
\wt{F}'(\wt{x}(k),\wt{x}'(k))=d\wt{A}'(\wt{x}(k),\wt{x}'(k))+[\wt{A}'(\wt{x}(k)),\wt{A}'(\wt{x}'(k))],
{\rm ~~for~} x,x'\in\glk {\rm ~and~} k\in K 
$$
(so $\wt{x}(k),\wt{x}'(k)\in T_k K$, cf. \S \ref{subsub:tilde}). Then 
$d\wt{A}'(\wt{x},\wt{x}')=\wt{x}(\wt{A}'(\wt{x}'))-\wt{x}'(\wt{A}'(\wt{x}))
-\wt{A}'([\wt{x},\wt{x}'])$. But $\wt{A}'(\wt{x})=A'(x)$ is constant,
so $\wt{x}'(\wt{A}'(\wt{x}))=0$, and analogously $\wt{x}(\wt{A}'(\wt{x}'))=0$. 
On the other hand $[\wt{x},\wt{x}']=[x,x']\wt{~~}$. Therefore 
$$
\wt{F}'(\wt{x},\wt{x}')=[\wt{A}'(\wt{x}),\wt{A}'(\wt{x}')]-\wt{A}'([\wt{x},\wt{x}'])=[A'(x),A'(x')]-A'([x,x']).
$$ 
This shows that $\wt{F}'(\wt{x},\wt{x}')=0$ for $x\in\glj$ or
$x'\in\glj$, i.e. $\wt{F}'$ is horizontal, as we already knew from
part (a), and for $x,x'\in\glr$, it is given by $\wt{F}'(\wt{x},\wt{x}')=
-A'([x,x'])$ which is (\ref{eq:connection-value}). Let us denote by
the same symbol the complexification of $F'$. Given $e,e'\in\glu_\pm$, we have
$[e,e']\in\glu_\pm$, since $\glu_\pm$ are Lie algebras, hence its
projection by $A':\glg=\glk_\CC
=\gll\oplus\glu\oplus\bar{\glu}\ra\gll=\glj_\CC$ is zero, i.e. 
$F'([e,e'])=-A'([e,e'])=0$. Hence $F'\in C^{1,1}(\glj)$. 
\qed

Let $V$ be a $J$-module. 
Since the curvature $\wt{F}'$ is of type $(1,1)$, the induced 
connection on any homogeneous vector bundle $E=K\times_J V$ defines a
holomorphic structure, which is obviously $G$-invariant.
Let us define the linear maps $d:C^r(V)\ra C^{r+1}(V)$ between the
linear spaces $C^r(V):=(\wedge^r\glr^*\otimes V)^J$ by
\begin{equation}\label{eq:computing-inv-cov-derivative}
d a(x_0,\ldots ,x_r)=\sum_{0\leq i<j\leq r}(-1)^{i+j}
a([x_i,x_j],x_0,\ldots ,\widehat{x_i},\ldots ,\widehat{x_j},\ldots
,x_r), \quad {\rm for~~} x_0,\ldots ,x_r\in\glr.
\end{equation}

\begin{lemma}\label{lemma:computing-inv-cov-derivative}
Let $A$ be the $K$-invariant connection on $E=K\times_J V$ induced by the
$K$-invariant connection one-form $\wt{A}'$ on $K\ra K/J$. The
following diagram is commutative: 
$$
\begin{CD}
C^{r}(V)                        @>{d}>>         C^{r}(V)        \\
@V{\cong}VV                     @VV{\cong}V                     \\
\Omega^{r}(E)^K                 @>{d_A}>>       \Omega^{r}(E)^K
\end{CD}
$$
Here $C^r(V)\stackrel{\cong }{\ra }\Omega^r(E)^K$ and
$C^{r+1}(V)\stackrel{\cong }{\ra }\Omega^{r+1}(E)^K$ are the 
isomorphism appearing in \S \ref{subsub:trivial-but-useful}.
\end{lemma}

\proof 
First we assume that $r=0$, so let $v\in C^0(V)= V^J$,
$s\in\Omega^0(E)^K$ and $\wt{v}\in\Omega^0_b(V)$ be
related by the isomorphisms of \S \ref{subsub:trivial-but-useful}. By definition
$dv=0$, so we have to see that $d_{A}s=0$, or equivalently, that
$d_{\wt{A}}\wt{v}=0$, where $d_{\wt{A}}= d+\rho(\wt{A}'):
\Omega^0_b(V)\to\Omega^1_b(V)$. Since $\wt{v}:K\to V$
is a constant map, $d\wt{v}=0$, while
$\rho(\wt{A}')\cdot\wt{v}\in\Omega^1_N(E)$ is given by
$(\rho(\wt{A}')\cdot\wt{v})\cdot e=\rho(\wt{A}'\cdot e) \cdot v\in V$ 
for $e\in TK$, which is zero, since $v$ is $J$-invariant. This proves
$d_{A}s=0$. For $r\geq 1$, let $a\in C^r(V)$, and let
$\eta\in\Omega^r(E)^K$ and $\wt{a}\in\Omega_b^r(V)$ be related
the isomorphism of \S \ref{subsub:trivial-but-useful}. Let $x_0,\ldots
,x_r\in\glr$, and let $\wt{x}_0,\ldots, \wt{x}_r$, be as in \S
\ref{subsub:tilde}. Then 
\begin{multline*}
d_{\wt{A}}\wt{a}(\wt{x}_0,\ldots ,\wt{x}_r)=
\sum_{i=0}^r (-1)^i \iota(\wt{x}_i)d_{\wt{A}}(\wt{a}(\wt{x}_0,\ldots
,\widehat{\wt{x}_i},\ldots ,\wt{x}_r)))\\ +\sum_{0\leq i<j\leq r}
(-1)^{i+j}\wt{a}([\wt{x}_i,\wt{x}_j],\wt{x}_0, \ldots
,\widehat{\wt{x}_i},\ldots ,\widehat{\wt{x}_j},\ldots ,\wt{x}_r)),
\end{multline*}
where $\iota(\wt{x}):\Omega^1\ra\Omega^0$ is contraction with the
vector field $\wt{x}$. But $\wt{a}(\wt{x}_0,\ldots
,\widehat{\wt{x}_i},\ldots ,\wt{x}_r)=a(x_0,\ldots
,\widehat{x_i},\ldots ,x_r)$ is constant, so
$d_{\wt{A}}(\wt{a}(\wt{x}_0,\ldots ,\widehat{\wt{x}_i},\ldots
,\wt{x}_r))$ $=0$ as seen before, while
$[\wt{x}_i,\wt{x}_j]=[x_i,x_j]\wt{~~}$. This proves the assertion.
\qed

\bfsubsubsection{Invariant forms on $K/J$ and the linear maps
$\psi_{\mu\lambda}$, $\psi_{\mu\nu\lambda}$}
\label{subsub:invariant-forms}

Let $\{\eta_{\mu\lambda }^{(i)}|i=1,\ldots ,n_{\mu\lambda }\}$ and
$\{\xi_{\mu\lambda }^{(p)}|p=1,\ldots ,m_{\mu\lambda }\}$, 
for fixed $\lambda ,\mu\in Q_0$, be the bases of
$\Omega^{0,1}(\Hom(H_\lambda,H_\mu))^K$ and 
$\Omega^{0,2}(\Hom(H_\lambda,$ $H_\mu))^K$, that correspond
to the basis $\{ a^{(i)}_{\mu\lambda }|i=1,\ldots ,n_{\mu\lambda }\}$
of $A_{\mu\lambda }$ and to the basis $\{ b^{(p)}_{\mu\lambda
}|p=1,\ldots ,b_{\mu\lambda }\}$ of $B_{\mu\lambda }$ by the
isomorphisms \eqref{eq:omega-A-B}, respectively
(cf. \S \ref{subsub:notation-eta_a}). 
Let $c_{\mu\lambda }^{(i,p)}$ and $c^{(j,i,p)}_{\mu\nu\lambda }$ be
the coefficients defined in \S \ref{subsub:P-relations}.
\begin{lemma}
Let $A'_\lambda$ the unique $K$-invariant connection  of $H_{\lambda
}$, and let $A'_{\mu\lambda}$ be the connection induced by
$A'_\lambda$ and $A'_\mu$ on the vector bundle $\Hom(H_\lambda$,
$H_\mu)$. Then  
\begin{equation}\label{eq:eta-coefficients}
\eta^{(j)}_{\mu\nu }\wedge\eta^{(i)}_{\nu\lambda
}=\sum_{p=1}^{m_{\mu\lambda }} c^{(j,i,p)}_{\mu\nu\lambda
}\xi^{(p)}_{\mu\lambda },\quad \dbar_{A'_{\mu\lambda}}
(\eta^{(i)}_{\mu\lambda })=\sum_{p=1}^{m_{\mu\lambda }} c_{\mu\lambda
}^{(i,p)}\xi^{(p)}_{\mu\lambda }, 
\quad {\rm for~each~~} \lambda,\mu\in Q_0.
\end{equation}
\end{lemma}

\proof
The first equation follows immediately from the first in 
\eqref{eq:coefficientes-relations}. The second equation is obtained
from the second equation in \eqref{eq:coefficientes-relations}. 
Indeed,  $d_{A'_{\mu\lambda}} \eta^{(i)}_{\mu\lambda }\in \Omega^2_N(\Hom(H_\lambda,
H_\mu)^K$ corresponds, by the isomorphism \eqref{eq:omega-A-B}, to
$da^{(i)}_{\mu\lambda }\in B_{\mu\lambda}$ as given by
(\ref{eq:computing-inv-cov-derivative}). To evaluate its $(0,2)$-part
$\dbar a^{(i)}_{\mu\lambda }:=(da^{(i)}_{\mu\lambda })^{0,2}
\in\wedge^{0,2}\glr\otimes\Hom(M_\lambda ,M_\mu)$, let
$e,e'\in\glr^{0,1}=\glu$:
$$
da^{(i)}_{\mu\lambda }(e,e')=-a^{(i)}_{\mu\lambda
}([e,e'])=\psi_{\mu\lambda }(a)(e,e') =\sum_{p=1}^{m_{\mu\lambda
}}c^{(i,p)}_{\mu\lambda }b^{(p)}_{\mu\lambda }.
$$
Using the isomorphism \eqref{eq:omega-A-B} once more implies the
second equation in \eqref{eq:eta-coefficients}. 
\qed

From the point of view of $\dbar$-operators, the occurrence of the
coefficients $c_{\mu\lambda }^{(i,p)}$ and 
$c^{(j,i,p)}_{\mu\nu\lambda }$ in the relations of the quiver will 
appear, via equation \eqref{eq:eta-coefficients}, when demanding the
integrability condition $(\dbar_F)^2=0$ on the $\dbar$-operators
$\dbar_F$ defined on equivariant bundles over $X\times K/J$. 
Thus, the linear terms in the relations $r^{(p)}_{\mu\lambda}$,
corresponding to the coefficients $c_{\mu\lambda }^{(i,p)}$, can
be seen as a consequence of a non-holomorphic phenomenon
(i.e. $\dbar_{A'_{\mu\lambda}}(\eta_{\mu\lambda}^{(i)})$ may not be
zero). Note that the forms $\eta_a$ are holomorphic precisely when
the unipotent radical $U$ is abelian, e.g. for Grassmann varieties. It
is worth remarking that Hille already proved that the relations are
quadratic when the flag variety $G/P$ is a Grassmannian (cf. \cite[Corollary
2.2]{Hl2}), using a different method which involves his level function
(cf. \S \ref{sec:flags}), and which, therefore, in principle can only be
applied when $P$ is a parabolic subgroup.\\

\subsection{Proof of Proposition \ref{prop:corresp-inv-hol-str}}
(a) Let us fix a $\dbar$-operator $\dbar_{E_\lambda^\circ}$ on each
$E_\lambda$. We define $\dbar$-operators $\dbar_{F^\circ_\lambda}$ 
on each $F_\lambda$ by $\dbar_{F^\circ_\lambda }=
p^*\dbar_{E^\circ_\lambda }\otimes\id + \id\otimes q^*\dbar_{H_\lambda
}$. They are obviously $G$-invariant, so the $\dbar$-operator
$\dbar_{F^\circ}=\sum_{\lambda\in Q_0'}
\dbar_{F^\circ_\lambda }\circ\pi_\lambda$ on $F$, is 
also $G$-invariant. Thus, any
$G$-invariant $\dbar$-operator $\dbar_F$ on $F$ can be written as
$\dbar_F=\dbar_{F^\circ}+\theta$ for $\theta\in\Omega^{0,1}(\End(F))^G$, 
the $G$-invariant subset of $\Omega^{0,1}(\End(F)$.  What we
have to prove is that there is a one-to-one correspondence between
$\Omega^1(\End(F))^G$ and $\Theta\times\RRR(Q',\bmE)$, where
$\Theta:=\bigoplus_{\lambda\in Q_0'}\Omega^{0,1}(\End(E_\lambda))$.
It is clear that
$
\Omega^{0,1}(\End(F))^G\cong\bigoplus_{\lambda,\mu\in Q_0'}
\Omega^{0,1}(\Hom(F_\lambda,F_\mu))^G, 
$ 
where $T^{0,1~*}_\CC (X\times K/J)\cong p^*T^*_\CC X\oplus q^*T^{0,1~*}_\CC (K/J)$, so
\begin{align*}
\Omega^{0,1}(\Hom(F_\lambda,F_\mu))^G\cong &\Omega^{0,1}(\Hom(E_\lambda,E_\mu))
\otimes\Omega^0(\Hom(H_\lambda ,H_\mu))^G \\
& \quad\oplus \Omega^0(\Hom(E_\lambda ,E_\mu))\otimes\Omega^{0,1}(\Hom(H_\lambda ,H_\mu))^G. 
\end{align*}
By the isomorphisms \eqref{eq:C-ij(V)} and \eqref{eq:omega-A-B} and Schur's lemma, it follows that 
$$
\Omega^{0,1}(\End(F))^G\cong\bigoplus_{\lambda\in Q_0'}
\Omega^{0,1}(\End(E_\lambda)) \oplus\bigoplus_{\lambda,\mu\in
Q_0'} A_{\mu\lambda }\otimes\Omega^0(\Hom(E_\lambda ,E_\mu))
\cong\Theta\times\RRR(Q',\bmE).
$$
(b) Let $\phi_{\mu\lambda}^{(i)}=\phi_a$, $\beta_{\mu\lambda
}^{(i)}=\beta_a$ for each $a=a_{\mu\lambda }^{(i)}\in Q_1'$, let
$\dbar_{E_{\mu\lambda}}$ be the $\dbar$-operator induced by
$\dbar_{E_\lambda}$ and $\dbar_{E_\mu}$ on the vector bundle
$\Hom(E_\lambda,E_\mu)$, and let $\dbar_F\in\DDD^G$ be
given by (\ref{eq:corresp-hol-str}). Then
$$
(\dbar_F)^2     =\sum_{\lambda\in Q_0'}(\dbar_{F_\lambda
})^2\circ\pi_\lambda  +\sum_{\lambda ,\mu\in
Q_0'}\sum_{i=1}^{n_{\mu\lambda}}\dbar_{E_{\mu\lambda}}
(\beta^{(i)}_{\mu\lambda })\circ\pi_\lambda 
 +\sum_{\lambda ,\nu ,\mu\in Q_0'}\sum_{i=1}^{n_{\nu\lambda
}}\sum_{j=1}^{n_{\mu\nu }} (\beta^{(j)}_{\mu\nu
}\wedge\beta^{(i)}_{\nu\lambda })\circ\pi_\lambda ,
$$
where $(\dbar_{F_\lambda})^2=p^*(\dbar_{E_\lambda })^2\otimes\id
+\id\otimes q^*(\dbar_{H_\lambda })^2=p^*(\dbar_{A_\lambda
})^2\otimes\id$. Therefore (\ref{eq:eta-coefficients}) implies
\begin{multline*}
(\dbar_F)^2=\sum_{\lambda\in Q_0'}(p^*(\dbar_{A_\lambda
})^2\otimes\id)\circ\pi_\lambda  +\sum_{\lambda ,\mu\in
Q_0'}\sum_{i=1}^{n_{\mu\lambda }} (p^*\dbar_{E_{\lambda_\mu}}
(\phi^{(i)}_{\mu\lambda })\otimes q^*\eta^{(i)}_{\mu\lambda
})\circ\pi_\lambda \\ +\frac{1}{2}
\sum_{\lambda ,\mu\in Q_0'}\sum_{p=1}^{m_{\mu\lambda }}
p^*\left(\sum_{\nu\in Q_0'}\sum_{i=1}^{n_{\nu\lambda
}}\sum_{j=1}^{n_{\mu\nu }} c^{(j,i,p)}_{\mu\nu\lambda
}\phi^{(j)}_{\mu\nu }\circ\phi^{(i)}_{\nu\lambda
}+\sum_{i=1}^{n_{\mu\lambda }}c^{(j,i,p)}_{\mu\lambda
}\phi^{(i)}_{\mu\lambda }\right)\otimes q^*\xi^{(p)}_{\mu\lambda }.
\end{multline*}
So $(\dbar_F)^2=0$ if and only if  $(\dbar_{E_\lambda })^2=0$ for all $\lambda$,
$\dbar_{E_{\mu\lambda}}(\phi^{(i)}_{\mu\lambda })=0$ for all
$\lambda,\mu,i$, and the corresponding holomorphic $Q$-bundle
$\bcR=(\bcE ,\bphi)$ satisfies the relations in $\cK$.\\ 
(c) By Schur's lemma, $(\GGG^{c})^G$ is included in 
$$
\Omega^0(\End(F)))^G\cong\bigoplus_{\lambda,\mu\in
Q_0'}\Omega^0 (\Hom(F_\lambda ,F_\mu))^G 
\cong\bigoplus_{\lambda\in Q_0'}\Omega^0(\End(F_\lambda))^G\cong
\bigoplus_{\lambda\in Q_0'}\Omega^0(\End(E_\lambda)).
$$ 
The result is now immediate.\\ 
(d) This is trivial.
\qed  

\section{Gauge equations, stability, and dimensional reduction}
\label{sec:dim-red}

The goal of this section  is to study natural gauge equations and
stability criteria for equivariant bundles on $X\times G/P$, and to 
investigate their dimensional reduction to $X$. Proposition
\ref{prop:isotopical-equiv-fil-1} establishes an equivalence between
holomorphic equivariant bundles on $X\times G/P$ and their
$G$-equivariant filtrations. As a result, one can consider certain
{\em deformed} Hermite--Einstein equation and {\em deformed}
stability criteria on equivariant bundles, which take into account the
additional structure encoded in the filtration. As we shall see, these
deformations will depend on as many 
{\em stability parameters} $\sigma_0,\sigma_1,\ldots,\sigma_{m-1}$ 
as steps are in the filtration. This deformed notions were already
used in the case $G/P=\PP^1$ in \cite{AG1} (see also \cite{M}), where
the authors proved  a Hitchin--Kobayashi correspondence for
equivariant holomorphic filtrations. We review this correspondence in
\S \ref{subsection:HKC-fil}.

In Theorem \ref{thm:equivalence-categories} we proved an
equivalence between the category of $G$-equivariant holomorphic
vector bundles on $X\times G/P$ and the category of holomorphic
$(Q,\cK)$-bundles on $X$, where $(Q,\cK)$ is the quiver with
relations associated to $P$. In fact, in Proposition
\ref{prop:corresp-inv-hol-str} we described this equivalence in terms
of $\dbar$-operators. Thus, the deformed 
Hermite--Einstein equation and deformed stability condition lead, by
the so-called process of dimensional reduction, to new gauge equations
and stability conditions on the associated quiver bundle on $X$, while
the Hitchin--Kobayashi correspondence for equivariant holomorphic
filtrations on $X\times G/P$ lead to a Hitchin--Kobayashi
correspondence for holomorphic quiver bundles on $X$. 
By using Proposition \ref{prop:corresp-inv-hol-str}, in 
Sections \ref{sub:dim-red-eq} and \ref{sub:dim-red-stab}
we describe the resulting equations and stability criteria
for quiver bundles on $X$, which will also depend on   
certain stability parameters. We relate the stability parameters
for the equivariant filtration on $X\times G/P$ with the stability
parameters for the corresponding quiver bundle.
In these two subsections we notice that there is some
freedom in the relationship between the parameters, that can be traced
back to the fact that the $K$-invariant symplectic form
$\omega_{\bvarepsilon}$ depends on as many positive parameters
$\varepsilon_\alpha$ as simple non-parabolic roots $\alpha\in\Sigma$ of $G$ which
are not roots of the Levi subgroup $L$ of $P$. Another effect of the
process of dimensional reduction is the appearance of purely
group-theoretic multiplicity factors $n_\lambda=\dim M_\lambda$, for 
$\lambda\in\Lambda^+_P$ integral dominant weights of $P$, in the 
relation between the parameters. 

\subsection{Hitchin--Kobayashi correspondence for equivariant
holomorphic filtrations} 
\label{subsection:HKC-fil}

Given a holomorphic vector bundle $\cF$ on a \kah\ manifold $(M,\omega)$, we define
the {\em degree} and {\em slope} of $\cF$ as 
\begin{equation}\label{eq:conventions-deg-slope}
\deg(\cF)=\frac{1}{\Vol(M)}\int_M\tr(\imag\Lambda F_A)\frac{\omega^{n}}{n!},
\quad \mu(\cF)=\frac{\deg(\cF)}{\rk(\cF)},
\end{equation}
respectively, where $n=\dim_{\CC}(M)$, $\Vol(M)$ is the volume of $M$,
$F_A$ is the curvature of a connection $A$ on $E$, $\Lambda$ is
contraction with $\omega$, and $\rk(\cF)$ is the rank of $\cF$. 
The degree of a torsion-free coherent sheaf  on $M$ is normalised
with these conventions as well. The {\em Hermite--Einstein equation}
for a hermitian metric $h$ on a holomorphic vector bundle $\cF$ is
$\imag\Lambda F_h=\mu(\cF) I$, where $F_h$ is the curvature of the
Chern connection associated to the metric $h$ on $\cF$.

\bfsubsubsection{Deformed Hermite--Einstein equation}
\label{subsub:THE}
Throughout this subsection, $(M,\omega)$ is a compact \kahler\
$K$-manifold, where $K$ is acompact Lie group. We assume then that the $K$-action leaves $\omega$ invariant, and 
extends to a unique holomorphic action of $G$ (the complexification of $K$) on $M$. 
Let $\bcF$, given by 
\begin{equation}\label{eq:hfil-2}
\bcF:\,\holfil 
\end{equation}
be a $G$-equivariant holomorphic filtration
over $M$. The {\em deformed Hermite--Einstein equation} involves as
many parameters $\tau_0,\tau_1,\ldots ,\tau_m\in\RR$ as steps are in
the filtration, and has the form  
\begin{equation}\label{eq:tau-metric-HEE}
\imag\Lambda F_h=\begin{pmatrix}
                \tau_0 I_0      &               &        &              \\ 
                                & \tau_1 I_1    &        &              \\ 
                                &               & \ddots &              \\ 
                                &               &        & \tau_m I_m
\end{pmatrix},
\end{equation}
where the RHS is a diagonal matrix, written in blocks corresponding
to the splitting which a hermitian metric $h$ defines in the
filtration $\bcF$. If $\tau_0=\cdots=\tau_m$,
\eqref{eq:tau-metric-HEE} reduces to the Hermite--Einstein equation. 
Taking traces in \eqref{eq:tau-metric-HEE} and integrating
over $M$, we see that there are only $m$ independent parameters, since
they are constrained by 
\begin{equation}\label{eq:constraint-tau}
\sum_{i=0}^m\tau_i\rk(\cF_i/\cF_{i-1})=\deg(\cF).
\end{equation}
The group $K$ acts in the space $\Met$ of hermitian metrics on $\cF$,
in a natural way, by $K\times\Met\to\Met$, $(\gamma,h)\mapsto
\gamma\cdot h=(\gamma^{-1})^* h$. We are interested in the
$K$-invariant solutions of \eqref{eq:tau-metric-HEE}.

\begin{definition}\label{def:THE}
Let $\btau=(\tau_0,\ldots,\tau_m)\in\RR^{m+1}$, and let $\bcF$, as in
\eqref{eq:hfil-2}, be a holomorphic filtration on $M$. We say that a
hermitian metric $h$ on $\cF$ is a $K$-invariant solution of the
$\btau$-{\em Hermite--Einstein equation} on $\bcF$ if it is
$K$-invariant and it satisfies the $\btau$-Hermite--Einstein equation
\eqref{eq:tau-metric-HEE}. If such an $K$-invariant
$\btau$-Hermite--Einstein metric $h$ exists on $\bcF$, we say that
$\bcF$ is a $K$-invariantly $\btau$-{\em Hermite--Einstein
holomorphic filtration}.
\end{definition}

\bfsubsubsection{Deformed stability}

As in the ordinary Hermite--Einstein equation, the existence of
invariant solutions to the $\btau$-Hermite--Einstein equation on an
equivariant holomorphic filtration is related to a stability condition
for the equivariant holomorphic filtration.

\begin{definition}\label{def:stab-sheaf-fil}
Let $\bsigma=(\sigma_0,\ldots,\sigma_{m-1})\in\RR^{m-1}$, and let 
$\bcF$, as in (\ref{eq:hfil-2}), be a $G$-equivariant sheaf filtration
on $M$. We define its $\bsigma$-{\em degree} and $\bsigma$-{\em slope}
respectively by 
$$
\deg_\bsigma(\bcF)=\deg(\cF)+\sum_{i=0}^{m-1} \sigma_i\rk(\cF_i),
\quad
\mu_\bsigma (\bcF)=\frac{\deg_\bsigma(\bcF)}{\rk(\cF)}.
$$
We say that $\bcF$ is $G$-invariantly $\bsigma$-{\em (semi)stable} if
for all $G$-invariant proper sheaf subfiltrations $\bcF'\hra\bcF$, we
have $\mu_\bsigma(\bcF')<(\leq)\mu_\bsigma(\bcF)$. 
A $G$-invariantly $\bsigma$-{\em polystable} sheaf filtration is a
direct sum of $G$-invariantly $\bsigma$-stable sheaf filtrations, all
of them with the same $\bsigma$-slope.
\end{definition}

\bfsubsubsection{Hitchin--Kobayashi correspondence}

\begin{theorem}\label{thm:HKC-fil}
Let $\bcF$ be a $G$-equivariant holomorphic filtration on a compact \kahler\ $K$-manifold 
$M$. Let $\btau=(\tau_0,\cdots,
\tau_m)\in\RR^{m+1}$ be related  by \eqref{eq:constraint-tau} and
let $\bsigma=(\sigma_0,\ldots,\sigma_{m-1})\in\RR^m$ be defined by
$\sigma_i=\tau_{i+1}-\tau_i$, for $0\leq i\leq m-1$, so that
$\sigma_i>0$. Then $\bcF$ admits a $K$-invariant 
$\btau$-Hermite--Einstein metric if and only if it is $G$-invariantly  $\bsigma$-polystable.
\end{theorem}

\proof 
This theorem was proved in \S 2.3 of \cite{AG1} when $G$ is trivial
and also for $SL(2,\CC)$-equivariant
holomorphic filtrations on $X\times\PP^1$. 
Since Proposition \ref{prop:isotopical-equiv-fil-1}
establishes that the $G$-equivariant coherent sheaves admit
similar $G$-equivariant filtrations, the same proof given there can be
applied to this more general situation.
\qed

\subsection{Dimensional reduction and equations}
\label{sub:dim-red-eq}

We study now  the dimensional reduction of the deformed 
Hermite--Einstein equation for a $G$-equivariant holomorphic
filtration on $X\times G/P$, as defined in \S \ref{subsub:THE}.
This subsection is organised as follows.
In \S \ref{sub:inv-kahler}  we study $K$-invariant
\kah\ form on $G/P$, since the 
\kahler\  metric on $X\times  G/P$ plays an important role
in the Hermite--Einstein equation for the equivariant holomorphic
filtration.
In \S \ref{subsub:dim-red-eq}
we define natural  gauge equations for quiver bundles,
which actually make sense for any 
quiver, not necessarily associated to a Lie group $P$
\cite{AG2}.

The main result
is then stated in Theorem \ref{thm:corresp-equations}, which is proved in
\S \ref{subsub:proof-corresp-equations} after some preliminaries about
connections on homogeneous bundles, which are covered in \S
\ref{subsub:connection-one-form-c}.

\bfsubsubsection{Invariant \kah\ structures on $K/J$}
\label{sub:inv-kahler}

The following lemma, which we adapt to the notation used throughout
this paper, is standard and can be found e.g. in \cite{Be}. Let
$\kappa(\cdot ,\cdot)$ be the Killing form on $\glg$, given by 
$\kappa(e,e')=\tr(\ad(e)\circ\ad(e'))$, $e,e'\in\glg$.
Let $\cS$ be a system of simple
roots of $\glg$ with respect to $\glh$, such that all the 
negative roots of $\glg$ are roots of $\glp$,
and let $\Sigma$ be the set of non-parabolic simple roots, defined as
in \S \ref{subsub:notation-S-Sigma}.  
Further, let $\glt$ be as in \S \ref{subsub:normal-cx-str}.
For each $\alpha\in\cS$, let $h_\alpha$ be dual to the co-root
$\alpha^\vee=2\alpha/\kappa(\alpha,\alpha)$ w.r.t. the Killing form
(so $\{ h_\alpha|\alpha\in\cS\}$ is a basis of $\glh$). Thus, $\{\imag
h_\alpha|\alpha\in\cS\}$ is a basis of $\glt$, and $-\kappa(\cdot
,\cdot)$ restricted to $\glt$ is an inner product. 

\begin{lemma} \label{lemma:inv-kah-forms}
There is a bijection between the set of $K$-invariant \kah\ forms on
$K/J$ compatible with the canonical complex structure and the set
$\RR^\Sigma_+$ of collections of positive numbers indexed by
$\Sigma$. The bijection associates to a collection $\bvarepsilon$ of
positive numbers $\varepsilon_\alpha$, for $\alpha\in\Sigma$, the
unique $K$-invariant \kah\ form $\omega_{\bvarepsilon}$ whose value
at the base point $o=J\in K/J$ of $K/J$ is the $J$-invariant 2-form
$\omega_{o \bvarepsilon}\in\wedge^2\glr^*$ given by
\begin{equation}\label{eq:inv-kahler}
\omega_{o\bvarepsilon }(x,x')=-\kappa(t,[x,x'])
\end{equation}
for $x,x'\in\glr$. Here $t\in\glt$ is given in terms of
$\bvarepsilon$ by $-\kappa(t,\imag h_{\alpha })=\varepsilon_{\alpha }$ 
for $\alpha\in\Sigma$, and $-\kappa(t,\imag h_{\alpha })=0$ for
$\alpha\in\cS\setminus\Sigma$.
\end{lemma}

\proof It is well known (see e.g. \cite[Proposition 8.83]{Be}) that
the set of $K$-invariant \kah\ forms on $K/J$ compatible with the
canonical complex structure is in bijection with the set of vectors
$t\in\glt$ such that $-\kappa(t,\imag h_\alpha)>0$ for the vectors $\imag
h_\alpha$ spanning the center $\glz=\oplus_{\alpha\in\Sigma} \RR\imag
h_\alpha$ of $\glj$, and $\kappa(t,\imag h_\alpha)=0$ for the vectors $\imag
h_\alpha$ in the semisimple part $[\glj,\glj]$ of $\glj$.
In other words, the set of $K$-invariant \kah\ forms on $K/J$,
compatible with the canonical complex structure, are in one-to-one
correspondence with the set of vectors $\bvarepsilon\in
\RR^\Sigma_+$. Formula \eqref{eq:inv-kahler} follows for instance from
\cite[Proposition 8.83]{Be}; actually, the stabiliser of the vector 
$t$ for the adjoint action of $K$ on $\glt$ is $J$, so the adjoint
orbit $K\cdot t$ is isomorphic to $K/J$, and by identifying $\glt$
with $\glt^*$ by means of the Killing form, $\omega_{o\bvarepsilon }$,
as given by \eqref{eq:inv-kahler}, transforms into be the well-known 
Kirillov--Kostant symplectic form on the coadjoint orbits. \qed

\bfsubsubsection{Statement of the main result}\label{subsub:dim-red-eq}
Let $\cS$
and $\Sigma$ be defined 
as in \ref{subsub:notation-S-Sigma}, and let  $\Delta_+(\glr)$ be the
set of positive roots of $\glg$ with respect to $\glh$ and $\cS$ which
are not roots of $\gll$ (cf. \S \ref{subsub:normal-cx-str}). 
Let $\{\lambda_\alpha |\alpha\in\cS\}$ be the basis of fundamental
weights of $\glh^*$, i.e. it is the dual basis of $\{ h_\alpha
|\alpha\in\cS\}$ given in \S \ref{sub:inv-kahler}. 
Let $\bvarepsilon$ be a collection of
positive real numbers $\varepsilon_\alpha$, for $\alpha\in\Sigma$.

Let $\btau$ and  $\btau'$  be a collection of real parameters
$\tau_\lambda$ and  $\tau_\lambda'$ 
for each $\lambda\in Q_0$, related to 
$\bvarepsilon$  by 
\begin{equation}\label{eq:dim-red-tau'-tau}
\tau'_\lambda=n_\lambda\tau_\lambda - n_\lambda 
\sum_{\alpha\in\Delta_+(\glr)}\varepsilon^{-1}_\alpha\kappa(\lambda ,\alpha^\vee), \quad {\rm ~for~}\lambda\in Q_0,
\end{equation}
where $n_\lambda=\dim_\CC M_\lambda$, and 
$\varepsilon_\alpha$ is defined, for $\alpha\not\in\Sigma$, by
\begin{equation}\label{eq:thm-def-sigma}
\varepsilon_{\alpha }:=\sum_{\beta\in\Sigma }\varepsilon_\beta \kappa(\lambda_\beta ,\alpha^\vee).
\end{equation}

\begin{remark}{\rm
The  numbers
$\varepsilon_\alpha$ in \eqref{eq:thm-def-sigma} do not depend on
the choice of $\kappa$, i.e. any $\gll$-invariant metric $\kappa$ on
$\gll$ gives the same $\varepsilon_\alpha$. Actually, if we multiply
the Killing form $\kappa$ by a positive constant $c>0$, then we obtain
the same $\varepsilon_\alpha$ in \eqref{eq:thm-def-sigma} (for
$\kappa\mapsto c \kappa$ gives $\alpha^\vee\mapsto
c^{-1}\alpha^\vee$), and clearly this transformation can be made 
separately for the factors of the Killing form corresponding to the
different simple factors of $\gll$.}
\end{remark}

\begin{theorem}\label{thm:corresp-equations}
Let $\cF$ be a $G$-equivariant holomorphic vector bundle on $X\times
G/P$. Let  $\bcF$ be the $G$-equivariant holomorphic filtration associated to  $\cF$
 and
$\bcR=(\bcE ,\bphi)$ be its corresponding  holomorphic $(Q,\cK)$-bundle on
$X$, where $(Q,\cK)$ is  the quiver with relations associated to $P$.
Then $\bcF$ has a $K$-invariant $\btau$-Hermite--Einstein
metric, with respect to the \kah\ form $p^*\omega +q^*\omega_\bvarepsilon$, 
if and only if the vector bundles $\cE_\lambda$ in $\bcR$
admit hermitian metrics $k_\lambda$ on $\cE_\lambda$,
for each $\lambda\in Q_0$ with $\cE_\lambda\neq 0$, satisfying
\begin{equation}\label{eq:vortex-connection}
\imag n_\lambda \Lambda F_{k_\lambda}
+\sum_{a\in h^{-1}(\lambda)}\phi_a\circ\phi_a^*-\sum_{a\in t^{-1}(\lambda)}\phi_a^*\circ\phi_a 
=\tau_\lambda'\id_{\cE_\lambda}, \quad\quad 
\end{equation}
 where $F_{k_\lambda}$ is the
curvature of the Chern connection $A_{k_\lambda}$ associated to the metric
$k_\lambda$ on the holomorphic vector bundle $\cE_\lambda$, for each
$\lambda\in Q_0$ with $\cE_\lambda\neq 0$, and
$n_\lambda=\dim_\CC(M_\lambda)$ is the multiplicity of the irreducible
representation $M_\lambda$, for each $\lambda\in Q_0$.
\end{theorem}

To prove this theorem, we need some preliminaries about 
connections on irreducible homogeneous vector bundles. 

\bfsubsubsection{Hermite--Einstein connections on irreducible
homogeneous vector bundles}
\label{subsub:connection-one-form-c}
In this subsection we evaluate the slope
$\mu_{\bvarepsilon}(\cO_\lambda)$ of any irreducible homogeneous
vector bundle, with respect to the invariant \kah\ form $\omega_{\bvarepsilon}$
defined in Lemma \ref{lemma:inv-kah-forms}. To do this, we reprove a
well-known fact that these bundles are Hermite--Einstein (a result 
originally due to Kobayashi \cite{Ko}), hence stable (as originally
proved by Ramanan \cite{Ra}; see also Umemura \cite{U}). 
Let $A'_\lambda$ be the $K$-invariant connection induced by the
connection one-form $\wt{A}'$ on $H_\lambda$
(cf. 
\S \ref{subsub:notation-H-lambda}, \eqref{eq:omega-A-B}, 
\S \ref{subsub:an-invariant-connection}), 
which is unitary with respect to the unique (up to a constant)
$K$-invariant hermitian metric $k'_\lambda$ on $H_\lambda$. 
Let $\End(H_\lambda,k'_\lambda)$ be the vector bundle of
anti-hermitian endomorphisms of $(H_\lambda,k'_\lambda)$. 
Contraction with the \kah\ form $\omega_{\bvarepsilon}$
(cf. \S \ref{lemma:inv-kah-forms}) is denoted by
$\Lambda_{\bvarepsilon}$. 

\begin{lemma}\label{lemma:connection-one-form-c}
The connection $A'_\lambda$ is the unique $K$-invariant connection 
on $H_\lambda$. It defines the unique $G$-invariant holomorphic 
structure $\dbar_{H_\lambda}$ on $H_\lambda$, and it is unitary with 
respect to the unique (up to scale) $K$-invariant hermitian metric
$k'_\lambda$ on $H_\lambda$. Moreover,
$A'_\lambda$ is Hermite--Einstein with respect to the \kah\ form
$\omega_\bvarepsilon$ (i.e. $\imag\Lambda_\bvarepsilon F_{A'_\lambda
}=\mu_{\bvarepsilon}(\cO_\lambda)\id)$, and the slope of $\cO_\lambda$
with respect to the \kah\ form $\omega_{\bvarepsilon}$ is
\begin{equation}\label{eq:HE-homogeneous}
\mu_{\bvarepsilon}(\cO_\lambda)
=\sum_{\alpha\in\Delta_+(\glr)}\varepsilon^{-1}_\alpha\kappa(\lambda
,\alpha^\vee).
\end{equation}
Here, $\varepsilon_\alpha$ is defined by \eqref{eq:thm-def-sigma} for
$\alpha\not\in\Sigma$. 
\end{lemma}

\proof By construction, $A'_\lambda$ is $K$-invariant and unitary. 
Its curvature $F_{A'_\lambda}$\-$\in
\Omega^{1,1}(\End(H_\lambda,k'_\lambda))$ is given, in terms of the
isomorphic element $F'_\lambda\in C^{1,1}(\End
(M_\lambda),k'_{\lambda})$, by $F'_\lambda(x,x')=-\rho_\lambda
(A'([x,x']))$, where $\rho_\lambda:K\to U(M_\lambda)$ is the unitary
representation associated to the dominant integral weight 
$\lambda$. Thus, $A'_\lambda$ defines a holomorphic 
structure on $H_\lambda$. It is easy to see that $A'_\lambda$ is the
only $K$-invariant connection on $H_\lambda$, because any other 
would be $d_{A'_\lambda}+\theta$ with $\theta\in\Omega^1_N
(\End(H_\lambda))^K\cong(\glr^*_\CC\otimes\End(M_\lambda))^J
=A_{\lambda\lambda }+A^*_{\lambda\lambda }=0$ (cf. \S  
\ref{eq:trivial-but-useful} and $\glr_{\CC}=\bar{\glu}\oplus\glu$). 
Analogously one proves that $\dbar_{A'_\lambda}$ is the only
$G$-invariant $\dbar$-operator on $H_\lambda$. Since
$\omega_{\bvarepsilon }$ and $F_{A'_\lambda}$ are both
$K$-equivariant, $\imag\Lambda_\bvarepsilon F_{A'_\lambda
}$ is $K$-equivariant as well, so by Schur's lemma,
$\imag\Lambda_\bvarepsilon F_{A'_\lambda
}=\mu_{\bvarepsilon}(\cO_\lambda)\id$ for some constant 
$\mu_{\bvarepsilon}(\cO_\lambda)$, which of course is the slope
$\mu_{\bvarepsilon}(\cO_\lambda)$ of $H_\lambda$
w.r.t. $\omega_{\bvarepsilon}$.  
To evaluate $\mu_{\bvarepsilon}(\cO_\lambda)$, first 
we compute (the complexification of) $\omega_{o\bvarepsilon}
\in\wedge^2\glr_{\CC}$. For every root $\alpha$ of $\glg$, let
$e_\alpha$ be the corresponding Chevalley generator, and for each 
pair of roots $\alpha ,\beta$, let $N_{\alpha \beta }\in\ZZ$ be 
the coefficients defined by the commutation relations
$[e_\alpha ,e_\beta ]=N_{\alpha\beta }e_{\alpha +\beta }$ 
if $\alpha+\beta$ is a root as well, and $N_{\alpha\beta }=0$
otherwise. Note that $\bar{e_\alpha}=-e_{-\alpha}$ for $\alpha
\in \Delta_+(\glr)$. Let $\alpha,\beta$ be two roots of
$\glr_{\CC}=\glu\oplus\bar{\glu}$; if $\alpha\neq\beta$, then 
$\omega_{o\bvarepsilon}(e_\alpha,\bar{e_\beta})=-\kappa(t,
-N_{\alpha,-\beta} e_{\alpha-\beta})=0$ since the center $\glz$ of
$\glj$ is orthogonal to $e_{\alpha-\beta}$, while if $\alpha=\beta$,
then $[e_\alpha,e_{-\alpha}]= h_\alpha$ implies $\omega_{o\bvarepsilon}
(e_\alpha,\bar{e_\alpha})=-\kappa(t,h_\alpha)$. To evaluate this
number, first we expand $\alpha^\vee=\sum_{\beta\in\cS}\kappa(\alpha^\vee,
\lambda_\beta)\beta^\vee$ and take into account that $-\imag\beta^\vee(t)
=-\kappa(t,\imag h_\beta)$ for any root $\beta$, and that this is
$\varepsilon_\beta$ if $\beta\in\Sigma$ and zero if
$\beta\in\cS\setminus\Sigma$. Thus,  
$$
\omega_{o\bvarepsilon} (e_\alpha,\bar{e_\alpha})
=-\imag\alpha^\vee(t) 
=-\imag\sum_{\beta\in\cS}\kappa(\alpha^\vee,\lambda_\beta)\beta^\vee(t)
=-\imag\sum_{\beta\in\Sigma}\kappa(\alpha^\vee,\lambda_\beta)\varepsilon_\beta
=\varepsilon_\alpha ,
$$
where $\varepsilon_\alpha$ is as given in \eqref{eq:thm-def-sigma} for
$\alpha\notin\Sigma$. Therefore 
\begin{equation}\label{eq:omega-varepsilon}
\omega_{o\bvarepsilon }=\imag\sum_{\alpha\in\Delta_+(\glr)}g_{\alpha\overline{\beta }
}e^{\alpha }\wedge \overline{e^{\beta }} \quad
{\rm with~} g_{\alpha\overline{\beta }}:=\delta_{\alpha\overline{\beta }}
\varepsilon_\alpha.
\end{equation}
If $\alpha\in\Delta_+(\glr)$ then $[e_\alpha ,\overline{e_\alpha }]=
-[e_\alpha ,e_{-\alpha }]=-h_\alpha\in\glh\subset\gll$ so
$F'(e_\alpha ,\overline{e_{\alpha
}})=A'(h_\alpha)=h_\alpha$. Therefore, contraction of
$\omega_{o\bvarepsilon }$ with $F'\in C^{1,1}(\glj)$ is
$$
\imag\Lambda_{o\bvarepsilon }F'
=\sum_{\alpha\in\Delta_+(\glr)}g^{\alpha\overline{\beta}}
F'(e_\alpha ,\overline{e_\beta })
=\sum_{\alpha\in\Delta_+(\glr)}\varepsilon^{-1}_\alpha h_\alpha, 
$$
so
$
\imag\Lambda_{o\bvarepsilon }F'_\lambda=\rho_\lambda(\imag\Lambda_{o\bvarepsilon }F')
=\sum_{\alpha\in\Delta_+(\glr)}\varepsilon^{-1}_\alpha\rho_\lambda(h_\alpha)$.
Since $\imag\Lambda_{o\bvarepsilon }F'_\lambda$ is $J$-invariant and
$M_\lambda$ is irreducible, $\imag\Lambda_{0\bvarepsilon
}F'_\lambda=c_\lambda\id$ for some number $c_\lambda$ which can be
computed by evaluating $\rho_\lambda(h_\alpha)$ at a highest weight
vector $v^+$ of $M_\lambda$:
$\rho_\lambda(h_\alpha)v^+=\lambda(h_\alpha)v^+$, where
$\lambda(h_\alpha)=\kappa(\lambda ,\alpha^{\vee })$. Thus, 
$\imag\Lambda_{o\bvarepsilon }F'_\lambda
=\left(\sum_{\alpha\in\Delta_+(\glr)}\varepsilon^{-1}_\alpha\kappa(\lambda
,\alpha^\vee)\right)\id_{M_\lambda }$ which implies
\eqref{eq:HE-homogeneous}.
\qed

\bfsubsubsection{Proof of Theorem \ref{thm:corresp-equations}}
\label{subsub:proof-corresp-equations}
The smooth vector bundle $F$ underlying $\cF$ has a $K$-equivariant
decomposition \eqref{eq:equiv-decomp}. Let $Q'$ be the subquiver of
$Q$ defined as in \S \ref{sub:dbar-cxgauge-grp-quivers}.
Proposition \ref{prop:corresp-inv-hol-str} determines the $\dbar$-operator
$\dbar_F$ associated to the holomorphic structure on $\cF$ in terms of
the $\dbar$-operators $\dbar_{E_\lambda}$ associated to 
the holomorphic structures on $\cE_\lambda$, for $\lambda\in Q'_0$,
and the maps $\phi^{(i)}_{\mu\lambda}$. Any $K$-invariant hermitian
metric on $F$ has a $K$-invariant orthogonal decomposition 
$k=\oplus_{\lambda } \wt{k}_\lambda$, where the sum is over
$\lambda\in Q'_0$, $\wt{k}_\lambda:=p^*k_\lambda\otimes q^*k'_\lambda$
is a $K$-invariant metric on $F_\lambda$, $k_\lambda$ is a metric
on $E_\lambda$, and $k_{\lambda }'$ is defined as in \S 
\ref{subsub:connection-one-form-c}. Let $A$ (resp. $A_\lambda$) be
the Chern connection associated to such a $K$-invariant hermitian
metric $k$ (resp. $k_\lambda$) on the holomorphic vector bundle $\cF$
(resp. $\cE_\lambda$). From \eqref{eq:corresp-hol-str},  
$$
\dbar_A=\sum_{\lambda\in Q_0'} \dbar_{\wt{A}_\lambda }\circ\pi_\lambda 
+\sum_{\lambda ,\mu\in Q_0'}\sum_{i=1}^{n_{\mu\lambda
}}\beta_{\mu\lambda }^{(i)}\circ\pi_\lambda ,\quad
\partial_A=\sum_{\lambda\in Q_0'} \partial_{\wt{A}_\lambda }\circ\pi_\lambda 
-\sum_{\lambda ,\mu\in Q_0'}\sum_{j=1}^{n_{\lambda\mu
}}\beta_{\lambda\mu }^{(j)*}\circ\pi_\lambda .
$$
Therefore
$F_A=F_A^{1,1}=\partial_A\circ\dbar_A+\dbar_A\circ\partial_A$ is given
by
\begin{equation}\label{eq:proof-corresp-equations}
\begin{split}
F_A &=\sum_{\lambda }F_{\wt{A}_\lambda}\circ\pi_\lambda
+\sum_a\partial_{\wt{A}_a}(\beta_a)\circ\pi_{ta}
+\sum_a\dbar_{\wt{A}^t_a}(\beta_a^*)\circ\pi_{ha}\\ &\quad\quad\quad
-\sum_{\lambda ,\nu ,\mu  }\sum_{i,j}(\beta^{(j)*}_{\nu\mu
}\wedge\beta^{(i)}_{\nu\lambda })\circ\pi_\lambda -\sum_{\lambda ,\nu
,\mu  }\sum_{i,j}(\beta^{(i)}_{\mu\nu }\wedge\beta^{(j)*}_{\lambda\nu
})\circ\pi_\lambda .
\end{split}\end{equation}
Here $\wt{A}_a$ (resp. $\wt{A}_a^t$) is the connection induced by
$\wt{A}_{ta}$ and $\wt{A}_{ha}$ on the bundle $\Hom(F_{ta},F_{ha})$
(resp. $\Hom(F_{ha},F_{ta})$). Let $A_a$ and
$A'_a$ (resp. $A_a^t$ and $A_a^{\prime t}$) be the connections induced
by $A_{ta}$, $A_{ha}$ on $\Hom(E_{ta},E_{ha})$ and by $A'_{ta}$,
$A'_{ha}$ on $\Hom(H_{ta},H_{ha})$ (resp. on
$\Hom(E_{ha},E_{ta})$ and $\Hom(E(H_{ha},H_{ta})$), respectively. To
evaluate $\imag\Lambda F_A$ in \eqref{eq:proof-corresp-equations}, we
need:\\ 
\begin{enumerate}
\item[(i)] $\imag\Lambda F_{\wt{A}_\lambda }=p^*\imag\Lambda F_{A_\lambda }\otimes q^*\id +p^*\id\otimes q^*\imag\Lambda_{\bvarepsilon
}F_{A'_\lambda }$;\\
\item[(ii)] $\imag\Lambda\partial_{\wt{A}_a}(\beta_a)
=\imag\Lambda(p^*\partial_{A_a}(\phi_a)\otimes q^*\eta_a
+p^*\phi_a\otimes q^*\partial_{A'_a}(\eta_a))$;\\
\item[(iii)] $\imag\Lambda\dbar_{\wt{A}^t_a}(\beta_a^*)
=\imag\Lambda(p^*\dbar_{A^t_a}(\phi_a^*)\otimes q^*\eta_a^*
+p^*\phi_a^*\otimes q^*\dbar_{A^{\prime t}_a}(\eta_a^*))$;\\
\item[(iv)] $\imag\Lambda (\beta^{(j)*}_{\nu\mu }\wedge\beta^{(i)}_{\nu\lambda })=p^*(\phi^{(j)*}_{\nu\mu }
\circ\phi^{(i)}_{\nu\lambda })\otimes q^*\imag\Lambda_\bvarepsilon (\eta^{(j)*}_{\nu\mu }\wedge\eta^{(i)}_{\nu\lambda })$;\\
\item[(v)] $\imag\Lambda (\beta^{(i)}_{\mu\nu }\wedge\beta^{(j)*}_{\lambda\nu })=p^*(\phi^{(i)}_{\mu\nu }
\circ\phi^{(j)*}_{\lambda\nu })\otimes q^*\imag\Lambda_\bvarepsilon (\eta^{(i)}_{\mu\nu }\wedge\eta^{(j)*}_{\lambda\nu })$.
\end{enumerate}
\n 
The expression (i) is given by lemma
\ref{lemma:connection-one-form-c}: $\imag\Lambda_\bvarepsilon F_{A'_\lambda }
=\mu_{\bvarepsilon}(\cO_\lambda) \id$. 
Let us prove that (ii) is zero. First,
$\imag\Lambda(p^*\partial_{A_a}(\phi_a)\otimes q^*\eta_a$ 
is zero because the direct sum $T(X\times K/J)=p^*TX\oplus q^*T(K/J)$
is orthogonal with respect to the metric associated to the \kah\ form
$p^*\omega +q^*\omega_{\bvarepsilon }$. To see that the
second term is zero, we note that
$\imag\Lambda\partial_{A'_a}(\eta_a)$ corresponds to
$\imag\Lambda_{o\bvarepsilon }\partial a$ by Lemma
\ref{lemma:computing-inv-cov-derivative}, and 
\eqref{eq:omega-varepsilon} gives
$$
\imag\Lambda_{o\bvarepsilon }\partial
a=\sum_{\alpha\in\Delta_+(\glr)}\varepsilon^{-1}_\alpha (\dbar a)(e_\alpha
,\overline{e_\alpha }).
$$ 
Let $\lambda=ta$ and $\mu=ha$, so $\partial a$ is in $B_{\mu\lambda}
=C^{0,2}(\Hom(M_\lambda ,M_\mu))$, which is a subspace of
$\wedge^2\glr^*_\CC\otimes\Hom(M_\lambda ,M_\mu)$ and this space is
injected in $\wedge^2\glk^*_\CC\otimes\Hom(M_\lambda ,M_\mu)$ by the 
monomorphism induced by the projection $\glk=\glj\oplus\glr
\to\glk/\glj\cong\glr$. Therefore, $\partial a$ is zero in $\glj$,
so $(\partial a)(e_\alpha ,\overline{e_\alpha })=-\partial a(e_\alpha
,e_{-\alpha })=a([e_\alpha ,e_{-\alpha }])=a(h_\alpha)=0$.
Thus, (ii) is zero. Similarly, (iii) is zero.
To evaluate (iv) (and (v)), we note that $\omega_\bvarepsilon$ is
$K$-invariant and  $\eta^{(j)*}_{\nu\mu }\wedge\eta^{(i)}_{\nu\lambda
}$ is $K$-invariant, so $\imag\Lambda_\bvarepsilon
(\eta^{(j)*}_{\nu\mu }\wedge\eta^{(i)}_{\nu\lambda })$ is
$K$-invariant as well. By Schur's lemma, it is zero if $\lambda\neq\mu$. 
If $\lambda=\mu$, we choose the  basis
$\{ a^{(i)}_{\nu\lambda }| i=1,\ldots ,n_{\mu\lambda }\}$ of
$A_{\lambda\nu }$ which is orthonormal with respect to
the metric induced  by the hermitian metric associated to
$\omega_{o\bvarepsilon }$ and the canonical complex structure, and 
by the $J$-invariant hermitian metrics
$k'_{o\lambda },k'_{o\nu }$ on $M_\lambda$ and $M_\nu$, i.e. they are 
normalised so that 
$\imag\Lambda_{o\bvarepsilon }\tr(a^{(j) *}_{\nu\lambda
}\wedge a^{(i)}_{\nu\lambda })=\delta^{ij}$ (so  
$\imag\Lambda_{o\bvarepsilon }\tr(a^{(i)}_{\mu\nu }\wedge a^{(j)*}_{\lambda\nu })=-
\delta^{ij}$). By Schur's lemma 
$$
\imag\Lambda_{o\bvarepsilon}(a^{(j) *}_{\nu\lambda }\wedge a^{(i)}_{\nu\lambda })
=\frac{1}{n_\lambda}\delta^{ij}\id_{M_\lambda},
\quad
\imag\Lambda_{o\bvarepsilon }(a^{(i)}_{\mu\nu }\wedge a^{(j)*}_{\lambda\nu })
=-\frac{1}{n_\lambda}\delta^{ij}\id_{M_\lambda}, 
$$
(since $n_\lambda=\dim_\CC(M_\lambda)$). Therefore (iv) and (v)
are given by
$$
\imag\Lambda_{\bvarepsilon }(\eta^{(j) *}_{\nu\mu }\wedge
\eta^{(i)}_{\nu\lambda })=\frac{1}{n_\lambda }\delta_{\lambda\mu }\delta^{ij}\id ,
\quad
\imag\Lambda_{\bvarepsilon }(\eta^{(i)}_{\mu\nu }\wedge \eta^{(j)*}_{\lambda\nu })=-\frac{1}{n_\lambda }\delta_{\lambda\mu }\delta^{ij}\id .
$$
Combining \eqref{eq:proof-corresp-equations} and the results for
(i)-(v), we can evaluate $\imag\Lambda F_A$: 
\begin{align*}
\imag\Lambda F_A=\sum_\lambda &
 \left( p^*\left( \imag F_{A_\lambda} + \mu_{\bvarepsilon}(\cO_\lambda)\id_{\cE_\lambda} - \frac{1}{n_\lambda }\sum_\mu\sum_i\phi^{(i)*}_{\mu\lambda}\circ\phi^{(i)}_{\mu\lambda}\right. \right.\\
& \left. \left. +\frac{1}{n_\lambda }\sum_\mu\sum_i \phi^{(i)}_{\mu\lambda}\circ\phi^{(i)*}_{\mu\lambda} \right)\otimes  q^*\id_{H_\lambda} \right)\circ \pi_\lambda .
\end{align*}
The relation between $\tau_\lambda$ and $\tau'_\lambda$ given in
\eqref{eq:dim-red-tau'-tau} is $\tau'_\lambda=n_\lambda(\tau_\lambda -
\mu_{\bvarepsilon}(\cO_\lambda))$, due to \eqref{eq:HE-homogeneous}. 
The theorem is now straightforward. 
\qed

\subsection{Dimensional reduction and stability}
\label{sub:dim-red-stab}

In this subsection we introduce a notion of stability for 
quiver sheaves on $X$ and prove that this is precisely the stability
criteria which appears by dimensional reduction for $G$-invariant
stability of $G$-equivariant sheaf filtrations, by the correspondences in 
\S \ref{sec:equivb-quivers}.
We also compute the relations among the stability parameters for the
equivariant filtration and the quiver sheaf, and the way that the
$K$-invariant \kah\ form on $K/J$ enters in these relations.

\bfsubsubsection{Stability for quiver bundles}
\label{subsub:Q-stability}

Let $Q$ be the  quiver associated to $P$,  and  let  $\btau$ be
collections of real numbers $\tau_\lambda$, 
for each $\lambda\in Q_0$. Let $\bcR=(\bcE,\bphi)$ be a $Q$-sheaf on
$X$. Let $n_\lambda=\dim_\CC(M_\lambda)$, for each $\lambda\in Q_0$. 

\begin{definition}\label{def:stab-Q-sheaf}
The $\btau$-{\em degree} and $\btau$-{\em slope}
of $\bcR$ are
$$
\deg_{\btau }(\bcR)=\sum_{\lambda\in Q_0}
\left(n_\lambda\deg(\cE_\lambda)-\tau_\lambda\rk(\cE_\lambda)\right) , 
\quad
\mu_{\btau }(\bcR)=\frac{\deg_{\btau }(\bcR)}{\sum_{\lambda\in Q_0}
n_\lambda\rk(\cE_\lambda)}, 
$$
respectively. The $Q$-sheaf $\bcR$ is called $\btau$-{\em (semi)stable} if
for all proper $Q$-subsheaves $\bcR'$ of $\bcR$,
 $\mu_{\btau}(\bcR')<(\leq)\mu_{\btau }(\bcR)$.
A $\btau$-{\em polystable} $Q$-sheaf is a direct sum of
$\btau$-stable $Q$-sheaves, all of them with the same $\btau$-slope. 
\end{definition}

\begin{remark}\label{rem:stab}{\rm 
If $\cK$ is the set of relations of $Q$ and $\bcR$ is a
$(Q,\cK)$-sheaf, then all its $Q$-subsheaves are $(Q,\cK)$-sheaves,
so the $\btau$-(semi)stability criteria does not depend on
$\cK$.\\  
}\end{remark}

\begin{theorem}\label{thm:corresp-stability}
Let $\cF$ be a $G$-equivariant sheaf on 
$X\times G/P$. Let $\bcF$ be its associated $G$-equivariant sheaf filtration,
and  $\bcR=(\bcE ,\bphi)$ be its corresponding   $(Q,\cK)$-sheaf on
$X$, where $(Q,\cK)$ is  the quiver with relations associated to $P$.
Let $Q_0'=\{\lambda_0,\ldots,\lambda_m\}$ be the set of vertices $\lambda
\in Q_0$ such that $\cE_\lambda\neq 0$, listed in ascending order as
$\lambda_0<\cdots<\lambda_m$.
Let $\bvarepsilon$ be a collection of positive real numbers
$\varepsilon_\alpha$, for each $\alpha\in\Sigma$, and let
$\varepsilon_\alpha$ be defined, for
$\alpha\in\Delta_+(\glr)\setminus\Sigma$, by \eqref{eq:thm-def-sigma}.
Let $\bsigma=(\sigma_0,\ldots,\sigma_{m-1})$ with $\sigma_s>0$ for
each $0\leq s\leq m-1$, and let $\tau'_\lambda$, for each $\lambda\in
Q'_0$, be given by
\begin{equation}\label{eq:dim-red-sigma-tau}
\tau'_{\lambda_s}=n_{\lambda_s} \sum_{s'=0}^{s-1} \sigma_{s'} - n_{\lambda_{s}}
\sum_{\alpha\in\Delta_+(\glr)}\varepsilon^{-1}_\alpha\kappa(\lambda_s,\alpha^\vee), 
\end{equation}
for $0\leq s\leq m-1$, where $n_\lambda=\dim_\CC (M_\lambda)$.
Then $\bcF$ is $G$-invariantly $\bsigma$-(semi)stable  with respect to
the \kah\ form $p^*\omega +q^*\omega_\bvarepsilon$ if and only if
$\bcR$ is $\btau'$-(semi)stable.
\end{theorem}

\proof 
To simplify the notation, let us denote $\sigma_s$ by
$\sigma_{\lambda}$ when $\lambda=\lambda_s\in Q'_0$. 
Using \eqref{eq:HE-homogeneous}, we can write
\eqref{eq:dim-red-sigma-tau} as
\begin{equation}\label{eq:dim-red-tau-sigma}
\frac{\tau'_\lambda}{n_\lambda}=\sum_{\mu<\lambda}\sigma_\mu-\mu_\bvarepsilon
(\cO_\lambda),
\end{equation}
for all $\lambda\in Q_0'$.
Let now $\bcF'$, given by \eqref{eq:hsfil},
be a $G$-invariant sheaf subfiltration of $\bcF$, and let
$\bcR'=(\bcE',\bphi')$ be the corresponding $(Q,\cK)$-subsheaf of
$\bcR$. Then
\begin{align*}
\sum_{s=0}^{m-1}\sigma_{\lambda_s}\rk(\cF'_s)
&=\sum_\lambda\sum_{\mu\leq\lambda}\sigma_\lambda\rk(p^*\cE'_\mu\otimes q^*\cO_\mu)
=\sum_\lambda\sum_{\mu\leq\lambda} \sigma_\lambda n_\mu \rk(\cE'_\mu) \\
&=-\sum_\lambda\sum_{\mu>\lambda} \sigma_\lambda n_\mu \rk(\cE'_\mu) 
+\sum_{\lambda}\sum_{\mu} \sigma_\lambda n_\mu \rk(\cE'_\mu) \\
&=-\sum_\lambda n_\lambda\left(\sum_{\mu<\lambda}\sigma_\mu\right)\rk(\cE'_\lambda)
+\left(\sum_\lambda\sigma_\lambda\right) \left(\sum_\lambda n_\lambda\rk(\cE'_\lambda)\right) ,
\end{align*}
since $\rk(\cO_\lambda)=n_\lambda$ (sums are only over
$\lambda,\mu$ in $Q'_0$), while
\begin{align*}
\deg(\cF')      &=
\sum_\lambda\deg(p^*\cE'_\lambda\otimes q^*\cO_\lambda)
=\sum_\lambda\left(\rk(\cO_\lambda)\deg(\cE'_\lambda)+\rk(\cE'_\lambda)\deg_\bvarepsilon(\cO_\lambda)\right)\\
&\quad\quad =\sum_\lambda \left(n_\lambda\deg(\cE'_\lambda)+n_\lambda\mu_\bvarepsilon(\cO_\lambda)\rk(\cE'_\lambda)\right)
\end{align*}
and
$$
\rk(\cF')=\sum_\lambda\rk(p^*\cE'_\lambda\otimes q^*\cO_\lambda)
=\sum_\lambda n_\lambda\rk(\cE'_\lambda) .
$$
Therefore
\begin{align*}
\deg_{\bsigma}(\bcF')   &=\sum_\lambda n_\lambda 
                \left(\deg(\cE'_\lambda)-\left(\sum_{\mu<\lambda}\sigma_\mu -\mu_\bvarepsilon(\cO_\lambda)\right)\rk(\cE'_\lambda)\right)
+\left(\sum_\lambda\sigma_\lambda\right)\left(\sum_\lambda n_\lambda\rk(\cE'_\lambda)\right)\\
&=\sum_\lambda  \left(n_\lambda\deg(\cE'_\lambda)-\tau'_\lambda\rk(\cE'_\lambda)\right)
+\left(\sum_\lambda\sigma_\lambda\right)\left(\sum_\lambda n_\lambda\rk(\cE'_\lambda)\right),
\end{align*}
so
$$
\mu_{\bsigma}(\bcF')=\frac{\deg_{\bsigma}(\bcF')}{\sum_\lambda
n_\lambda\rk(\cE'_\lambda)} 
=\mu_{\btau'}(\bcR')+\sum_\lambda\sigma_\lambda,
$$
and we are done.
\qed

\subsection{A \HKC\ for quiver bundles}
Combinining Theorems \ref{thm:HKC-fil}, \ref{thm:corresp-equations}
and \ref{thm:corresp-stability}, we obtain as a corollary a Hitchin-Kobayashi
correspondence for the quiver bundles that arise from dimensional reduction.
More precisely.

\begin{theorem}\label{thm:HKC-quivers}
Let $(Q,\cK)$ be   the quiver with relations associated to $P$.
Let $\bcR=(\bcE ,\bphi)$ be a  holomorphic $(Q,\cK)$-bundle on
$X$ corresponding to a $G$-equivariant holomorphic vector bundle on $X\times
G/P$. Let $\btau'$ be  a collection of real numbers
$\tau'_\lambda$, for each $\lambda\in Q_0$. Then the vector bundles $\cE_\lambda$ in $\bcR$
admit hermitian metrics $k_\lambda$ on $\cE_\lambda$,
for each $\lambda\in Q_0$ with $\cE_\lambda\neq 0$, satisfying
\begin{equation}\label{eq:vortex-connection-2}
\imag n_\lambda \Lambda F_{k_\lambda}
+\sum_{a\in h^{-1}(\lambda)}\phi_a\circ\phi_a^*-\sum_{a\in t^{-1}(\lambda)}\phi_a^*\circ\phi_a 
=\tau'_\lambda\id_{\cE_\lambda}, \quad\quad 
\end{equation}
if and only if  $\bcR=(\bcE ,\bphi)$ is $\btau'$-polystable.
\end{theorem}

\proof
First, by Theorem \ref{thm:corresp-equations}, the
bundles $\cE_\lambda$ have hermitian metrics satisfying
\eqref{eq:vortex-connection-2} if and only if the corresponding
$G$--equivariant holomorphic filtration $\bcF$ on $X\times G/P$ has a
$G$-invariant $\btau$-Hermite--Einstein metric, where $\btau$ is a
collection of parameters $\tau_\lambda$ related to $\btau'$ by
\eqref{eq:dim-red-tau'-tau}, or equivalently, by 
\begin{equation}\label{eq:dim-red-tau'-tau-bis}
\frac{\tau'_\lambda}{n_\lambda}=\tau_\lambda - \mu_{\bvarepsilon}(\cO_\lambda). 
\end{equation}
Similarly, by Theorem \ref{thm:corresp-stability}, the quiver bundle
$\bcR$ is $\btau'$-stable if and only if $\bcF$ is $G$-invariantly
$\bsigma$-stable, where $\bsigma=(\sigma_0,\ldots,\sigma_{m-1})$ is
related to $\btau'$ by \eqref{eq:dim-red-sigma-tau}, or equivalently,
by \eqref{eq:dim-red-tau-sigma} (we are using the notation of the
proof of Theorem \ref{thm:corresp-stability}).
By Theorem \ref{thm:HKC-fil}, we only have to check that
$\sigma_s=\tau_{s+1}-\tau_s$. But this follows from the previous two
equations \eqref{eq:dim-red-tau'-tau-bis} and
\eqref{eq:dim-red-tau-sigma}: 
$$
\sigma_s=\sum_{\mu<\lambda_{s+1}}\sigma_\mu - \sum_{\mu<\lambda_{s}}\sigma_\mu
=\left(\frac{ \tau'_{\lambda_{s+1}}}{n_{\lambda_{s+1}}} + \mu_{\bvarepsilon}(\cO_{\lambda_{s+1}}) \right)
-\left( \frac{\tau'_{\lambda_{s}}}{n_{\lambda_{s}}} + \mu_{\bvarepsilon}(\cO_{\lambda_{s}}) \right)
=\tau_{s+1} - \tau_s.
$$
\qed

\begin{remark}{\rm 
It is clear that equations \eqref{eq:vortex-connection-2} make also sense
for a $(Q,\cK)$-bundle, where  $(Q,\cK)$ is  not necessarily associated
to a parabolic group $P$.
 Moreover,  the $n_\lambda$ appearing in the equations
could actually be arbitrary positive numbers. All these generalizations,
not related {\em a priori} to dimensional reduction,  had been
carried out in \cite{AG2}.
}\end{remark}

\section{Examples of dimensional reduction}
\label{sec:examples}

In \S\S \ref{subsub:hb-prod-P1} and \ref{subsub:hb-P2} we described
explicitly the quiver when $G/P$ is $(\PP^1)^N$ or $\PP^2$. We
now apply Theorem \ref{thm:corresp-equations} to obtain, in
these examples, the dimensional reduction of the deformed
$\btau$-Hermite--Einstein equation, for an equivariant holomorphic
filtration $\bcF$, on $X\times G/P$, i.e. we evaluate the
multiplicities $n_\lambda$ and the stability parameters $\btau'$,
appearing in the quiver vortex equations, in terms of $\btau$ and the
parameters $\bvarepsilon$ parametrising the invariant \kah\ form
$\omega_{\varepsilon}$ on $G/P$. In particular, we shall recover the
vortex equations corresponding to holomorphic triples \cite{G2,BG}
and to their generalizations, the holomophic chains \cite{AG1}, which
are obtained by dimensional reduction when $G/P=\PP^1$. 

\subsection{Dimensional reduction for $\PP^2$}

We begin with $G/P=\PP^2$. The dimensional reduction of the
$\btau$-Hermite--Einstein equation for an equivariant holomorphic
filtration on $X\times\PP^2$ gives 
$$
\imag n_x \Lambda F_{k_x}
+\phi_x^{(1)}\circ\phi_x^{(1)*}+\phi_x^{(2)}\circ\phi_x^{(2)*}
-\phi_x^{(1)*}\circ\phi_x^{(1)}-\phi_x^{(2)*}\circ\phi_x^{(2)} 
=\tau'_x\id, \quad {\rm for~} x\in Q_0, 
$$
where we are using the notation in \S
\ref{subsub:hb-P2}, and $\phi_x^{(i)}:=\phi_a$, for $a=a_x^{(i)}\in
Q_1$ and $i=1,2$. The relations give $\phi_{x-3\epsilon_1}^{(2)}
\circ\phi_x^{(1)}=\phi_{x-3\epsilon_2}^{(1)}\circ\phi_x^{(2)}$ for
$x\in Q_0$. The multiplicities can be found e.g. in \cite[(15.17)]{FH}: 
$n_\lambda:=\dim_\CC(M_\lambda)=1+\lambda_1-\lambda_2$, or using the
notation of \S \ref{subsub:hb-P2}, 
\begin{equation}\label{eq:sigma-sl}
n_x:=1+\frac{1}{3}(x_1-x_2)
\end{equation} 
(they take all possible positive integer values,
since $(x_1-x_2)/3=\lambda_1-\lambda_2\geq 0$). Note that any arrow
$a^{(1)}_x:x\to y$ decreases $n_y=\sigma_x-1$, for
$(y_1,y_2)=(x_1,x_2+3)$ (resp. $a^{(2)}_x:x\to y$ increases
$n_y=\sigma_x+1$, for $(y_1,y_2)=(x_1+3,x_2)$). 
There is only one non-parabolic simple root, i.e. 
$\Sigma=\{\alpha_2\}$ (cf. \S \ref{subsub:hb-P2}), so the family of
invariant \kah\ forms $\omega_{\bvarepsilon}$ on $\PP^2$ only depends
on a positive real number $\varepsilon>0$ (with the notation of Lemma
\ref{lemma:inv-kah-forms}, this is $\varepsilon_{\alpha_2}$). 
According to Theorem \ref{thm:corresp-equations}, $\btau'$ is given by 
$\tau'_\lambda=n_\lambda (\tau_\lambda+\mu_{\varepsilon}(\cO_\lambda))$,
where $\mu_{\varepsilon} (\cO_\lambda)$ is given by Lemma
\ref{lemma:connection-one-form-c}. 
To evaluate $\mu_{\varepsilon} (\cO_\lambda)$, we take into account
that $\Delta_+(\glr)=\{-\gamma_1,-\gamma_2\}$ (cf. \S
\ref{subsub:dim-red-eq}), so \eqref{eq:thm-def-sigma} gives 
$\varepsilon_k:=\varepsilon_{\gamma_k}=\varepsilon
\kappa(\lambda_{\alpha_2},-\gamma_k^\vee)$, and
\eqref{eq:HE-homogeneous} gives 
$$
\mu_{\varepsilon} (\cO_\lambda)
=\sum_{k=1}^2 \varepsilon_k^{-1} \kappa(\lambda,-\gamma_k^\vee)
=\varepsilon^{-1} \sum_{k=1}^2 \frac{\kappa(\lambda,-\gamma_k^\vee)}{\kappa(\lambda_{\alpha_2},-\gamma_k^\vee)}
=\varepsilon^{-1} \sum_{k=1}^2 \frac{\kappa(\lambda,-\gamma_k)}{\kappa(\lambda_{\alpha_2},-\gamma_k)}.
$$
The Killing form is a multiple of $\kappa(\lambda,\lambda')=\sum_{i=1}^{3}
\lambda_i\lambda'_i -\frac{1}{3}\sum_{i,j=1}^{3} \lambda_i\lambda'_j$, for
$\lambda=\sum_{i=1}^{3}\lambda_i L_i$, $\lambda'
=\sum_{i=1}^{3}\lambda'_i L_i$ (cf. e.g. \cite[(15.2)]{FH}).
Now, $\gamma_k=L_k-L_{3}$ and $\lambda_{\alpha_2}=-L_{3}$, so
\begin{equation}\label{eq:slope-homvb-P^n}
\mu_{\varepsilon} (\cO_\lambda)=\frac{\lambda_1+\lambda_2}{\varepsilon} 
=-\frac{x_1+x_2}{\varepsilon}.
\end{equation}
Thus, the $\btau'$-parameters are related to the $\btau$-parameters
appearing in the $\btau$-Hermite--Einstein equation for equivariant
holomorphic filtrations by 
$$
\tau'_x=n_\lambda(\tau_\lambda+\mu_{\varepsilon}(\cO_\lambda))
=\left(1+\frac{1}{3}(x_1-x_2)\right)~\frac{\tau_\lambda -(x_1+x_2)}{\varepsilon}, 
\quad {\rm for~} x\in Q_0.
$$

\subsection{Dimensional reduction for $(\PP^1)^N$}

Let us now  consider $G/P=(\PP^1)^N$. We make use of
the results and notation of \S \ref{subsub:hb-prod-P1}.
The Levi subgroup $L\cong(\CC^*)^N$ is abelian, so its irreducible
representations are one-dimensional, i.e. $n_\lambda:=\dim_\CC(M_\lambda)=1$, 
for $\lambda\in Q_0$. Thus, the dimensional reduction of the
$\btau$-Hermite--Einstein equation for an equivariant holomorphic
filtration on $X\times(\PP^1)^N$ gives 
\begin{equation}\label{eq:vortex-(P^1)^N}
\imag\Lambda F_{k_\lambda} +
\sum_{i=1}^N\phi^{(i)}_{\lambda-2L_i}\circ \phi^{(i)~*}_{\lambda-2L_i}
-\sum_{i=1}^N\phi^{(i)~*}_{\lambda}\circ \phi^{(i)}_{\lambda}
=\tau'_\lambda\id_{\cE_\lambda}, \quad {\rm for~} \lambda\in Q_0.
\end{equation}
The relations lead to 
$\phi^{(j)}_{\lambda-2L_i}\circ\phi^{(i)}_\lambda=\phi^{(i)}_{\lambda-2L_j}
\circ\phi^{(j)}_\lambda$. All the positive roots 
$\alpha_i=2L_i$, $1\leq i\leq N$, are simple, so the invariant \kah\ form 
$\omega_{\bvarepsilon}$ on $(\PP^1)^N$ depends on $N$ 
parameters $\varepsilon_i>0$, $1\leq i\leq N$. Since the
$\bvarepsilon$-slope of $\cO_\lambda$ is
$\mu_\bvarepsilon(\cO_\lambda)=\sum_{i=1}^N \lambda_i/\varepsilon_i$,
for $\lambda=\sum_{i=1}^N\lambda_i L_i\in 
Q_0$, \eqref{eq:dim-red-tau'-tau} and \eqref{eq:dim-red-sigma-tau} give 
\begin{equation}\label{eq:tau-sigma-(P^1)^N}
\tau'_\lambda =\tau_\lambda - \sum_{i=1}^N\frac{\lambda_i}{\varepsilon_i},
\quad 
\tau'_{\lambda_s}=\sum_{s=0}^{s-1}\sigma_{s'} - \sum_{i=1}^N\frac{\lambda_{s,i}}{\varepsilon_i},
~~ {\rm for~}\lambda\in\ZZ^N, ~ {\rm and~} 0\leq s\leq m.
\end{equation}
Here $\lambda_0<\lambda_1<\cdots<\lambda_m$ are the weights appearing
in the filtration \eqref{eq:equi-hol-fil} with components
$\lambda_{s,i}\in\ZZ$ given by $\lambda_s=\sum_{i=1}^N\lambda_{s,i}L_i$.
The arrow $a^{(i)}_\lambda$ takes $\lambda_s$ into
$\lambda_s-2L_i$, so if $\cF$ is indecomposable, then necessarily
$\lambda_{s+1,i}-\lambda_{s,i}=2$ for some $1\leq i\leq N$. 
In particular, we can apply the preceding discussion to the 
complex projective line $G/P=\PP^1$: the quiver has two connected
components $Q^{(h)}$, for $h=0,1$, each one isomorphic to the quiver
with vertex set $\ZZ$ (the weight of $i\in\ZZ$ being $\lambda=2i+h$),
and arrows $a_i:i\to i-1$, for $i\in\ZZ$. The process of dimensional
reduction for any equivariant bundle on $X\times\PP^1$ was already
studied in \cite{AG1}. Each indecomposable holomorphic quiver bundle
on $X$ is given by a sequence of weights $\lambda_0,\lambda_1,\ldots,
\lambda_m \in Q_0^{(h)}$, with $\lambda_i-\lambda_{i-1}=2$, and by a
sequence of morphisms $\phi_i:\cE_i\to\cE_{i-1}$ among holomorphic
vector bundles $\cE_0,\ldots,\cE_m$, where $\cE_i$ corresponds to the
weight $\lambda_i$. Such a quiver bundle, or {\em holomorphic chain}
\cite{AG1}, $\bcC=(\bcE,\bphi)$, is specified by a diagram 
\begin{equation}
\label{eq:hol-chain}
\bcC :\, \cE_m\stackrel{\phi_m}{\lra }\cE_{m-1}\stackrel{\phi_{m-1}}{\lra}
\cdots\stackrel{\phi_1}{\lra}\cE_0.
\end{equation}
(A {\em sheaf} chain is analogously defined, by using coherent sheaves
instead of holomorphic bundles (see \cite{AG1} for more details).) 
By making a translation $\lambda\mapsto \lambda-\lambda_0$, we can always
assume that the weights are $\lambda_i=2i$, for $0\leq i\leq m$. The
equivariant holomorphic filtration on $X\times\PP^1$ corresponding to
$\bcC$, by Theorem \ref{thm:equivalence-categories} and
Proposition \ref{prop:isotopical-equiv-fil-1}, is given by 
\begin{equation}\begin{gathered}\label{eq:hfil-P^1}
\bcF:\, \holfil ,\\
\cF_i/\cF_{i-1}\cong p^*\cE_i\otimes q^*\cO(2i),
\quad 0\leq i\leq m.
\end{gathered}\end{equation}
Translating $\lambda\mapsto \lambda-\lambda_0$ corresponds to
twisting $\cF$ by $q^*\cO(-\lambda_0)$. 
The quiver vortex equations, or {\em chain $\btau$-vortex equations}
\cite{AG1}, for metrics $k_i$ on the bundles $\cE_i$, are
(cf. \eqref{eq:vortex-(P^1)^N}) 
\begin{align*}
\sqrt{-1}\Lambda F_{k_0}+\phi_{1}\circ\phi^*_{1} &=\tau_0\id_{\cE_0},\\ 
\sqrt{-1}\Lambda F_{k_i}+\phi_{i+1}\circ\phi^*_{i+1}-\phi^*_i\circ\phi_i
&=\tau_i\id_{\cE_i}, \qquad (1\leq i\leq m-1),\\
\sqrt{-1}\Lambda F_{k_m}-\phi^*_{m}\circ\phi_{m} &=\tau_m\id_{\cE_m}.
\end{align*}
The story for stability follows similarly: a coherent sheaf chain
$\bcC$ is ${\btau'}$-{\em (semi)stable} if for all proper sheaf
subchains $\bcC'\hra\bcC$, $\mu_{\btau'}(\bcC')<(\leq)
\mu_{\btau'}(\bcC)$, where  
$\mu_{\btau'} (\bcC)={\deg_{\btau'}(\bcC)}/{\sum_{i=0}^m \rk(\cE_i)}$,
$\deg_{\btau'}(\bcC)=\sum_{i=0}^m\deg(\cE_i)-\sum_{i=0}^m\tau'_i\rk(\cE_i)$. 
The relation among $\bsigma,\btau,\btau'$ is
given by \eqref{eq:tau-sigma-(P^1)^N} for $N=1$. Thus, if
the vertices of $\bcC$ are $\lambda_i=2i$, for $0\leq i\leq 2m$, then
$\tau'_i=\tau_i-\mu_\varepsilon(\cO_\lambda)=\tau_i-\frac{2i}{\varepsilon},~
\sigma_i=\tau'_{i+1}-\tau'_i+\frac{2}{\varepsilon}$.
These results were already obtained in \cite{AG1} when 
$\varepsilon=1$. The case $m=1$ has been deeply studied in
\cite{G1,G2,BG}; then the sheaf filtration \eqref{eq:hfil-P^1} is
just an invariant extension of equivariant holomorphic vector bundles
on $X\times\PP^1$, 
$$
\begin{CD}
        0   @>>> p^*\cE_0 @>>> \cF @>>> p^*\cE_1\otimes q^*\cO(2) @>>> 0,
\end{CD}
$$
and the holomorphic chain \eqref{eq:hol-chain} is a triple (cf. \cite{BG}),
i.e. a representation $\phi:\cE_1\to \cE_0$ of the quiver
$\bullet\to\bullet$. Assume now that $\sigma_0=0$,
i.e. $\tau_0=\tau_1$. Then the 
$\btau$-Hermite--Einstein equation \eqref{eq:tau-metric-HEE}, 
for a holomorphic filtration $\bcF$, is the usual Hermite--Einstein equation
for $\cF$, and the $\bsigma$-stability condition is the usual Mumford--Takemoto 
stability condition for a (torsion free) coherent sheaf on $X\times\PP^1$. 
So the parameter $\varepsilon$, which was encoded in the invariant
\kah\ form $p^*\omega+q^*\omega_{\varepsilon}$, allows us to get the
moduli space of $(\tau'_0,\tau'_1)$-stable triples on $X$ as the set of
fixed points of the moduli space of stable sheaves on $X\times\PP^1$, under
the induced action of $SL(2,\CC)$. The constraint in Remark
\ref{rem:stab}(c) is $\deg(\cE_0)+\deg(\cE_1)
=\tau'_0\rk(\cE_0)+\tau'_1\rk(\cE_1)$. Thus, the relation
$\tau'_0-\tau'_1={2}/{\varepsilon}$, among $\tau'_0,\tau'_1$, and
$\varepsilon$, turns out to be 
$$
\frac{2}{\varepsilon}=\frac{\rk(\cE_0)+\rk(\cE_1)\tau_0
-(\deg(\cE_0)+\deg(\cE_1))}{\rk(\cE_1)},
$$ 
which is precisely as in \cite{G2}. It provides the value of
$\varepsilon$ in terms of the stability parameter $\tau'_0$. 
Thus, the mechanism of dimensional reduction, when applied to the
well-known Hitchin--Kobayashi correspondence for the usual
Hermite--Einstein equation, proved by Donaldson, Uhlenbeck and Yau,
for equivariant bundles on $X\times\PP^1$, gives a Hitchin--Kobayashi
correspondence for the coupled vortex equations for holomorphic
triples (cf. \cite{BG}). However, when $m>1$ in \eqref{eq:hol-chain},
the parameter $\varepsilon$, encoded in the polarisation
$\omega_{\varepsilon}$, is not enough to obtain a Hitchin--Kobayashi
correspondence for the chain $\btau'$-vortex equations by dimensional
reduction of the usual Hermite--Einstein equation ---we need the
deformed $\btau$-Hermite--Einstein equation
\eqref{eq:tau-metric-HEE}. Of course, the same occurs for other flag
varieties.

\end{document}